\documentclass[preprint,10pt]{elsarticle}  %

\usepackage[margin=2.5cm]{geometry}
\usepackage{amsmath,amssymb,amsfonts,mathtools,bm}
\usepackage{graphicx}
\usepackage{booktabs}
\usepackage{enumitem}
\usepackage{hyperref}
\usepackage{physics}
\usepackage{xcolor}
\usepackage{caption}
\usepackage{subcaption}
\usepackage{amsthm}
\usepackage{float}
  \usepackage{xcolor}
\usepackage[linesnumbered,ruled]{algorithm2e}

\SetCommentSty{mycommfont}

\floatstyle{boxed}
\newfloat{BOX}{!ht}{box}[section]
\floatname{BOX}{Box}
\newtheorem{remark}{Remark}

\newcommand{\rowdos}[2]{\begin{bmatrix} #1  & #2  \end{bmatrix}}

\newcommand{\qM}{\bm{q}_M}  
\newcommand{\qMlin}{{q}_M^{lin}}  
\newcommand{\qMnon}{{q}_M^{non}}  

\newcommand{\qMone}{{q_{M}}_1}
\newcommand{\qMtwo}{{q_{M}}_2}

\newcommand{\zero}{\boldsymbol{0}}

\newcommand{\PKone}{\boldsymbol{P}}

\newcommand{\normF}[1]{\left\lVert #1 \right\rVert_F}

\newcommand{\coldos}[2]{\begin{bmatrix} #1 \\ #2 \end{bmatrix}}

\newcommand{\dred}{\bm{d}}  

\newcommand{\dFE}{\bm{d}^{\mathrm{FE}}}  

\newcommand{\dBOUND}{\bm{\hat{d}}}  

\newcommand{\dredLIN}{\bm{d}_L^{lin}}  
\newcommand{\dredNON}{\bm{d}_L^{non}}  

\newcommand{\Tlift}{\bm{T}} 

\newcommand{\Umacro}{\bar{\bm{U}}}
\newcommand{\Gmacro}{\bm{G}^{\mathrm{macro}}}
\newcommand{\GGmacro}{{G}^{\mathrm{macro}}}

\newcommand{\Gmacrox}{\GGmacro_x}      
\newcommand{\Gmacroxy}{\GGmacro_{xy}}  

\newcommand{\muvec}{\bm{\mu}}

\newcommand{\Bmat}{\mathbf{B}}
 \newcommand{\Fext}{\mathbf{F}^{\mathrm{ext}}} 
 \newcommand{\FextD}{\mathbf{F}^{\mathrm{ext}}_{\mathcal{D}}} 

  \newcommand{\Rfe}{\mathbf{R}^{\mathrm{FE}}}

 \newcommand{\FextONE}{\bm{\hat{f}}^{\mathrm{ext}} }
  \newcommand{\FextONEtau}{\bm{\hat{f}}_{\mathrm{ext}}^{\tau} }

\newcommand{\Fint}{\mathbf{F}^{\mathrm{int}}} 
\newcommand{\FintD}{ \displaystyle  \mathbf{F}^{\mathrm{int}}_{\mathcal{D}}} 

\newcommand{\Le}{\mathbf{L}_e}

\newcommand{\Dsnap}{\bm{D}}      
\newcommand{\PhiROM}{\bm{\Phi}}    

\newcommand{\qrom}{\bm{q}}         
\newcommand{\nL}{N}              
\newcommand{\nDOFS}{N_{dof}}              

\newcommand{\nsnap}{P}           
\newcommand{\nsnapNON}{\hat{P}}           

\newcommand{\rROM}{n_{\Phi}}              

 \newcommand{\Bred}[1][]{\bm{B}^{\Phi}_{ #1}} 
 \newcommand{\BredT}[1][]{\bm{B}^{\Phi^T}_{ #1}} 
\DeclareMathOperator{\supp}{supp} 

\newcommand{\wFE}{\bm{w}^{\mathrm{FE}}}           
\newcommand{\rFE}[1]{\bm{r}^{\mathrm{FE}}_{#1}}
\newcommand{\ngp}{M_{\mathrm{gp}}}      
\newcommand{\wFEg}[1]{w^{\mathrm{FE}}_{#1}}   

\newcommand{\Ucub}{\bm{V}}                     
\newcommand{\bcub}{\bm{b}^{\mathrm{FE}}}                     
\newcommand{\Agen}{\bm{A}}               
\newcommand{\AgenLIN}{\bm{A}^{\textrm{lin}}}               
\newcommand{\AgenNON}{\bm{A}^{\textrm{non}}}               

\newcommand{\onevec}{\bm{1}} 
\newcommand{\PKmacro}{\bm{P}^{\mathrm{macro}}}  

\newcommand{\nstress}{s_{\sigma}}        
\newcommand{\nel}{N_{\mathrm{el}}}      
\newcommand{\ngpel}{m_{\mathrm{gp}}^{e}} 

\newcommand{\errdisp}{e_d}              
\newcommand{\errPK}{e_P}                
\newcommand{\nECM}{m_{\mathrm{ECM}}}    
\newcommand{\Pmu}{\mathcal{P}}    
\newcommand{\PmuTRAIN}{\mathcal{P}^{\textrm{tr}}}    
\newcommand{\PmuM}{\mathcal{P}^{\mathcal{M}}}    

\newcommand{\Decod}{\mathcal{\boldsymbol{D}}}      
\newcommand{\Encod}{\mathcal{\boldsymbol{E}}}      
\newcommand{\rD}{n}                  
\newcommand{\PhiD}{\bm{\Phi}_{\mathcal{D}}}        
\newcommand{\rDec}[1]{\bm{r}^{\mathcal{D}}_{#1}}   
\newcommand{\derpar}[2]{\dfrac{\partial #1}{\partial #2}} 

\newcommand{\wDg}[1]{w^{\mathcal{D}}_{#1}}      
\newcommand{\wD}{\bm{w}^{\mathcal{D}}}           
\newcommand{\ZD}{\mathcal{Z}_{\mathcal{D}}}      

\newcommand{\nmu}{n_{\mu}}                       

\newcommand{\PhiM}{\PhiROM_{\!M}}              
\newcommand{\PhiS}{\PhiROM_{\!S}}              

\newcommand{\qS}{\bm{q}_{S}}                   
\newcommand{\Tm}{\bm{H}_{\!M}}                   
\newcommand{\TmNON}{\bm{H}_{\!M}^{non}}                   

\newcommand{\varphiM}{\varphi}                   

\newcommand{\Am}{\bm{Z}_{M}}                   
\newcommand{\Nslave}{\boldsymbol{\mathcal{N}}_{S}}          
\newcommand{\NslaveONE}{\boldsymbol{\mathcal{\widetilde{N}}}_{S}}          

\newcommand{\NslaveONEder}{\boldsymbol{\mathcal{\widetilde{N}}}_{S}'}          

\newcommand{\Qrom}{\bm{Q}_{ROM}}              
\newcommand{\Mtrain}{\bm{M}_{\mu}}             

\newcommand{\Col}{\operatorname{col}}           
\newcommand{\ident}{\bm{I}}                     
\newcommand{\SigM}{\bm{\Sigma}_{M}}             
\newcommand{\VM}{\bm{V}_{M}}                    

\newcommand{\ellRBF}{\bm{\ell}}                 

\newcommand{\nM}{P_{\mathcal{M}}}

\newcommand{\KtanD}{\bm{K}_{\mathcal{D}}}           
\newcommand{\KstdD}{\bm{K}_{\mathcal{D}}^{\mathrm{std}}} 
\newcommand{\KcurvD}{\bm{K}_{\mathcal{D}}^{\mathrm{curv}}} 
\newcommand{\KFEg}[1]{\bm{k}^{\mathrm{FE}}_{#1}} 
\newcommand{\Rdec}{\bm{R}_{\mathcal{D}}} 

\newcommand{\Zomega}{\mathcal{Z}_{\omega}}

\newcommand{\nomega}{m}
\newcommand{\mLower}{m_{\textrm{min}}} 

\newcommand{\omegag}[1]{\omega_{#1}}

\newcommand{\Uinv}{\boldsymbol{U}_{\mathrm{inv}}}

\newcommand{\wDZ}{\wD_{\ZD}}

\newcommand{\mInit}{m_0}

\DeclareMathOperator{\ncol}{ncol} 

\newcommand{\omegaInit}{\bm{\omega}_0}

\newcommand{\ncond}{n_c}
\newcommand{\ncondQ}[1]{\bar{n}_c^{#1}} 

\newcommand{\defeq}{\mathrel{\mathop:}=}

\newcommand{\ninv}{n_{\mathrm{inv}}}

\newcommand{\rj}{r_j}

\newcommand{\Ubarj}[1]{\overline{\bm{U}}_{#1}}

\newcommand{\Uj}[1]{\bm{U}_{#1}}

\newcommand{\bj}[1]{\bm{b}_{#1}}

\newcommand{\omegaj}[1]{\bm{\omega}_{#1}}

\newcommand{\qMj}[1]{\qM^{(#1)}}

\newcommand{\Wad}{\boldsymbol{\mathcal{W}}}

\newcommand{\Reg}{\boldsymbol{\mathcal{R}}}

\newcommand{\Iact}{\mathcal{I}_{\mathrm{act}}}

\newcommand{\Irem}{\mathcal{I}_{\mathrm{rem}}}
\newcommand{\IremLOC}{\mathcal{\bar{I}}_{\mathrm{rem}}}

\newcommand{\wold}[1]{\bm{\omega}^{(#1)}_{\mathrm{old}}}

\newcommand{\wnew}[1]{\bm{\omega}^{(#1)}_{\mathrm{new}}}

\newcommand{\UTact}[1]{{\Uj{#1}}_{\Iact}^{T}}

\newcommand{\UremDEFt}[1]{({\Uj{#1}}_{\Irem})^T}

\newcommand{\KG}{\boldsymbol{\mathcal{K}}_{\mathcal{G}}}

\newcommand{\Glat}{\mathcal{G}_{\mathcal{M}}}

\newcommand{\Wold}{\Wad_{\mathrm{old}}}

\newcommand{\Wnew}{\Wad_{\mathrm{new}}}

\newcommand{\Woldrem}{\overline{\Wad}_{\mathrm{old}}}

\newcommand{\alphaG}{\alpha_{\mathcal{G}}}

 \newcommand{\woldremDEF}[1]{{{\wold{#1}}_{\IremLOC}}}
\newcommand{\woldrem}[1]{\bm{\bar{\omega}}^{(#1)}_{\mathrm{old}}}

\newcommand{\KGs}{\boldsymbol{\mathcal{K}}_{\mathcal{G}}}
\newcommand{\KGsM}{\boldsymbol{\mathcal{\hat{K}}}_{\mathcal{G}}}
 \newcommand{\EdgesM}{\mathcal{E}_{\mathcal{M}}}

 \newcommand{\dwj}[1]{\Delta\bm{\omega}^{(#1)}}

\newcommand{\znew}{\bm{z}_{\mathrm{new}}}

\newcommand{\zold}{\bm{z}_{\mathrm{old}}}

\newcommand{\zpart}{\bm{z}_{\mathrm{p}}}

\newcommand{\Hquad}{\bm{H}}

\newcommand{\gquad}{\bm{g}}

\newcommand{\Nmat}{\bm{N}}

\newcommand{\yred}{\bm{y}}

\newcommand{\Dj}[1]{\mathcal{D}^{(#1)}}

\newcommand{\Fj}[1]{\mathcal{F}^{(#1)}}

 \newcommand{\wpfull}[1]{\bm{\omega}_{\mathrm{p,full}}^{(#1)}}
 \newcommand{\Nj}[1]{\bm{N}^{(#1)}}

 \newcommand{\SUCCESS}{\chi_{\mathrm{feas}}}
\newcommand{\TRUE}{\mathtt{true}}
\newcommand{\FALSE}{\mathtt{false}}
\newcommand{\ENFORCE}{\chi_{\mathrm{enf}}}

 \newcommand{\Psij}[1]{\boldsymbol{\Psi}^{#1}}
\newcommand{\Sigmaj}[1]{\boldsymbol{\Sigma}^{#1}}
\newcommand{\Vj}[1]{\mathbf{\bar{V}}^{#1}}
\newcommand{\SVD}{\operatorname{SVD}}
\newcommand{\tolSVD}{\varepsilon^A}
\newcommand{\tolSVDfixed}{\tolSVD}

\newcommand{\dwFj}[1]{\Delta\mathbf{\tilde{w}}_{F}^{#1}}

\newcommand{\AND}{\;\mathtt{AND}\;}

\newcommand{\kact}{k_{\mathrm{act}}}
\newcommand{\ploc}{ \bar{p} }

\newcommand{\loc}{\operatorname{loc}}
\newcommand{\Uactj}[1]{\widetilde{\mathbf{U}}^{#1}}

\newcommand{\wFj}[1]{\mathbf{w}_{F}^{#1}}
\newcommand{\btildej}[1]{\widetilde{\mathbf{b}}^{#1}}
\newcommand{\diag}{\operatorname{diag}}
\newcommand{\Sdir}{S_{\mathrm{Dir}}} 

\newcommand{\ntry}{n_{\mathrm{try}}} 
\newcommand{\Iord}{\mathcal{I}_{\mathrm{ord}}} 

\newcommand{\Jcand}{\mathcal{J}_{p}} 
\newcommand{\Jrem}{\mathcal{J}_{I}} 

\newcommand{\Qw}{\mathcal{Q}_{W}} 
\newcommand{\Rsdir}{\mathcal{R}_{\mathrm{Dir}}} 

\newcommand{\Wadc}{\Wad} 

\newcommand{\fload}{\mu}
\newcommand{\floadNON}[1]{\fload_{#1}^{\mathrm{non}}} 
\newcommand{\floadmax}{\mu^{\mathrm{max}}}

\newcommand{\stress}{\bm{\sigma}}           
\newcommand{\strain}{\bm{\varepsilon}}      

\newcommand{\youngD}{E}                     
\newcommand{\poissonD}{\nu}                 

\newcommand{\strengthD}{\sigma_{\mathrm c}} 
\newcommand{\hardeningD}{H}                 

\newcommand{\Celas}{\bm{\mathbb C}}

\newcommand{\tolSVDd}{\varepsilon^D} 
\newcommand{\DsnapLIN}{\Dsnap^{\mathrm{lin}}}      
\newcommand{\DsnapNON}{\Dsnap^{\mathrm{non}}}      
\newcommand{\DsnapNONorth}{ {\bm{\tilde{D}}}_L^{\mathrm{non}}}      

\newcommand{\PhiROMlin}{\boldsymbol{\Phi}^{\mathrm{lin}}}      
\newcommand{\PhiROMnon}{\boldsymbol{\Phi}^{\mathrm{non}}}    

 \newcommand{\Rn}[1]{\mathbb{R}^{#1}}    
  \newcommand{\RRn}[2]{\mathbb{R}^{#1 \times #2}}  
\newcommand{\damage}{d_{\epsilon}} 
 \newcommand{\rDAMAGE}{r_{\epsilon}} 
 \newcommand{\qDAMAGE}{q_{\epsilon}} 
\newcommand{\matcdos}[4]{\begin{bmatrix} #1 & #2   \\ #3 & #4  \end{bmatrix}}

\newcommand{\qMnonrel}{\widehat{q}^{\mathrm{\,non}}_{M}} 

\newcommand{\qNONall}{\bm{q}^{\mathrm{non}}} 
\newcommand{\qNONallrel}{\bm{\hat{q}}^{\mathrm{non}}} 

\newcommand{\QNONrel}{\widehat{\bm{Q}}_{\mathrm{non}}} 
 \newcommand{\normd}[1]{\Vert #1 \Vert}  

\newcommand{\Qnon}{\bm{Q}_{\mathrm{non}}}
\newcommand{\Qlin}{\bm{Q}_{\mathrm{lin}}}

\newcommand{\gSnon}{\bm{g}_{S}^{\mathrm{non}}}
\newcommand{\ggSnon}[1]{{g}_{S,#1}^{\mathrm{non}}}

\newcommand{\TAU}{\bm{\tau}}

\newcommand{\RTauD}{\bm R_{\tau}^{D}}
\newcommand{\JTAU}{\bm{J}_{\tau}}

 \newcommand{\Uel}{\mathbf{U}_{\mathrm{el}}} 
\newcommand{\KweigD}{\bm{K}_{\mathcal{D}}^{\omega}} 

\newcommand{\nOUT}{n_{{y}}}

\newcommand{\yOUT}{\bm{y}_{{y}}}

\newcommand{\rOUT}{\bm{r}^{{y}}}

\newcommand{\rOUTg}[1]{\bm{r}^{{y}}_{#1}}
\newcommand{\rGENg}[1]{\bm{r}_{#1}} 

\newcommand{\PROMPT}[1]{ }
\newcommand{\MATLAB}[1]{ }

\makeatletter
\newcommand{\CIMNEaffiliationnote}{\gdef\@fnmark{\(\dagger\)}}
\makeatother

\begin{document}
\begin{frontmatter}


\title{Dimensional hyperreduction of nonlinear finite element models via
empirical cubature with manifold-adaptive weights}

 \author[rvt,els]{Joaqu\'{i}n A. Hern\'{a}ndez\corref{cor1}}
\ead{jhortega@cimne.upc.edu}

\author[rvt,stan]{S. Ares de Parga\CIMNEaffiliationnote}

\author[rvt,elss]{Riccardo Rossi}

\nonumnote{\textsuperscript{\(\dagger\)} This work was conducted while S. Ares de Parga was affiliated with CIMNE.}
\cortext[cor1]{Corresponding author}

\address[rvt]{Centre Internacional de M\`{e}todes Num\`{e}rics en Enginyeria (CIMNE), Barcelona, Spain}

\address[elss]{Universitat Polit\`{e}cnica de Catalunya, Department of Civil and Environmental Engineering (DECA), Barcelona, Spain}

\address[els]{Universitat Polit\`{e}cnica de Catalunya, Departament de Resist\`{e}ncia de Materials i Estructures a l'Enginyeria, E.S. d'Enginyeries Industrial, Aeroespacial i Audiovisual de Terrassa (ESEIAAT), Terrassa, Spain}

\address[stan]{Department of Aeronautics and Astronautics, Stanford University, Stanford, CA 94305, USA}

\begin{abstract}
Nonlinear-manifold reduced-order models for parametrized finite
element problems can achieve substantial compression both in the number
of generalized (latent) coordinates and, through sampling-and-weighting
hyperreduction, in the number of sampled elements/integration points.
Yet current sampling-and-weighting approaches employ weights that remain
fixed over the entire solution manifold. We contend that this restriction
leaves part of the hyperreduction potential untapped: allowing the weights
to vary continuously and nonlinearly with the latent coordinates can
further decrease the number of spatial entities that must be sampled.  To exploit this possibility, we propose the Manifold-Adaptive-Weight
Empirical Cubature Method (MAW--ECM). Starting from a feasible fixed-weight ECM
rule, a greedy pruning strategy progressively removes sampled entities
through a sequence of convex quadratic weight-redistribution problems
that enforce the local   conditions and positivity.
The methodology is assessed on two nonlinear benchmarks: homogenization of a metamaterial unit cell exhibiting negative
incremental stiffness, and a history-dependent continuum-damage problem. In both cases, the nonlinear manifold is constructed from an
initial linear compression followed by an input-informed identification of
the latent coordinates as general linear combinations of the retained
modal amplitudes, incorporating graph information when relevant to seek
the intrinsic dimensionality of the solution manifold. We show that  combining the nonlinear-manifold representation with MAW--ECM reduces
the number of sampled integration points by more than two orders of
magnitude relative to the corresponding standard linear reduced model.
 Relative to the fixed-weight manifold models alone, the adaptive weights
eliminate approximately \(80\%\) of the remaining  points in
the homogenization benchmark and more than \(97\%\) in the damage
benchmark, while essentially preserving their accuracy.

\end{abstract}
%


\end{frontmatter}



\section{Introduction}



\subsection{Manifold HROMs using standard sampling-and-weighting hyperreduction}

Recent years have witnessed rapid development of projection-based
hyperreduced-order models built on nonlinear solution manifolds,
hereafter referred to as manifold HROMs. In these
models, the high-dimensional state is represented as a nonlinear
function of a small set of latent coordinates.
Galerkin formulations enforce the governing equations on the tangent
space of the manifold
\cite{lee2019model,jain2019hyper}, whereas
least-squares Petrov--Galerkin   formulations determine the
reduced state through residual minimization
\cite{lee2019model,barnett2022quadratic,barnett2023neural,chmiel2025unified,
ares2026nonlinear}. These approaches constitute a natural evolution
of classical projection-based HROMs \cite{farhat2014dimensional,hernandez2017dimensional}, in which the state is confined
to a fixed linear space, typically obtained from solution snapshots
through proper orthogonal decomposition (POD), computed in turn via the     singular value decomposition (SVD).
Among the various strategies encompassed by the term
``hyperreduction'',   we focus here on
  \emph{sampling-and-weighting} procedures, in which the full-order evaluation
  of the nonlinear terms entering the reduced equations is
replaced by a positive linear combination of a   subset of local
contributions. This class includes the
Energy-Conserving Sampling and Weighting method (ECSW)
\cite{farhat2014dimensional,farhat2015structure},
the Empirical Cubature Method (ECM)
\cite{hernandez2017dimensional,hernandez2020multiscale}, and the
Empirical Quadrature Procedure (EQP)
\cite{yano2019lp,yano2019discontinuous}.

 Of the two levels of reduction involved in manifold HROMs---the
nonlinear parametrization of the state and the subsequent
hyperreduction by sampling-and-weighting---research efforts have so far
been concentrated almost entirely on the former. Representative examples (the list is not exhaustive) include
methods that learn the nonlinear manifold directly in the
high-dimensional state space via   autoencoders
\cite{lee2019model,fresca2021comprehensive}, as well as approaches that
first compress the state linearly, typically by POD, and subsequently
construct a nonlinear parametrization in the resulting reduced
coordinate space
\cite{fresca2022pod,barnett2022quadratic,barnett2023neural,
 ares2026nonlinear,koike2026sparse}. Despite their different architectures, all these approaches share the same objective: to identify the smallest set possible  of latent variables and reconstruct with maximum accuracy the high-dimensional training states from them.
%

 By contrast, the second compression stage has so far been carried out
essentially with the same sampling-and-weighting machinery developed
for standard HROMs; in particular, the manifold HROM formulations in
Refs.~\cite{jain2019hyper,barnett2023neural,chmiel2025unified,
ares2026nonlinear} all rely on ECSW. Neither specific algorithms nor
refinements of existing ones have been devised to exploit the new
nonlinear-manifold framework. ECSW (and, in essence, the same applies to ECM and EQP) may be
viewed as seeking a sparse vector of nonnegative weights whose local
contributions reproduce, within a prescribed accuracy, those of the
high-dimensional model over the training set. These contributions are
assembled into a matrix that collects, for every candidate spatial
entity, the projected contributions generated over the sampled states.    Passing from a linear HROM
to a manifold HROM reduces the number of such conditions generated at
each state (and hence the size of the matrix), because projection is now performed onto a much smaller
state-dependent tangent space; this reduction in the number of
conditions per state helps explain the smaller sampled meshes reported
for manifold HROMs using ECSW, in comparison with standard HROMs, in
Refs.~\cite{jain2019hyper,barnett2023neural,ares2026nonlinear,
chmiel2025unified}.
Yet the underlying approach remains essentially unchanged: a
single weight vector is required to satisfy the conditions from all
training states simultaneously. Since the constraints are linear,  the resulting rule should be able to reproduce, not only the conditions obtained from the training data, but the entire linear span of such conditions, including combinations that need not correspond to any state on the solution manifold.

\subsection{Motivation for using adaptive weights}
 It is precisely this global linear character of the standard
sampling-and-weighting construction that we challenge in the present paper.
 Our central claim is that this standard approach does not fully exploit the
nonlinear-manifold structure underlying the reduced model.
 If the state is parameterized on a curved manifold, and
quantities entering the reduced equations ---such as the test functions
in a Galerkin projection--- change as the solution moves over that
manifold, then the positive weights used for hyperreduction should
change as well, rather than remain fixed over all latent states. The
expected consequence of this adaptiveness is a further reduction in
the number of spatial entities, such as elements, control volumes, or
integration points, that must be retained to reproduce the
full-order contributions entering the reduced equations.

 The mechanism by which adaptive weights can reduce the number of
required spatial entities can be readily illustrated by a simple
quadrature problem\footnote{Sampling-and-weighting hyperreduction can be
interpreted, in essence, as a quadrature/cubature problem, as made
explicit by the ECM and EQP
\cite{hernandez2017dimensional,yano2019lp}, as well as by the pioneering
work of Ref.~\cite{an2009optimizing}}. Consider the one-dimensional family
\(f_q(x)=x^q\), with \(x\in[0,1]\) and
\(q\in\{0,1,\ldots,P\}\), whose exact integral is
\(I(q)=\int_0^1 f_q(x)\,dx=1/(q+1)\). Here, the monomial exponent
\(q\) may be regarded as a discrete latent coordinate, while any
\(x\in[0,1]\) is a potential sampling point. It can be readily seen that, if a single fixed-weight rule is required
to integrate every latent state exactly,  $I(q) = \sum_{g=1}^{m} f_q(x_g)\,\omega_g$,
then it must also integrate
their complete linear span, namely, the polynomial space
of degree at most \(P\). Even the optimal Gauss--Legendre rule therefore
requires \(m\geq (P+1)/2 \) sampling points---the
required number of integration points grows   linearly with the number
of latent states. By contrast, if one retains a single sampling point,
say \(x_1=1\), and allows its weight to depend on the latent coordinate,
then \(I(q)=\omega_1(q)f_q(x_1)\), with \(\omega_1(q)=1/(q+1)\); thus, one
point suffices for every \(q\), irrespective of \(P\).  The gain in efficiency (from \(m\geq (P+1)/2 \) points to just $m=1$ point ) comes thus
  from allowing the weight to adapt to the latent state, as asserted.

In a more practical context, further evidence that adaptive weights can
substantially reduce the number of sampled spatial entities is provided
by our previous work~\cite{bravo2024subspace} on
local HROMs. In local HROMs, the solution manifold is represented
piecewise by a collection of local linear subspaces, each associated
with a different region of the solution space. Sampling-and-weighting
hyperreduction in this setting was first explicitly addressed in Ref.~\cite{grimberg2021mesh}. Building on this setting, and
using the Empirical Cubature Method, ECM, as the basic point-selection
procedure, our proposal in Ref.~\cite{bravo2024subspace} was to endow
each subspace with its own nonnegative weights while retaining a common
set of cubature (i.e., integration) points; we termed the resulting
approach the Subspace-Adaptive Weights ECM (SAW--ECM). Numerical
results  in Ref.~\cite{bravo2024subspace} showed that, as the local representation is progressively
refined, the required common support decreases markedly, approaching
in the limiting regime a size governed by the  dimension of the parameter space rather than by the dimension of the global span of all
subspaces.

\subsection{Manifold-adaptive weights ECM}
\label{sec:mawECMintro}

The main novel contribution of the present work is the extension of the
SAW--ECM described above from a discrete collection of subspaces to a
continuously parameterized nonlinear manifold. In direct analogy with
the terminology of Ref.~\cite{bravo2024subspace}, we term the resulting
method the \emph{Manifold-Adaptive Weights ECM} (MAW--ECM).

In passing from a discrete collection of subspaces in the SAW--ECM to a
continuous manifold in the proposed MAW--ECM, the subspace index associated to each local subspace is replaced by the latent
coordinates, so that each spatial entity is endowed with a
positive weight field over the latent space. This transition introduces, however, an essential additional
requirement: the corresponding weight fields should be as smooth as
possible over the latent space.  Such smoothness is important because the discrete weight samples must
ultimately be represented by regression over the latent domain; strong
local oscillations would compromise both the accuracy of this
regression and the robustness of the nonlinear solution procedure
(Newton--Raphson in the examples considered here).   To enforce this smoothness requirement, we propose to modify the sparsification-type optimization problem of
Ref.~\cite{bravo2024subspace} by adding a   graph-based
 regularization term into the objective function, thereby penalizing irregular variations
of the weight fields.

As in Ref.~\cite{bravo2024subspace}, and also in the standard ECSW,
ECM, and EQP, the resulting sparsification problem is NP-hard
\cite{natarajan1995sparse} because of the presence of the \(\ell_0\) pseudo-norm,
which   measures the number of nonzero entries of the weight vectors. Practical approaches to \(\ell_0\)-minimization
typically rely either on suboptimal greedy or active-set procedures,
as in ECM and ECSW, or on suitable convex relaxations, most commonly
obtained by replacing the \(\ell_0\) objective with an
\(\ell_1\)-based surrogate promoting sparsity, as in EQP \cite{yano2019discontinuous,yano2019lp}. The latter
strategy was explored by the authors in
Ref.~\cite{bravo2024subspace} for the SAW--ECM. The main conclusion
drawn there, which carries over to the problem at hand, is that this
\(\ell_1\)-based convexification route is not viable in practice, since
it requires stacking all adaptive weight vectors into a single global
optimization variable from the outset and, in doing so, fails to
exploit the fundamentally local structure of the constraints
associated with each sampled latent state.

Accordingly, the present work does not pursue this \(\ell_1\)-based route, but
instead adopts a greedy iterative strategy that avoids treating the entire
global optimization variable simultaneously at all iterations. The underlying idea is inspired by
the \emph{pruning} strategies employed in the derivation of generalized
quadrature rules \cite{bremer2010nonlinear,xiao2010numerical,hernandez2024cecm},
where the goal
is likewise to identify the smallest possible set of integration points,
including both their positions and weights.   The analogy with the present MAW--ECM setting is the following. In both
cases, sparsification starts from an already feasible cubature rule and
proceeds through successive pruning steps in which integration points are
removed while the remaining weights are redistributed so as to preserve the
constraints.
The essential difference lies in the nature of the design variables. In
generalized quadrature, such as in the  continuous version of the ECM proposed by the authors in Ref.~\cite{hernandez2024cecm},   both the positions and the weights of the
cubature points evolve continuously during pruning. In the present MAW--ECM
setting, by contrast, the candidate entities are fixed in space  (i.e., Gauss points inherited
from the FE discretization), and only the manifold-dependent weights are
modified.   At each tentative elimination, redistributing the remaining weights leads
to a convex quadratic program with linear equality and positivity
constraints. A key ingredient of the proposed procedure is an algorithm
that exploits the local structure of the state-wise constraints and
activates the graph regularization---the most computationally demanding
component of the method---only when required to preserve positivity.

  \subsection{Input-informed construction of the nonlinear manifold}
 \label{sec:inputinformed}
The proposed methodology is assessed on two nonlinear finite element (FE) solid-mechanics
benchmarks: a large-strain hyperelastic homogenization problem (Section~\ref{sec:metamaterial}) and a
history-dependent damage problem (Section~\ref{sec:DamageProblem}).
A prerequisite for applying the
MAW--ECM is the identification of suitable latent coordinates and the
corresponding nonlinear function between high-dimensional state and such latent variables (we call it here the  \emph{decoder}).   As in Refs.~\cite{barnett2023neural,ares2026nonlinear,koike2026sparse},
we first apply POD to obtain a linear compression of the state and then
construct the nonlinear decoder in terms of the resulting modal
coefficients. Our approach differs in the definition of the latent
coordinates.  Whereas
Refs.~\cite{barnett2023neural,ares2026nonlinear} identify them with the
leading POD coefficients, and Ref.~\cite{koike2026sparse} selects an
informative subset not necessarily following the singular-value
ordering, here the latent coordinates are allowed to be arbitrary
linear combinations of the retained modal coefficients.
Our proposal is that the identification of these linear combinations should also exploit
the associated problem inputs, as well as, when relevant,  their relative arrangement,
encoded through a graph over the training states. More specifically, in the
hyperelastic problem of Section~\ref{sec:metamaterial}, the input parameters
are used directly in a least-squares identification, while in the damage
problem of Section~\ref{sec:DamageProblem}, the nonlinear component of the
latent vector is identified from the damage-evolving states through a
graph-based eigenvalue problem informed by the loading history. In both cases, we manage to identify the intrinsic
dimensionality of the solution manifold, i.e., the minimum number of
independent latent coordinates required for its representation.

 To the best of the authors' knowledge, this input-informed approach to
constructing the POD-based decoder constitutes an original contribution
of the present work, particularly for history-dependent deformations,
where the loading history, its graph structure, and constitutive
information guide the identification of the latent coordinates. A related recent work  \cite{zhang2026unified}   has considered    manifold HROMs for
history-dependent problems, with POD-based  decoder  and fixed-weight ECSW hyperreduction;  however, they do not seek an
intrinsically minimal latent representation (  latent variables are chosen as leading POD coefficients), nor do they construct the
decoder using problem-specific graph-based, input information   (the nonlinear component is   quadratic ).

%


%
%
%

\subsection{Scope}

We choose parametrized, quasi-static nonlinear solid finite element
problems subject to affine kinematic constraints as the vehicle for
presenting and assessing the proposed methodology. More specifically, the
formulation is cast in terms of interior Gauss-point sampling (i.e., as a
true cubature problem), with the projected internal forces and, when
applicable, volumetric quantities of interest as the terms to be
hyperreduced. These choices are not intrinsic to MAW--ECM, which can in
principle be extended to other scenarios, including dynamic solid
problems, other governing equations, and discretizations other than
finite elements, provided that the corresponding local contributions and
invariant constraints to be preserved by the adaptive rule are defined
appropriately. Likewise, within the finite element setting, Gauss-point
selection is only one possible realization: MAW--ECM may also be posed as
a mesh-sampling procedure, as already done with fixed-weight ECM for
standard HROMs in Refs.~\cite{ares2023hyper,bravo2024geometrically}.

As for the input-informed POD-based decoders, the particular
constructions used in the two benchmarks are tailored to the problems
considered, and no claim is made here regarding their generality as
manifold-learning procedures. Finally, given the exploratory character of this study, the numerical
assessment is deliberately restricted to low-dimensional parameter and
latent spaces, for which the adaptive weight fields  can be directly visualized and interpreted.
Scalability to substantially higher-dimensional spaces, particularly
regarding the construction of the latent graph used for regularization, lies
outside the scope of the present work.

%


%
%
%

\section{Full-order finite element model}
\label{sec:FOM}



Let us  consider a general,  parametrized, quasi-static nonlinear solid FE problem
subject to affine kinematic constraints, in which  the solution
  depends on an input parameter vector
\(\muvec\in\Pmu\subset\mathbb{R}^{\nmu}\), which may
represent   prescribed displacements  and/or  applied forces.   In history-dependent problems, this
dependence is understood incrementally, with the current response also
depending on the previously converged internal variables. The FE nodal displacement vector $\dFE \in \Rn{\nDOFS}$ is written as
\begin{equation}
\label{eq:dvecreconstruction}
 \dFE = \Tlift \dred  + \dBOUND ,
\end{equation}
where \(\dred\in\mathbb{R}^{\nL}\)  ($\nL < \nDOFS$) collects the independent degrees of freedom,
\(\Tlift\) is the lifting operator associated with the homogeneous constrained
space, and \(\dBOUND=\dBOUND(\muvec)\) accounts for the affine prescribed
kinematic data. We define the full finite element residual $\Rfe \in \Rn{\nDOFS}$ as
\begin{equation}
\label{eq:Rfe_definition}
    \Rfe(\dred;\muvec)
    =
    \Fint(\dred;\muvec)-\Fext(\muvec),
\end{equation}
where \(\Fint\) and \(\Fext\) are the nodal internal and external force
vectors, respectively (the latter is assumed to be independent of $\dred$).  The constrained parametrization~\eqref{eq:dvecreconstruction} induces the
admissible virtual displacements
\begin{equation}
\label{eq:dvecvariations}
    \delta\dFE=\Tlift\,\delta\dred .
\end{equation}
Substitution of the preceding expression  into the discrete virtual-work statement
\begin{equation}
\label{eq:ResFE}
 {\delta\dFE}^T\Rfe=0,
\end{equation}
and exploiting the arbitrariness of $\delta \dred$, leads to the FE equilibrium conditions
\begin{equation}
\label{eq:FOM_equilibrium_lifted}
         \Tlift^T (\Fint - \Fext)  = \bm{0}.
 \end{equation}

\subsection{Internal forces}
\label{sec:FOM_gauss_decomposition}

By hypothesis,  the only term in the above equilibrium equations that
depends on the unknown displacement vector $\dred$ is the nodal internal force $\Fint$. This is
therefore the term to be approximated by hyperreduction via the proposed adaptive  empirical cubature rule.  Since  such a  rule will be constructed from Gauss-point
contributions, it is necessary to derive the expression for $\Fint$ so that these
contributions appear explicitly.  We start from its total Lagrangian expression
\begin{equation}
\label{eq:Fint_continuous_general}
 \Fint
 =
 \int_{\Omega_0}
 \Bmat^T\,\PKone\,d\Omega,
\end{equation}
where \(\Omega_0\) denotes the reference domain, \(\PKone \in \Rn{\nstress}\) is the first
Piola--Kirchhoff stress vector  ($\nstress = 9$ for 3D and $\nstress = 4$ for 2D ), and \(\Bmat \in \RRn{\nstress}{\nDOFS}\) is the global operator
relating nodal displacement variations to the corresponding variations of the
deformation gradient. In the damage
problem considered in Section~\ref{sec:DamageProblem}, the formulation is
small-strain; consequently, in this case,  \(\PKone\) is to be replaced by the Cauchy stress
vector, and \(\Bmat\) becomes the standard finite element
strain--displacement matrix. We next expand the integral in Eq.~\refeq{eq:Fint_continuous_general} as the sum of the   contributions of each finite element, and then approximate  such element contributions by Gauss quadrature; this leads to:
\begin{equation}
\label{eq:Fint_FE_general}
    \Fint
    =
    \sum_{e=1}^{\nel}
   \Le^T
    \left(
    \sum_{i=1}^{\ngpel}
    \Bmat_{e,i}^T \, \PKone_{e,i} \, w_{e,i}
    \right).
\end{equation}
Here \(\Le\) denotes the Boolean assembly operator of element \(e\),
\(\Bmat_{e,i}\) is the element-level counterpart of \(\Bmat\) evaluated at
Gauss point \(i\), \(\PKone_{e,i}\) is the corresponding stress value, and
\(w_{e,i}\) is the physical quadrature weight, including both the parent-domain
Gauss weight and the Jacobian determinant of the isoparametric mapping.

Finally, we reindex the element/Gauss-point pairs \((e,i)\) by a
single global index \(g=1,\ldots,\ngp\) (here \(\ngp\) denotes the total number of Gauss points).   This numbering defines a
bijection between \(g\) and the corresponding pair
\((e(g),i(g))\), where \(e(g)\) identifies the element containing the
point and \(i(g)\) its local Gauss-point index. The constrained internal
force appearing in Eq.~\refeq{eq:FOM_equilibrium_lifted} can then be written as
\begin{equation}
\label{eq:projected_internal_force_gp}
    \Tlift^T\Fint
    =
    \sum_{g=1}^{\ngp}
    \rFE{g}\,\wFEg{g},
\end{equation}
where
\begin{equation}
\label{eq:rFE_general}
    \rFE{g}(\dred;\muvec)
    :=
    \Tlift^T
    \mathbf{L}_{e(g)}^T
    \Bmat_{e(g),i(g)}^T
    \PKone_{e(g),i(g)}(\dred;\muvec)
\end{equation}
stands for the internal-force density associated with global Gauss
point \(g\), whereas
\(\wFEg{g}=w_{e(g),i(g)}>0\) denotes its integration weight.

\subsection{External actions}
\label{sec:FOM_external_forces}

The vector of external actions
\(\Fext\in\mathbb{R}^{\nDOFS}\) in Eq.~\refeq{eq:Rfe_definition} contains the contributions   induced by
prescribed tractions on the Neumann boundary and by body forces. As already pointed out, these external actions are
assumed here to be independent of the unknown displacement vector \(\dred\), and therefore will not be subjected to hyperreduction via the proposed  cubature scheme. Since in the problem under consideration the input parameters control the amplitude of the external forces, we can write the projection of this term into the constrained space as a linear function of $\muvec$, that is:
\begin{equation}
\label{eq:fextlift_linear}
    \Tlift^T\Fext(\muvec)
    =
    \FextONE\,\muvec ,
\end{equation}
 the columns of \(\FextONE \in \RRn{\nL}{\nmu}\) being  the nodal patterns
associated with unit values of the components of \(\muvec\).

%

Finally, substituting Eqs.~\refeq{eq:projected_internal_force_gp}
and~\refeq{eq:fextlift_linear} into
Eq.~\refeq{eq:FOM_equilibrium_lifted} gives the full-order equilibrium
equations in a form that makes explicit the loop over all Gauss points of the
mesh:
\begin{equation}
\label{eq:FOM_equilibrium_general}
      \overbrace{\sum_{g=1}^{\ngp}
    \rFE{g}(\dFE;\muvec)\,\wFEg{g}}^{\textrm{Internal forces}}
    -
    \overbrace{\FextONE\,\muvec}^{\textrm{External forces}}
    =
    \bm{0}.
\end{equation}
 This is a system of \(\nL\) nonlinear algebraic equations for the independent
unknowns \(\dred\in\Rn{\nL}\), to be solved for prescribed values of
\(\muvec\). In history-dependent problems, the solution also depends on the
internal variables from the previous pseudo-time step; for simplicity, this dependence is
omitted in the formulation that follows.

%

\subsection{Volumetric outputs of interest}
\label{sec:FOM_outputs}

 Some outputs of interest,   depending nonlinearly on the unknown
\(\dred\), may also involve volumetric integration and may therefore
benefit from the same cubature viewpoint. We denote by
\(\yOUT\in\mathbb{R}^{\nOUT}\) a generic output of interest of this type
and write its full finite element quadrature form as
\begin{equation}
\label{eq:generic_output_FE_quadrature}
    \yOUT
    =
    \int_{\Omega_0}
    \rOUT \,d\Omega
    =
    \sum_{g=1}^{\ngp}
    \rOUTg{g} \,\wFEg{g},
\end{equation}
where \(\rOUT\in\mathbb{R}^{\nOUT}\) is the volumetric integrand defining
the output of interest, and \(\rOUTg{g}(\dred;\muvec)\) denotes its
evaluation at Gauss point \(g\). When such an output is present, the
cubature rule used for hyperreduction is also applied to
\eqref{eq:generic_output_FE_quadrature}, so that both the internal-force
term in Eq.~\refeq{eq:FOM_equilibrium_general} and the output functional
are approximated over the same reduced set of integration points. This is
the case of the metamaterial benchmark of Section~\ref{sec:metamaterial},
where the output of interest is the volumetric average of the first
Piola--Kirchhoff stress \(\PKone\).

\section{Manifold hyperreduced-order model}

The goal in this section is to construct the hyperreduced-order counterpart of the equilibrium equation    ~\refeq{eq:FOM_equilibrium_general}. This involves two approximations. The first is a kinematic compression:  we approximate the unknown displacement $\dred \in \Rn{\nL}$ as a nonlinear function of a vector of  latent coordinates\footnote{The subscript \(M\) in $\qM$ stands for ``master'', anticipating the
master--slave partition of modal basis coordinates    used later to
construct the benchmark-specific decoders.  This
terminology is borrowed from nonlinear structural dynamics
\cite{touze2021model,vizzaccaro2020comparison}. Other terms found in the related literature refer to the master
coordinates as active~\cite{zuo2026nonlinear,callaham2022role} or
primary~\cite{ares2026nonlinear}, and to the remaining coordinates as
passive \cite{zuo2026nonlinear} or secondary \cite{ares2026nonlinear}.      } $\qM \in \Rn{\rD}$:
\begin{equation}
\label{eq:dredNON}
\dred \approx \Decod(\qM).
\end{equation}
Ideally, the number of latent coordinates $\rD$ should be as close as possible to  the number of input parameters  $\nmu$.  The second approximation is the hyperreduction per se: the loop over the \(\ngp\) Gauss points in
Eq.~\refeq{eq:FOM_equilibrium_general} is to be replaced by a loop over a much
smaller set of selected integration points, with tailored,   positive weights.   The adaptive-weight
formulation introduced later in Section~\ref{sec:MAWECM} will modify only this second reduction step, by
allowing the cubature weights to depend on the latent coordinate $\qM$.
\subsection{Decoder and encoder}
\label{sec:decoderENCODER}

We shall refer to \(\Decod:\Rn{\rD}\rightarrow\Rn{\nL}\) in
Eq.~\refeq{eq:dredNON} as the \emph{decoder}, and to the associated approximate
inverse map
\begin{equation}
\label{eq:inversemapping}
\qM \approx \Encod(\dred),
\end{equation}
as the \emph{encoder}. These maps may be represented through general-purpose autoencoder
architectures, as in Refs.~\cite{lee2019model,fresca2021comprehensive};
here, however, the terminology is used in a broader descriptive sense and
does not prescribe any particular architecture.   The only condition imposed on the decoder is that it be twice
differentiable with respect to \(\qM\). This regularity is required
because the first derivative, represented by the Jacobian matrix:
\begin{equation}
\label{eq:PhiD_definition}
\PhiD
:=
\derpar{\Decod}{\qM}
\in \mathbb{R}^{\nL\times\rD},
\end{equation}
defines
the tangent space used in the projected equilibrium equations, while the
second derivative   enters the consistent Newton--Raphson linearization required by the solution procedure.

\subsection{Factorized decoder and identification of latent coordinates}
\label{sec:fact}
The decoder is constructed offline from displacement snapshots associated
with a representative training set
\begin{equation}
\label{eq:pmutrain}
  \PmuTRAIN = \{\muvec_i\}_{i=1}^{\nsnap}\subset\Pmu.
\end{equation}
We collect these snapshots in
the matrix
\begin{equation}
\label{eq:manifold_snapshot_matrix}
\Dsnap =
\begin{bmatrix}
\dred(\muvec_1) &
\dred(\muvec_2) &
\cdots &
\dred(\muvec_{\nsnap})
\end{bmatrix}
\in\mathbb{R}^{\nL\times\nsnap}.
\end{equation}
 As already pointed out in Section \ref{sec:inputinformed},    we adopt a two-level decoder
representation, combining a preliminary linear compression with a
nonlinear coefficient map. With
this representation, the decoder is expressed in the factorized form
\begin{equation}
\label{eq:decoder_tau_generic}
\Decod(\qM)
=
\PhiROM\TAU(\qM),
\end{equation}
where  \(\PhiROM\in\RRn{\nL}{\rROM}\)  is    an orthonormal basis matrix for a
low-rank approximation of the column space of \(\Dsnap\), and  \(\TAU:\Rn{\rD}\rightarrow\Rn{\rROM}\) maps the latent
coordinates to the coefficients in this linear basis.
As customary, the number of modes $\rROM$ in $\PhiROM$ is controlled by a relative tolerance $\tolSVDd \in [0,1)$ such that
\begin{equation}
\label{eq:PhiROM_truncation}
\frac{
\normF{\Dsnap-\PhiROM\PhiROM^T\Dsnap}
}{
\normF{\Dsnap}
}
\leq
\tolSVDd, \qquad  \PhiROM^T \PhiROM = \ident,
\end{equation}
where $\normF{}$ stands for the Frobenius norm   ($\PhiROM$ may be obtained as the left singular vectors of  a truncated  SVD of \(\Dsnap\), for instance).
If the latent coordinates \(\qM\) have already been identified, the map
\(\TAU\) can be fitted by standard regression procedures. By contrast, if the latent coordinates are not known a priori, they may instead be
learned jointly with the coefficient map through an autoencoder-based
construction, as in Ref.~\cite{fresca2022pod}.

In the examples presented later, we follow the former approach. More specifically,  let
\begin{equation}
\label{eq:qrom_definition}
    \qrom
    \defeq
    \PhiROM^T\dred
\end{equation}
denote the modal coefficients of a displacement state in the linear
basis \(\PhiROM\). Following  \cite{barnett2022quadratic,barnett2023neural,ares2026nonlinear,koike2026sparse}, we seek latent
coordinates   obtained as
linear combinations of  \(\qrom\):
\begin{equation}
\label{eq:qMdefq}
    \qM=\Tm\qrom.
\end{equation}
The problem of selecting the latent coordinates is thus recast as that
of determining
\(\Tm\in\RRn{\rD}{\rROM}\). We assume that \(\Tm\) has full row rank \(\rD\), so that the latent coordinates are linearly independent. References ~\cite{barnett2022quadratic,
barnett2023neural,ares2026nonlinear} use the SVD to determine $\PhiROM$, and  adopt $\Tm=[\,\ident_{\rD}\ \bm{0}\,]$,
where \(\ident_{\rD}\) is the \(\rD\times\rD\) identity matrix; the
latent coordinates are therefore the leading entries of \(\qrom\),
associated with the largest singular values.    Reference~\cite{koike2026sparse}, on the other hand,  allows  \(\Tm\) to select an arbitrary sparse
subset of the modal coefficients.


The identification adopted in the examples discussed later departs from this Boolean selection
strategy, in the sense that the matrix \(\Tm\) is not constrained to
select individual entries of \(\qrom\), but may be dense, allowing each
latent coordinate to combine all the modal coefficients. Its
identification is instead input-informed, in the sense discussed in
Section~\ref{sec:inputinformed}. Further details on the identification of $\Tm$ will be given in Sections \ref{sec:meta_manifold_representation} and \ref{sec:id2}.

%
%

Regardless of the particular expression for $\Tm$, on the approximation manifold, the relation in Eq.~\refeq{eq:qMdefq} together with Eq.~\ref{eq:qrom_definition} defines the following master--slave representation:
\begin{equation}
\label{eq:decoder_master_slave}
    \dred =    \PhiM \Am \qM
    +
    \PhiS \Nslave(\qM).
\end{equation}
 Here, \(\PhiM\in\RRn{\nL}{\rD}\) denotes the matrix of master
modes, whereas
\(\PhiS\in\RRn{\nL}{(\rROM-\rD)}\) contains the slave modes. They are chosen so that
\begin{equation}
\label{eq:orto}
    \PhiM^T\PhiM=\ident,
    \qquad
    \PhiS^T\PhiS=\ident,
    \qquad
    \PhiM^T\PhiS=\zero
\end{equation}
( \(\PhiS\) is  therefore an orthonormal basis for the
orthogonal complement of \(\Col(\PhiM)\) within \(\Col(\PhiROM)\), and hence \([\PhiM\ \PhiS]\) spans \(\Col{(\PhiROM)}\)). On the other hand, $\Am \in \RRn{\rD}{\rD}$ is a scaling matrix. The matrices \(\PhiM\) and \(\Am\) are obtained from the compact SVD of $\PhiROM\Tm^T$; in particular \(\PhiM\) contains the corresponding left singular
vectors, i.e.:
\begin{equation}
\label{eq:master_basis_svd}
    \PhiROM\Tm^T
    =
    \PhiM\SigM\VM^T.
\end{equation}
 whereas
\begin{equation}
\label{eq:Am_definition}
    \Am
    =
    \SigM^{-1}\VM^T.
\end{equation}
The latter expression follows by premultiplying
Eq.~\refeq{eq:decoder_master_slave} by \(\PhiM^T\), using the
orthogonality conditions in Eq.~\refeq{eq:orto}, and combining the
resulting relation with Eqs.~\refeq{eq:qrom_definition},
\refeq{eq:qMdefq} and the SVD above.    The slave amplitudes associated with the training snapshots
are then computed as
\[
    {\qS}_j=\PhiS^T\dred(\muvec_j), \qquad j = 1,2 \ldots \nsnap
\]
and the nonlinear closure relation
\begin{equation}
\label{eq:closure}
     \qS\approx\Nslave(\qM)
\end{equation}
is learned from the training pairs $({\qM}_j,{\qS}_j)$  ($j=1,2 \ldots \nsnap$).

With a slight abuse of notation, we henceforth redefine the reduced
basis as $\PhiROM \leftarrow [\PhiM,\PhiS]$.  This  change of basis leaves
\(\Col{(\PhiROM)}\) unchanged and yields the following particular
expression for the coefficient map in Eq. \refeq{eq:decoder_tau_generic}:
\begin{equation}
\label{eq:tau_master_slave}
    \TAU(\qM)
    =
    \begin{bmatrix}
        \Am\qM\\
        \Nslave(\qM)
    \end{bmatrix},
\end{equation}
and for its Jacobian:
\begin{equation}
\label{eq:JTAU}
\JTAU
:=
\derpar{\TAU}{\qM} = \begin{bmatrix}
                             \Am  \\
                             \derpar{\Nslave(\qM)}{\qM}
                     \end{bmatrix}
\end{equation}.

Lastly, we particularize the encoder introduced in
Eq.~\refeq{eq:inversemapping} for the master--slave decoder in
Eq.~\refeq{eq:decoder_master_slave}. We adopt the explicit
projection-based map
\begin{equation}
\label{eq:encoder_explicit}
    \Encod(\dred)
    \defeq
    \Am^{-1}\PhiM^T\dred.
\end{equation}
 By virtue of the orthogonality conditions in
Eq.~\refeq{eq:orto},  this map is an
exact left inverse of the decoder on the approximation manifold, i.e.: $\Encod\bigl(\Decod(\qM)\bigr)=\qM$.
For states outside the approximation manifold, however, the coordinates
returned by Eq.~\refeq{eq:encoder_explicit} need not coincide with those
obtained from the nonlinear least-squares encoder used, for instance, in
Ref.~\cite{ares2026nonlinear}. We have nevertheless verified numerically
that this discrepancy has no appreciable effect on the results reported
here, and therefore adopt the simplified encoder in
Eq.~\refeq{eq:encoder_explicit} henceforth.

\subsection{Tangent-space projection of the equilibrium equations}

The nonlinear approximation in Eq.~\refeq{eq:dredNON} restricts the
admissible variations of the independent displacement vector to the
tangent space of the manifold:
\begin{equation}
\label{eq:manifold_virtual_variations}
\delta\dred
=
\PhiD(\qM)\,\delta\qM,
\end{equation}
where $\PhiD$ is the Jacobian of the decoder, introduced in Eq.~\refeq{eq:PhiD_definition}. Substituting Eq.~\refeq{eq:manifold_virtual_variations} into
Eq.~\refeq{eq:dvecvariations}, and the resulting expression into the
virtual work statement~\refeq{eq:ResFE}, yields, by virtue of the arbitrariness of
\(\delta\qM\), the reduced-order equilibrium equations:
\begin{equation}
\label{eq:manifold_ROM_equilibrium_FE_weights}
\overbrace{
\sum_{g=1}^{\ngp}
\rDec{g}(\qM;\muvec)\,\wFEg{g}}^{\displaystyle\FintD:\textrm{Internal forces}}
-
\overbrace{
\PhiD(\qM)^T\FextONE\,\muvec}^{\displaystyle\FextD: \textrm{External forces}}
=
\bm{0},
\end{equation}
where
\begin{equation}
\label{eq:rDec_definition}
\rDec{g}(\qM;\muvec)
:=
\PhiD(\qM)^T
\rFE{g}\bigl(\qM;\muvec\bigr)
\in\mathbb{R}^{\rD}
\end{equation}
is the  projected internal force density associated with Gauss
point \(g\).

In the form given above, both terms in
Eq.~\refeq{eq:manifold_ROM_equilibrium_FE_weights} appear, in principle, to require  online operations whose cost depends on the size
of the underlying finite element model. This   dependence arises, on the one hand,
from the products involving the tangent basis
\(\PhiD\in\RRn{\nL}{\rD}\), in the expressions of  both $\FintD$ and $\FextD$, and, on the other hand,  from the summation over all
  Gauss points of the mesh in $\FintD$.

  The former dependence can be eliminated in both terms by
exploiting the factorized decoder in
Eq.~\refeq{eq:decoder_tau_generic}. Indeed,  in terms of this factorization, the
Jacobian matrix defined in Eq.~\refeq{eq:PhiD_definition} adopts the form
\begin{equation}
\label{eq:PhiD_tau_generic}
\PhiD(\qM)
=
\PhiROM\JTAU(\qM),
\end{equation}
where  $\JTAU$ is given in Eq.~\refeq{eq:JTAU}.
Combining Eqs.~\refeq{eq:rDec_definition},
\refeq{eq:rFE_general}, and \refeq{eq:PhiD_tau_generic} (and reverting
temporarily to the local Gauss-point indices \((e,i)\) for clarity), we obtain
\begin{equation}
\label{eq:rDec_reduced_operator}
\rDec{e,i}(\qM;\muvec)
=
\JTAU(\qM)^T
\BredT[e,i]
\PKone_{e,i}\bigl(\qM;\muvec\bigr),
\end{equation}
where
\begin{equation}
\label{eq:Bred_definition}
\Bred[e,i]
:=
\Bmat_{e,i}\Le\Tlift\PhiROM
\in\mathbb{R}^{\nstress\times\rROM}
\end{equation}
 maps variations of the auxiliary
coordinates (amplitude  of the basis vectors $\PhiROM$) to the  displacement gradient components at Gauss point
\(i\) of element \(e\). This matrix is independent of both the latent
coordinates and the input parameters and can therefore be precomputed
offline for every Gauss point. Thus, once \(\JTAU(\qM)\) is available,
evaluation of \(\rDec{e,i}\) requires only the constitutive update at
the corresponding Gauss point---whose kinematic input is likewise
obtained through \(\Bred[e,i]\)---together with matrix--vector products
in spaces of dimensions \(\rROM\) and \(\rD\). No operation involving
the ambient dimension \(\nL\) is required.

The same conclusion holds for the external-force
term \(\FextD\).  Substitution of Eq.~\refeq{eq:PhiD_tau_generic} into the expression for
\(\FextD\) in
Eq.~\refeq{eq:manifold_ROM_equilibrium_FE_weights} gives
\begin{equation}
\label{eq:FextD}
\FextD(\qM;\muvec)
=
\JTAU(\qM)^T\FextONEtau\muvec,
\end{equation}
where
\begin{equation}
\label{eq:FextONEtau}
\FextONEtau
:=
\PhiROM^T\FextONE
\end{equation}
is independent of both the latent coordinates and the parameter values
and can therefore be computed offline.

 The only remaining source of mesh-dependent online cost is thus   the
summation over the complete set of Gauss points required to evaluate
\(\FintD\) in Eq.~\refeq{eq:manifold_ROM_equilibrium_FE_weights} and,
when present, the volumetric output of interest in
Eq.~\refeq{eq:generic_output_FE_quadrature}. As previously mentioned, this bottleneck is addressed
by replacing the full finite element quadrature sums with sums over a
much smaller set of Gauss points, equipped with tailored positive
weights. We first consider a fixed-weight rule constructed by the
Empirical Cubature Method (ECM)
\cite{hernandez2017dimensional,hernandez2020multiscale}. Although its
derivation is similar to that of a standard HROM, it is presented here
because it provides both the conceptual basis and the initialization for
the manifold-adaptive rule  introduced in
Section~\ref{sec:MAWECM}.

%
%

\subsection{Fixed-weight cubature approximation  }
\label{sec:fixedweight_cubature}

We seek the smallest set of Gauss points $\ZD\subset\{1,\ldots,\ngp\}$
and positive weights \(\wDg{g}\), \(g\in\ZD\), such that the resulting
cubature rule reproduces, to prescribed accuracy, the projected
internal forces at the training states
\(\{({\qM}_j,\muvec_j)\}_{j=1}^{\nsnap}\):
\begin{equation}
\label{eq:approx_ECMfixed1}
\FintD({\qM}_j;\muvec_j)
=
\sum_{g=1}^{\ngp}
\rDec{g}({\qM}_j;\muvec_j)\,\wFEg{g}
\approx
\sum_{g\in\ZD}
\rDec{g}({\qM}_j;\muvec_j)\,\wDg{g},
\qquad
j=1,\ldots,\nsnap .
\end{equation}
When the   output of interest in
Eq.~\refeq{eq:generic_output_FE_quadrature} is present, the same rule is
also required to satisfy
\begin{equation}
\label{eq:approxOUTPUT}
\yOUT({\qM}_j;\muvec_j)
=
\sum_{g=1}^{\ngp}
\rOUTg{g}({\qM}_j;\muvec_j)\,\wFEg{g}
\approx
\sum_{g\in\ZD}
\rOUTg{g}({\qM}_j;\muvec_j)\,\wDg{g},
\qquad
j=1,\ldots,\nsnap .
\end{equation}

For the sake of generality, we formulate the above   requirements in
terms of a generic vector-valued  integrand
\(\rGENg{g}({\qM}_j;\muvec_j)\in\Rn{\ncond}\):
 \begin{equation}
\label{eq:approxGEN}
\sum_{g=1}^{\ngp}
\rGENg{g}({\qM}_j;\muvec_j)\,\wFEg{g}
\approx
\sum_{g\in\ZD}
\rGENg{g}({\qM}_j;\muvec_j)\,\wDg{g},
\qquad
j=1,\ldots,\nsnap .
\end{equation}
If the rule is employed
only to evaluate the projected equilibrium equations, then
\(\rGENg{g}=\rDec{g}\) and \(\ncond=\rD\). If equilibrium and output
evaluation are both hyperreduced using the same rule,
\(\rGENg{g}\) is formed by stacking \(\rDec{g}\) and
\(\rOUTg{g}\), so that \(\ncond=\rD+\nOUT\). A third case arises
when the cubature rule is required only for the
evaluation of the output of interest, i.e.,
\(\rGENg{g}=\rOUTg{g}\) and \(\ncond=\nOUT\). This   occurs when the
latent state is obtained directly from the input parameters,
without the need for solving the projected equilibrium equations.

To apply the ECM, Eq.~\refeq{eq:approxGEN} must first be cast in matrix
form. To this end, we   define
\begin{equation}
\label{eq:Agen_j_definition}
\Agen_j
:=
\begin{bmatrix}
\rGENg{1}({\qM}_j;\muvec_j)^T\\
\rGENg{2}({\qM}_j;\muvec_j)^T\\
\vdots\\
\rGENg{\ngp}({\qM}_j;\muvec_j)^T
\end{bmatrix}
\in\mathbb{R}^{\ngp\times\ncond}, \qquad  j = 1,2 \ldots \nsnap.
\end{equation}
The columns of the preceding matrix contain the values of the \(\ncond\)
scalar-valued integrands at the \(\ngp\) Gauss points of the finite
element mesh, for the $j$-th training parameter.
Concatenating these matrices over the complete training
set gives
\begin{equation}
\label{eq:AgenDEF}
\Agen
:=
\begin{bmatrix}
\Agen_1 &
\Agen_2 &
\cdots &
\Agen_{\nsnap}
\end{bmatrix}
\in
\mathbb{R}^{\ngp\times(\ncond\nsnap)}.
\end{equation}
Next we introduce the vectors of  finite element weights and reduced integration weights
\begin{equation}
\label{eq:weight_vectors_fixed_ECM}
\wFE
=
\begin{bmatrix}
\wFEg{1}\\
\vdots\\
\wFEg{\ngp}
\end{bmatrix},
\qquad
\wD
=
\begin{bmatrix}
\wDg{1}\\
\vdots\\
\wDg{\ngp}
\end{bmatrix}
\end{equation}
respectively. With the preceding definitions, the
integration conditions in Eq.~\refeq{eq:approxGEN} take the compact
matrix form
\begin{equation}
\label{eq:fixed_cubature_sampled_conditions}
\Agen^T\wFE
\approx
\Agen^T\wD.
\end{equation}
 The fixed-weight cubature problem therefore amounts to determining the
sparsest nonnegative vector \(\wD\) that satisfies the above conditions
to the prescribed accuracy.

\subsection{Construction of the fixed-weight ECM rule}
\label{sec:fixed_weight_cubature}

The formulation presented above is common to other  point-selection and weighting
procedures, notably the pioneering cubature method of
An et al.~\cite{an2009optimizing} and the  ECSW
 \cite{farhat2014dimensional}.  Unlike these approaches, however, the ECM does not
address the selection and weighting  directly in terms of the sampled
training matrix \(\Agen\).  Instead, it first constructs a low-dimensional orthonormal basis
\(\Ucub\in\RRn{\ngp}{n_{\mathrm U}}\)  of the column space of $\Agen$, controlled by a  truncation tolerance
 \(\tolSVD\in[0,1)\), i.e.:
\begin{equation}
\label{eq:Ucub_definition}
\dfrac{
\left\|
\Agen-\Ucub\Ucub^T\Agen
\right\|_F
}{
\|\Agen\|_F
}
\leq \tolSVD, \qquad \Ucub^T \Ucub = \ident,
\end{equation}
and then requires\footnote{For simplicity of notation, orthogonality is stated here
with respect to the Euclidean inner product. Other inner products may
be used; in particular, the natural choice for the approximation of
volume integrals is the   inner product induced by
\(\diag{(\wFE)}\), together with its associated norm. This is the
convention employed in the original ECM formulation in Refs.~ \cite{hernandez2017dimensional,hernandez2020multiscale,hernandez2024cecm}.} the reduced
cubature rule to integrate exactly the basis functions represented by
    \(\Ucub\). Together with the requirement that the
selected set \(\ZD\) be as small as possible, this condition leads
formally to the optimization problem
\begin{equation}
\label{eq:ECM_problem_manifold_fixed}
\begin{aligned}
\min_{\wD}\quad
& \|\wD\|_0,
\\
\text{\rm s.t.}\quad
& \Ucub^T\wD=\bcub,
\\
& \wD\geq\bm{0},
\end{aligned}
\end{equation}
where the \(\ell_0\) pseudo-norm \(\|\bullet\|_0\) counts the number
of nonzero entries of its argument, and
\begin{equation}
\label{eq:bcub_definition}
\bcub
:=
\Ucub^T\wFE
\end{equation}
contains the reference integrals of the retained basis functions,
evaluated using the full finite element integration rule. The selected  set of Gauss points is then given
by the support of \(\wD\), i.e.:
\begin{equation}
\label{eq:ZD_support}
\ZD
=
\supp(\wD)
:=
\left\{
g\in\{1,\ldots,\ngp\}
\,:\,
\wDg{g}>0
\right\}.
\end{equation}

A few remarks concerning the preceding optimization problem are in order here.  First, the entries of \(\Agen\) are constructed from the FE training
data as follows. For each FE displacement snapshot
\(\dred(\muvec_j)\), the associated latent coordinates are obtained
through the encoder in Eq.~\refeq{eq:inversemapping}, i.e.,
\(\qM^{train}(\muvec_j)=\Encod(\dred(\muvec_j))\), and are used to
evaluate \(\JTAU\) in Eq.~\refeq{eq:rDec_reduced_operator}. The stresses
entering the same equation, together with any other   quantities
required by the output integrand \(\rOUTg{g}\), are then reevaluated at all
Gauss points from the reconstructed\footnote{For path-dependent constitutive models, the snapshots cannot
be processed independently. Each snapshot must retain its position
within the training trajectory, because the constitutive update requires
the internal variables of the preceding state. Stress reconstruction
must therefore proceed in trajectory order.} displacement vector
\(\dred^{train}(\muvec_j) =
\Decod(\qM^{train}(\muvec_j))\).

 Second,  \(\Ucub\) may   be obtained from a standard,  truncated SVD of \(\Agen\).
For large problems, however, storing and factorizing the complete matrix
may be computationally prohibitive. In the examples shown later, we use the Sequential
Randomized SVD (SRSVD) proposed by the authors in   Ref.~\cite{hernandez2024cecm}, which processes the
column blocks \(\Agen_j\) sequentially and uses randomization to reduce
the cost of the basis construction. Parallel and distributed
randomized SVD algorithms, such as those in Ref.~\cite{ares2025parallel}, provide
an alternative when distributed resources are available.

 Third, to avoid  potential degeneracies for the case $\bcub\approx\bm{0}$, the ECM further imposes the  exact integration of the volume of the domain $|\Omega_0|$:
\begin{equation}
\label{eq:ECM_volume_constraint}
\onevec^T\wD
=
|\Omega_0|,
\end{equation}
where \(\onevec\in\mathbb{R}^{\ngp}\) is the all-ones vector. In order to preserve the orthogonality of $\Ucub$, the
constraint is incorporated by
augmenting \(\Ucub\) with the component of \(\onevec\) orthogonal to the column space of
$\Ucub$; see
Ref.~\cite{hernandez2024cecm} for more details.

Lastly, the \(\ell_0\)-minimization problem in
Eq.~\refeq{eq:ECM_problem_manifold_fixed} is introduced only as a
conceptual statement of the sparsification objective. Because its exact
solution is NP-hard and therefore computationally intractable, the ECM
does not attempt to solve the \(\ell_0\) problem directly. Instead, it
constructs a sparse nonnegative solution of the underdetermined system
\(\Ucub^T\wD=\bcub\). The authors show  in Ref.~\cite{bravo2024subspace} that the ECM invariably achieves a solution with a number of positive weights, denoted hereafter by $\mInit$, equal to the number of columns of $\Ucub$:
\begin{equation}
\label{eq:mInit_fixed_ECM}
\mInit=|\ZD|=\ncol(\Ucub)
\end{equation}
(this follows from Carathéodory's theorem for convex cones
\cite[\S~17]{rockafellar1996convex}, and the fact that the FE weights themselves are positive, and thus \(\bcub=\Ucub^T\wFE\)  lies  in the conical hull of the rows of \(\Ucub\)  ).    Ref.~\cite{bravo2024subspace} also shows  that the solution is generally nonunique:
distinct rules
can be elicited by varying the initial candidate set, thereby allowing
selected regions of the domain to be favored. We do not exploit this
possibility in the examples presented later, but instead adopt the
default choice of considering all Gauss points as candidates.

\section{Cubature with manifold-adaptive weights}
\label{sec:MAWECM}

\subsection{Problem statement and requirements}
\label{sec:motivation_MAW_ECM}
We  introduce next the main novel contribution of the present work: the cubature scheme with manifold-adaptive weights.  The idea is to select a
subset of the Gauss points provided by the preceding fixed-weight rule
and assign them positive, state-dependent weights, in such a way that the resulting cubature rule reproduces, with the same accuracy as the fixed-weight rule, the integration conditions in Eq.~\refeq{eq:approxGEN}, that is:
\begin{equation}
\label{eq:approxGEN2}
\sum_{g\in\ZD}
\rGENg{g}({\qM}^{j};\muvec^{j})\,\wDg{g}
=
\sum_{g\in\Zomega}
\rGENg{g}({\qM}^{j};\muvec^{j})\,
\omegag{g}({\qM}^{j}),
\qquad
j=1,\ldots,\nM,
\end{equation}
where
\begin{equation}
\label{eq:Zomega_subset_ZD}
\Zomega\subseteq\ZD,
\end{equation}
and
\begin{equation}
\label{eq:adaptive_weight_map}
\omegag{g}:\mathbb{R}^{\rD}\rightarrow\mathbb{R}_{+},
\qquad g\in\Zomega.
\end{equation}
Thus, the novelty lies in that  each retained Gauss point \(g\in\Zomega\) is endowed with a positive
scalar weight field over the latent space. This is the sense in which
the method is termed manifold-adaptive weights (MAW) cubature: the
locations of the integration points remain fixed in the reference
configuration $\Omega_0$, whereas their positive weights adapt to the current
state on the nonlinear solution manifold.

%
%

The   reduced
equilibrium equations in this new setting are obtained by replacing the
full finite element rule for internal forces in Eq.~\refeq{eq:manifold_ROM_equilibrium_FE_weights} with the manifold-adaptive rule:

\begin{equation}
\label{eq:MAW_ROM_equilibrium}
\overbrace{
\sum_{g\in\Zomega}
\rDec{g}(\qM;\muvec)\,
\omegag{g}(\qM)
}^{\displaystyle\FintD:\textrm{Internal forces}}
-
\overbrace{
\PhiD(\qM)^T\FextONE\,\muvec
}^{\displaystyle\FextD:\textrm{External forces}}
=
\bm{0}.
\end{equation}
  Because
the cubature weights now depend on \(\qM\), their variation must also
be included in the consistent Newton--Raphson linearization. The resulting tangent operator  is
derived in~\ref{app:MAW_linearization}.

To construct the adaptive weight fields, we first extract from the
input parameter set $\PmuTRAIN$, see Eq.~\refeq{eq:pmutrain}, a   subset
\begin{equation}
\label{eq:PmuM_definition}
\PmuM
=
\left\{\muvec^{j}\right\}_{j=1}^{\nM}
\subseteq
\PmuTRAIN,
\end{equation}
where, for simplicity and with a slight abuse of notation, the retained
parameters have been reindexed from \(1\) to \(\nM\). The corresponding
latent states are taken as the nodes of a graph \(\Glat\), with edges
encoding neighborhood relations on the latent manifold.  The Laplacian of
\(\Glat\) will later be used to measure and penalize distortions of the
weight fields between neighboring states.

The subsampling represented by expression~\refeq{eq:PmuM_definition}   may serve the usual
purpose of avoiding an unnecessarily fine representation of the weight
fields and the associated risk of overfitting. In some problems,
however,  restricting the
state-dependent conditions to \(\PmuM\) does not necessarily discard
the integration content associated with the eliminated states.
Integrand components that must be reproduced independently of the
current latent state are instead collected in an invariant matrix and
enforced for every \(\muvec^{j}\in\PmuM\), as described in the ensuing
subsection.   In the damage benchmark of Section~\ref{sec:DamageProblem}, for
instance, \(\PmuM\) contains only states from the damaged regime, while
the integrand space associated with the undamaged elastic response is
incorporated into the invariant matrix.

With this setting, the  selected integration points $\Zomega$ and the associated state-dependent weights  must satisfy four requirements:

\begin{enumerate}
\item \label{item:MAW_sparsity}
The selected set should be  as small as possible, i.e.,
\(\nomega=|\Zomega|\) is to be minimized.

\item \label{item:MAW_equa}
The adaptive rule must satisfy exactly the state-wise integration
conditions in Eq.~\refeq{eq:approxGEN2} for \(\muvec^{j}\in\PmuM\).

\item \label{item:positive}
The weights must remain nonnegative over the sampled latent domain, and also
preserve the total volume; that is, for a given ${\qM}$
\[
\omegag{g}(\qM) \ge 0,   \hspace{1cm} \sum_{g\in\Zomega}\omegag{g}(\qM)=|\Omega_0|.
\]
The second condition is the local counterpart of Eq.~\refeq{eq:ECM_volume_constraint}.

\item \label{item:smooth}
 Each scalar field
\(\omegag{g}=\omegag{g}(\qM)\), \(g\in\Zomega\), should vary as smoothly
as possible over the manifold.


\end{enumerate}

The preceding requirements are cast as a constrained optimization
problem over the next two subsections. The first subsection develops the
initialization and local constraints associated with
requirements~\ref{item:MAW_equa} and~\ref{item:positive}.
Section~\ref{sec:formulationGENERAL} then introduces the objective
function incorporating requirements~\ref{item:MAW_sparsity}
and~\ref{item:smooth}.

\subsection{Initialization and local integration constraints}
\label{sec:MAW_local_constraints}

As established in expression~\refeq{eq:Zomega_subset_ZD}, the reference
integration rule is not the original FE Gauss-point rule, but the
fixed-weight ECM rule discussed in Section~\ref{sec:fixed_weight_cubature}, characterized by the sparse weight
vector \(\wD\in\mathbb{R}^{\ngp}\) and its support \(\ZD\), of size $\mInit$. Accordingly, we set as design
variables of the optimization problem a  collection of
 weight vectors (one for each sampled state) of size $\mInit$
\[
\left\{
\omegaj{1},
\omegaj{2},
\ldots,
\omegaj{\nM}
\right\},
\]
where
\[
\omegaj{j}\in\mathbb{R}^{\mInit},
\qquad
j=1,\ldots,\nM,
\]
contains the cubature weights associated with the   $j$-th latent state, expressed on the initial candidate support \(\ZD\). We anticipate that the proposed solution algorithm initializes all
such state-dependent weight vectors with the initial weight vector
\begin{equation}
\label{eq:initialW}
\omegaInit
:=
\wDZ
\in\mathbb{R}^{\mInit},
\end{equation}
where \(\wDZ\) denotes the restriction of the sparse fixed-weight vector
\(\wD\) to its support \(\ZD\).

   For each
\(\muvec^{j}\in\PmuM\), let \(\Agen_j\) denote the corresponding local
integrand block introduced in Eq.~\refeq{eq:Agen_j_definition}.   The fixed ECM rule is constructed to integrate exactly those functions
whose discrete representations lie in the subspace spanned by
\(\Ucub\). Since \(\Ucub\) is obtained from the truncated SVD of the
global matrix \(\Agen\) (or a similar low-rank decomposition), the integrand information retained at the
\(j\)-th state is represented by the projected block
\begin{equation}
\label{eq:MAW_projected_local_block}
\widetilde{\Agen}_j
:=
\Ucub\Ucub^{T}\Agen_j
\in
\mathbb{R}^{\ngp\times\ncond},
\qquad
j=1,\ldots,\nM .
\end{equation}
Only in the absence of truncation, i.e., when \(\tolSVD=0\) in
expression~\refeq{eq:Ucub_definition}, does the basis \(\Ucub\) span the full
column space of \(\Agen\), in which case
\(\widetilde{\Agen}_j=\Agen_j\) for every \(j\).

As previously pointed out,  the local projected blocks are complemented  with state-independent
integration conditions. These are collected in the invariant matrix
\begin{equation}
\label{eq:invariants1}
\Uinv
\in
\mathbb{R}^{\ngp\times\ninv}.
\end{equation}
At minimum, \(\Uinv\) contains the all-ones vector
\(\onevec \in \Rn{\ngp}{}\), so that the volume constraint in
Eq.~\refeq{eq:ECM_volume_constraint} is enforced at every
state. To remove linearly dependent conditions, the invariant and projected
integrand blocks are combined and factorized by an SVD:
\begin{equation}
\label{eq:invariants}
\left[
\Uinv,\,
\widetilde{\Agen}_j
\right]
=
\Ubarj{j}\,
\boldsymbol{\Sigma}_j\,
\mathbf{V}_j^{T},
\qquad
j=1,\ldots,\nM ,
\end{equation}
where only the singular vectors associated with nonzero singular values
are retained. The columns of
\(\Ubarj{j}\in\RRn{\ngp}{\ncondQ{j}}\)  therefore form an orthonormal
basis for the combined invariant and state-dependent constraint space. Note that
\begin{equation}
\label{eq:nconstrainsBOUND}
  \ncondQ{j} \le  \ncond + \ninv
\end{equation}
for all $j=1,2,\ldots,\nM$.

Since the fixed-weight ECM rule is supported on \(\ZD\), only the rows
of \(\Ubarj{j}\) associated with the initial candidate points contribute
to the reference integrals. We accordingly define
\begin{equation}
\label{eq:local_constraint_system_definition}
\Uj{j}
:=
\left.\Ubarj{j}\right|_{\ZD}
\in
\mathbb{R}^{\mInit\times\ncondQ{j}},
\end{equation}
and evaluate the corresponding reference integrals using the initial
weight vector defined in Eq.~\ref{eq:initialW}:
\begin{equation}
\label{eq:local_reference_integrals}
\bj{j}
:=
\Uj{j}^{T}\omegaInit
\in
\mathbb{R}^{\ncondQ{j}}.
\end{equation}

The desired local equality constraints then follow directly: at every
sampled latent state, the adaptive weights must satisfy
\begin{equation}
\label{eq:constraintsLOCAL}
\Uj{j}^{T}\omegaj{j}
=
\bj{j},
\qquad
j=1,\ldots,\nM .
\end{equation}

%

%
%

\subsection{Objective function}
\label{sec:formulationGENERAL}

%

To accommodate    requirements ~\ref{item:MAW_sparsity} (minimization of the number of integration points) and~\ref{item:smooth}  (smoothness of the sampled
weight fields), we propose the following objective function:
\begin{equation}
\label{eq:MAW_objective}
\mathcal{J}
\defeq
\left\|
\sum_{j=1}^{\nM}
\omegaj{j}
\right\|_0
+
 \Reg(\Wad,\Glat),
\end{equation}
where
\begin{equation}
\label{eq:adaptive_weight_matrix}
\Wad
\defeq
\begin{bmatrix}
\omegaj{1} &
\omegaj{2} &
\cdots &
\omegaj{\nM}
\end{bmatrix}
\in
\mathbb{R}^{\mInit\times\nM}.
\end{equation}
The first term on the right-hand side of
Eq.~\eqref{eq:MAW_objective} is precisely the number of integration
points $\nomega$ (this follows from the fact that   nonnegativity prevents cancellations in the sum of the local
weight vectors).  The second term,
\(\Reg=\Reg(\Wad,\Glat)\), regularizes the sampled weight fields,
represented by the rows of \(\Wad\). It will be defined as a quadratic
functional involving, among other contributions, the graph Laplacian
associated with \(\Glat\), thereby penalizing irregular variations of
the weights over the sampled manifold. At this stage,
\(\Reg(\Wad,\Glat)\) is introduced only in abstract form; its specific
definition, together with the associated solution strategy, is
presented later.

\subsection{Optimization problem and solution strategy}
\label{sec:MAW_solution_strategy}


With all the ingredients introduced in the preceding subsections, the
manifold-adaptive weights ECM problem can now be formulated as the following
constrained optimization problem: find the  weight vectors
\(\{\omegaj{j}\}_{j=1}^{\nM}\), collected in \(\Wad\), as the solution
of
\begin{equation}
\begin{aligned}
\min_{\omegaj{1},\ldots,\omegaj{\nM}}
\quad &
\left\|
\sum_{j=1}^{\nM}
\omegaj{j}
\right\|_0
+
 \Reg(\Wad,\Glat)
\\
\text{subject to}
\quad &
\Uj{j}^{T}\omegaj{j}
=
\bj{j},
\qquad
j=1,\ldots,\nM,
\\
&
\omegaj{j}\ge0,
\qquad
j=1,\ldots,\nM .
\end{aligned}
\label{eq:MAW_global_l0}
\end{equation}

The solution determines simultaneously the state-dependent weights and
the final common support \(\Zomega\), identified by the nonzero rows of
\(\Wad\).  As discussed in
Section~\ref{sec:mawECMintro}, we   adopt a greedy pruning strategy. Starting from the feasible
fixed-weight ECM rule, integration points are removed successively while
the remaining manifold-dependent weights are redistributed so as to
preserve the local equality and positivity constraints and control the
regularity of the resulting weight fields. The resulting pruning
strategy is developed in the next section.

\section{Greedy pruning strategy}
\label{sec:Greedy}

To describe the proposed   greedy pruning procedure, we proceed
hierarchically from its most elementary operation to the final global
algorithm. We begin by isolating the fundamental building block of the
construction, namely a single tentative pruning step. Assume that one
integration point is selected for removal from the current active
support. The central question then becomes whether the remaining
adaptive weights can be redistributed so as to preserve the local
exactness constraints over the sampled manifold. This leads to the
\emph{local redistribution problem} introduced in
Section~\ref{sec:redistributionPR}, and  whose solution procedure is
summarized in Algorithm~\ref{alg:prune_step_optionB}.  Subsequently,
Section~\ref{sec:candidatesweepbranch} introduces
Algorithm~\ref{alg:graph_candidate_sweep}, whose purpose is to partially
mitigate the path dependence intrinsic to the greedy pruning strategy by
exploring several feasible candidate eliminations before committing to
one of them. Finally, the complete greedy pruning loop is described in
Section~\ref{sec:globalpruning} through
Algorithm~\ref{alg:MAW_global_pruning}, which repeatedly applies the
candidate-sweep procedure.  For convenience, the complete offline workflow is summarized in Section \ref{sec:summary},
Box~\ref{box:MAW_ECM_summary}. Readers primarily interested in the
practical implementation may jump directly to this summary and navigate
from there to the corresponding algorithms.

 \subsection{Redistribution problem}
\label{sec:redistributionPR}

%

Assume that, at a given pruning stage, the current active support is $\Iact \subset \{1,\ldots,\mInit\}$.
Here, active refers to the integration points that have not yet been
removed from the initial candidate support.   For each sampled latent state \(\qMj{j}\), the current adaptive
weights are denoted by $\wold{j}\in\mathbb{R}^{|\Iact|}$.  By construction, they satisfy the local linear equality constraints in Eq.~\refeq{eq:constraintsLOCAL}:
\begin{equation}
\UTact{j}\wold{j}
=
\bj{j},
\qquad
j=1,\ldots,\nM,
\end{equation}
where \({\Uj{j}}_{\Iact}\) denotes the restriction of \(\Uj{j}\) to
the active support.  A pruning trial begins by selecting one candidate integration point
(the selection strategy will be discussed later), $
p\in\Iact$, and attempts to remove it. The remaining active support is $\Irem =
\Iact\setminus\{p\}$. Equivalently, in the local numbering induced by the current active
support \(\Iact\), $\IremLOC = \{1,\ldots,|\Iact|\}\setminus\{\loc(p,\Iact)\}$.
Here, \(\loc(p,\Iact)\) denotes the position occupied by the global
index \(p\) within the ordered set \(\Iact\). The question is whether, for every sampled latent state, one can construct redistributed weights $\wnew{j}\in\mathbb{R}^{|\Irem|}$
such that
\begin{equation}
{\Uj{j}}^T_{\Irem}\wnew{j}
=
\bj{j},
\qquad
j=1,\ldots,\nM.
\end{equation}
If such weights exist and satisfy the required admissibility conditions,
the candidate removal is accepted.



\subsubsection{Regularized optimization problem}
\label{sec:regoptprobl}

 The  weights we seek are the values of $|\Irem|$ scalar fields defined at $\nM$ samples of the latent space:
\begin{equation}
\Wnew
=
\begin{bmatrix}
\label{eq:wnew}
\big| & \big| &        & \big| \\
\wnew{1} & \wnew{2} & \cdots & \wnew{\nM} \\
\big| & \big| &        & \big|
\end{bmatrix}
\in
\mathbb{R}^{|\Irem|\times\nM},
\end{equation}
(the column
\(\wnew{j}\in\mathbb{R}^{|\IremLOC|}\) contains the   weights
associated with the sampled latent state \(\qMj{j}\)). Similarly, the weights from the previous pruning stage, restricted to
the same remaining support, are denoted by
\[
\Woldrem
=
\begin{bmatrix}
\big| & \big| &        & \big| \\
\woldrem{1} & \woldrem{2} & \cdots & \woldrem{\nM} \\
\big| & \big| &        & \big|
\end{bmatrix}
\in
\mathbb{R}^{|\Irem|\times\nM},
\]
where $\woldrem{j} = \woldremDEF{j}$.  We propose to determine the redistributed weights $\Wnew$ as the solution of the
following constrained quadratic optimization problem:
\begin{equation}
 \label{eq:quadraticGENERAL}
\begin{aligned}
\min_{\Wnew}\quad
&
\frac12
\|
\Wnew-\Woldrem
\|_F^2
+
\frac{\alphaG}{2}
\operatorname{tr}
\left[
\left(\Wnew-\Woldrem\right)
\KGs
\left(\Wnew-\Woldrem\right)^T
\right]
\\
\text{subject to}\quad
&
\UremDEFt{j}\wnew{j}
=
\bj{j},
\qquad
j=1,\ldots,\nM,
\\
&
\wnew{j}\ge0,
\qquad
j=1,\ldots,\nM,
\end{aligned}
\end{equation}
where \(\operatorname{tr}(\cdot)\) denotes the matrix trace operator. 

The first contribution in the objective function is a least-change term written in Frobenius norm
form; notice it does not introduce any coupling between different sampled latent states, since
\[
\|
\Wnew-\Woldrem
\|_F^2
=
\sum_{j=1}^{\nM}
\|
\wnew{j}-\woldrem{j}
\|_2^2.
\]
Nevertheless, despite this decoupling, the local least-change
redistribution still possesses an important implicit smoothing effect.
Indeed, the redistribution is computed as the smallest admissible
departure from the previously accepted weight fields. Since the initial
adaptive weights are constant over the latent manifold, successive
least-change redistributions tend naturally to preserve this smooth
character and introduce only gradual local modifications.
The smooth evolution described above may eventually be altered by the
appearance of negative weights. In practice, this introduces localized irregularities, since solution algorithms    clamp these negative weights  to zero.   The second contribution in the objective function of
Problem~\eqref{eq:quadraticGENERAL} is a graph-based, quadratic
regularization term that naturally couples neighboring latent states and
penalizes these localized irregularities, with \(\alphaG\ge0\)
controlling the strength of the regularization.

The notion of neighborhood entering this
regularization is introduced through the already mentioned   graph structure
defined over the sampled latent manifold. More precisely,  let
\[
\Glat
=
\big(\{\qMj{j}\}_{j=1}^{\nM},\EdgesM\big)
\]
denote a graph whose nodes are the sampled latent states and whose edge
set \(\EdgesM\) encodes neighborhood relations over the sampled
manifold. The matrix
$
\KGs\in\mathbb{R}^{\nM\times\nM}
$
denotes the associated scalar graph operator defined on this graph,
typically chosen as a graph Laplacian associated with the connectivity
of the sampled latent states; see, e.g.,
Belkin and Niyogi~\cite{belkin2003laplacian}.




  The edge set \(\EdgesM\) may be defined, for
instance, from nearest-neighbor relations in latent space, as commonly
done in manifold-learning methods based on graph operators \cite{belkin2003laplacian}.
In the benchmarks
considered later, the sampled latent states admit particularly simple
graph constructions. In the metamaterial-cell problem of Section \ref{sec:metamaterial}, they lie on a
structured two-dimensional Cartesian grid. The resulting connectivity
defines an auxiliary quadrilateral mesh in latent space, and
\(\KGs\) is obtained from the standard finite element discretization
of an isotropic Laplacian using bilinear shape functions. In the
damage problem of Section \ref{sec:DamageProblem}, the adaptive weights depend on a single normalized
nonlinear coordinate. The ordered samples therefore define a
one-dimensional mesh, and \(\KGs\) is constructed from the
corresponding finite element discretization of the one-dimensional
Laplacian using linear shape functions.

 \subsubsection{Solution strategy}
\label{sec:stagegraph}


\begin{algorithm}[!ht]
\footnotesize
\SetCommentSty{tcpstyleFOOT}

 \DontPrintSemicolon

\KwData{
Current active support \(\Iact\); candidate point \(p\in\Iact\);
current weights \(\Wold\); local equality systems restricted to the
current active support,
\(\{(\Uactj{j},\bj{j})\}_{j=1}^{\nM}\),
graph operator \(\KGs\); regularization parameter \(\alphaG\ge0\).
}

\KwResult{Acceptance flag \(\SUCCESS\), remaining support \(\Irem\),
redistributed weights \(\Wnew\), and redistribution energy \(\Sdir \).}

\(\ploc\leftarrow\operatorname{loc}(p,\Iact)\)
\tcp{Local index of candidate point \(p\) in $\Iact$ \label{alg1:local}}
\(\Woldrem\leftarrow\Wold;\;\; (\Woldrem)_{\ploc,:}\leftarrow0\)
\tcp{Remove candidate point \(p\) by zeroing its weight field}

\(\Dj{j}\leftarrow\{\ploc\},\quad j=1,\ldots,\nM\)
\tcp{Initialize local active sets by tentatively clamping the candidate point \label{alg1:Dini}}

\(\SUCCESS\leftarrow\TRUE\) ; $\;$  \(\ENFORCE\leftarrow\FALSE\); $\;$ $\Sdir \leftarrow 0$ \tcp{Flags for feasibility/activation positivity. Energy initializ.\label{alg1:flagENF}}

\(\kact\leftarrow0\)
\tcp{Active-set iteration counter \label{alg1:counter}}

\While{active sets change}{

    \For{\(j=1,\ldots,\nM\)   \label{alg1:loopstates}}{

       \(\Fj{j}\leftarrow\{1,\ldots,|\Iact|\}\setminus\Dj{j}\)
\tcp{Free local rows at latent state \(j\)
\label{alg1:prune_free}}
 \(
(\Psij{j},\Sigmaj{j},\Vj{j})
\leftarrow
\SVD(({\Uactj{j}}_{\Fj{j}})^T)
\)
\tcp{Full SVD, see Eq. \refeq{eq:fullSVD}
\label{alg1:prune_svd}}

\(
\rj
\leftarrow
\operatorname{rank}(\Sigmaj{j})
\)
\tcp{Numerical rank associated with   nonzero singular values
\label{alg1:svdRANK}}

       \If{\(\rj<\ncondQ{j}\)}{
    \(\SUCCESS\leftarrow\FALSE\)
    \tcp{Full-row-rank acceptance condition fails  \label{alg1:nosuccess}}

         Reject candidate and set
\(\Irem\leftarrow\Iact,\ \Wnew\leftarrow\Wold\); $\;\;\;$ \textbf{stop}

        }

 Compute particular solution $\wpfull{j}$
\tcp{(see Eq.~\ref{eq:particularSOLUTION} )
\label{alg1:prune_particular}}

Compute null-space basis \(\Nj{j}\)
\tcp{ ( last \( |\Fj{j}|-\rj \) columns of \(\Vj{j}\)) )
\label{alg1:prune_null}}
        \If{\(\wpfull{j}\) contains negative entries}{
            \(\ENFORCE\leftarrow\TRUE\)
            \tcp{Positivity enforcement required  \label{alg1:enforceREQUIRE}}
        }

    }

\eIf{\(\ENFORCE=\FALSE\)}{

    \(\SUCCESS\leftarrow\TRUE\)
    \tcp{Successful pruning step, no need for regularization  \label{alg1:success3}}

    \(\Irem\leftarrow\Iact\setminus\{p\}\), assemble \(\Wnew\) from the local vectors \(\wpfull{j}\); $\;\;$
    restrict \(\Wnew\) to rows \(\{1,\ldots,|\Iact|\}\setminus\{\ploc\}\); $\;\;$ \textbf{stop}}

{

    \If{\(\kact=0\)}{

       \(\zold \leftarrow \operatorname{vec}(\Woldrem)\)
\tcp{Vectorize  by column stacking
\label{alg1:columns}}

\(
\Hquad
=
\ident_{|\Iact|\nM}
+
\alphaG
(\KGs\otimes\ident_{|\Iact|})
\)
\tcp{Graph-regularized operator (Eq. \ref{eq:Hoperator})
\label{alg1:prune_H}}

\(\gquad \leftarrow \Hquad \zold\)
\tcp{Reference state for the graph-regularized redistribution
\label{alg1:gCOMP}}

    }

}

%
%

Assemble global variables \(\zpart\) and \(\Nmat\)
\tcp{From \(\{\wpfull{j}\}_{j=1}^{\nM}\) and \(\{\Nj{j}\}_{j=1}^{\nM}\)
\label{alg1:prune_globalparam}}

    \(
    \left(\Nmat^T\Hquad\Nmat\right)\yred
    =
    \Nmat^T
    \left(\gquad-\Hquad\zpart\right)
    \)
    \tcp{Solve   for $\yred$ (see Eq. \ref{eq:solveY})
    \label{alg1:prune_reducedsolve}}

    \(
    \znew
    \leftarrow
    \zpart+\Nmat\yred
    \)
    \tcp{Reconstruct redistributed weights  \label{alg1:reconsZ}}

    Reshape \(\znew\) into \(\Wnew\) and   detect negative free entries of \(\Wnew\)

    \eIf{no negative free entries are detected}{

        $\SUCCESS\leftarrow\TRUE$
        \tcp{Successful positivity-enforced redistribution  \label{alg1:sss3333}}
\(\Sdir\leftarrow
\operatorname{tr}\!\left((\Wnew-\Woldrem)\KG(\Wnew-\Woldrem)^T\right)\)
\tcp{Graph Dirichlet energy \label{alg1:energy}}

    \(\Irem\leftarrow\Iact\setminus\{p\}\); $\;$  restrict \(\Wnew\) to   rows
\(\{1,\ldots,|\Iact|\}\setminus\{\ploc\}\); $\;\;$

\textbf{stop}

    }{

        Update local active sets \(\Dj{j}\)
        \tcp{Clamp negative weights to zero
        \label{alg1:prune_active}}

          \(\kact\leftarrow\kact+1\)
\tcp{Update active-set iteration counter  \label{alg1:counterU}}

    }
}

\caption{\textsc{PruneStep}: Active-set strategy for solving the quadratic graph-regularized
redistribution problem \ref{eq:quadraticGENERAL},  associated with the tentative pruning of one
integration point. This local redistribution procedure is
repeatedly invoked within
Algorithm~\ref{alg:graph_candidate_sweep} to evaluate different feasible
candidate eliminations.
\label{alg:prune_step_optionB}}
\end{algorithm}

Problem~\eqref{eq:quadraticGENERAL} is a standard convex quadratic
program with linear equality and positivity constraints, and could in
principle be solved using general-purpose quadratic optimization
algorithms; see, e.g., Refs.~\cite{boyd2004convex,Nocedal1999}.
However, off-the-shelf quadratic optimization procedures do not exploit
the particular structure of the present pruning process and therefore
eventually become inefficient. In particular, enforcing the
graph-regularization term (the one involving the graph matrix
\(\KGs\)) at every pruning step is both unnecessary and
computationally expensive. As discussed earlier, the graph-regularization
term becomes truly necessary only when positivity is violated.

The proposed solution strategy, summarized in Algorithm~\ref{alg:prune_step_optionB}, does exploit this  feature of the pruning
process, and distinguish two distinct    redistribution regimes (controlled by
the Boolean flag \(\ENFORCE\) in Line~\ref{alg1:flagENF}). The uncoupled regime consists in the solution of individual least-norm problems:
\begin{equation}
\label{eq:smallproblems}
\begin{aligned}
\min_{\dwj{j}}
\quad &
\frac12
\|
\dwj{j}
\|_2^2
\\
\text{subject to}
\quad &
\UremDEFt{j}\dwj{j}
=
\bj{j} - \UremDEFt{j}\woldrem{j},
\qquad
j=1,2,\ldots,\nM,
\end{aligned}
\end{equation}
(this readily follows from setting $\alphaG = 0$ in  Problem~\eqref{eq:quadraticGENERAL}). The
positivity-enforcement regime, on the other hand,  is  based on  a null-space active-set strategy in
which the local equality constraints are eliminated exactly through
local null-space parametrizations, and positivity  is enforced
through latent-state-dependent active sets.




The first steps of the algorithm are common to both redistribution
regimes and consist of a loop over the \(\nM\) sampled latent states
(Line~\ref{alg1:loopstates}). The key operation performed at each state $j$ in this loop
is the following  SVD:
(Line~\ref{alg1:prune_svd})
\begin{equation}
\label{eq:fullSVD}
({\Uactj{j}}_{\Fj{j}})^T
=
\Psij{j}\Sigmaj{j}(\Vj{j})^T,
\qquad
\Sigmaj{j}
=
\diag(\sigma_1^{j},\ldots,\sigma_{\rj}^{j},0,\ldots,0) \in
\mathbb{R}^{\ncondQ{j}\times|\Fj{j}|}.
\end{equation}
where
\(
{\Uactj{j}}_{\Fj{j}}
\in
\mathbb{R}^{|\Fj{j}|\times\ncondQ{j}}
\)
denotes the local equality system restricted to the current free rows
\(\Fj{j}\subseteq\{1,\ldots,|\Iact|\}\), defined in Line~\ref{alg1:prune_free} as the complementary set of the active rows $\Dj{j}$. Initially, \(\Dj{j}\) contains only the
local index \(\ploc\) associated with the tentatively removed
integration point (Line~\ref{alg1:Dini}), although additional rows may
subsequently be added during the positivity-enforcement iterations
(Line~\ref{alg1:prune_active}).

The SVD outputs (left and right singular vectors,$\Psij{j}$ and $\Vj{j}$, respectively, and singular values $\Sigmaj{j}$  ) play a threefold role within the algorithm: they are
used to determine the numerical rank of the restricted local system
(as the number of nonzero
singular values in \(\Sigmaj{j}\), Line~\ref{alg1:svdRANK}), construct a particular feasible
redistribution (Line~\ref{alg1:prune_particular}), and identify a basis
for the associated null space
(Line~\ref{alg1:prune_null}).  The full-row-rank condition ultimately determines whether the
tentative elimination of a candidate point can be accepted.
A violation at any sampled latent state
(Line~\ref{alg1:nosuccess}) immediately causes the tentative
elimination to be rejected. It should be noted that such a violation may
occur for two different reasons: either the restricted local system is
numerically rank deficient, or the active-set updates have clamped so
many weights that the number of remaining free design variables  $|\Fj{j}|$  becomes
smaller than the number of local equality constraints \(\ncondQ{j}\) ($j=1,2 \ldots \nM$).  Numerical experience indicates that the latter is the most frequent cause of rejection.  This observation also identifies the natural stopping limit of
the pruning process:  pruning cannot proceed below the maximum number of independent local
constraints, which, according to \refeq{eq:nconstrainsBOUND},  is given by
\begin{equation}
\label{eq:ncondmin}
\mLower =  \ncond + \ninv.
\end{equation}

%



The particular feasible redistribution computed in
Line~\ref{alg1:prune_particular} is obtained in an incremental
sense. Denoting by \(\dwFj{j}\in\mathbb{R}^{|\Fj{j}|}\) the redistribution
correction restricted to the current free rows, the local equality
constraints become
\[
({\Uactj{j}}_{\Fj{j}})^T \dwFj{j}
=
\btildej{j},
\]
where
\[
\btildej{j}
=
\bj{j}
-
({\Uactj{j}}_{\Fj{j}})^T
({\Woldrem}_{\Fj{j},j}),
\]
The minimum-norm correction satisfying this system is then obtained
through the Moore--Penrose pseudoinverse constructed from the SVD
outputs associated with the nonzero singular values,
\[
\dwFj{j}
=
\Vj{j}_{:,1:\rj}
\diag\!\left(
\frac{1}{\sigma_1^{j}},\ldots,
\frac{1}{\sigma_{\rj}^{j}}
\right)
(\Psij{j}_{:,1:\rj})^T
\,\btildej{j}.
\]
The updated free weights are subsequently reconstructed as
\begin{equation}
\label{eq:particularSOLUTION}
\wFj{j}
=
({\Woldrem})_{\Fj{j},j}
+
\dwFj{j}.
\end{equation}
At the first active-set iteration (\(\kact=0\)), the correction
\(\dwFj{j}\) coincides with the solution of the incremental
least-change problem~\refeq{eq:smallproblems}. Consequently, if the
vectors \(\wFj{j}\) remain nonnegative for all
\(j=1,\ldots,\nM\), the tentative pruning step already provides an
admissible redistribution and the algorithm terminates without
activating the graph-regularized stage
(Line~\ref{alg1:success3}). Otherwise, the same local construction is still used as the affine
origin for the constrained redistribution. For this purpose, each
\(\wFj{j}\) is lifted back to the full local ambient dimension by
inserting zeros on the active rows, leading to the vector
\(\wpfull{j}\in\mathbb{R}^{|\Iact|}\)   in Line~\ref{alg1:prune_particular}. On the other hand, the global vector \(\zpart\) in
Line~\ref{alg1:prune_globalparam} is obtained by stacking these lifted
particular solutions over all sampled latent states.

%
The null-space basis computed in
Line~\ref{alg1:prune_null}  is directly obtained from the matrix of right singular vectors  (specifically, the last  $|\Fj{j}| -\rj$ columns of  \(\Vj{j}\)).
 Then  it is subsequently
lifted back to the full local ambient dimension, producing  the basis matrix
$\Nj{j}\in\mathbb{R}^{|\Iact|\times(|\Fj{j}|-\rj)}$, whose columns span all admissible local homogeneous redistributions
compatible with the current clamping pattern. These local lifted bases  are then assembled into the sparse
block-diagonal matrix \(\Nmat\).

Using the global particular redistribution \(\zpart\) together with the
global null-space basis \(\Nmat\), every admissible redistribution
satisfying both the equality constraints and the current active-set
pattern can be parametrized as
\begin{equation}
\label{eq:nullp}
\znew
=
\zpart+\Nmat\yred,
\end{equation}
where \(\yred\) collects the reduced null-space coordinates associated
with all sampled latent states, whereas $\znew=\operatorname{vec}(\Wnew)$.

 Using standard
vectorization identities for Frobenius norms and trace operators, the
matrix-valued objective function in Problem~\eqref{eq:quadraticGENERAL} can then be expressed as an equivalent quadratic
functional in terms of \(\znew\). After expanding the resulting
expression and discarding terms independent of the optimization
variables, we obtain the reduced quadratic form
\begin{equation}
\label{eq:quadar}
\frac12\znew^T\Hquad\znew-\gquad^T\znew.
\end{equation}
with
\begin{equation}
 \label{eq:Hoperator}
\Hquad
=
\ident_{|\Iact|\nM}
+
\alphaG
(\KGs\otimes\ident_{|\Iact|}),
\qquad
\gquad
=
\Hquad\zold,
\end{equation}
(\(\Hquad\in\mathbb{R}^{(|\Iact|\nM)\times(|\Iact|\nM)}\) and
\(\gquad\in\mathbb{R}^{|\Iact|\nM}\) are introduced in
Line~\ref{alg1:prune_H}). The Kronecker-product
structure arises naturally from the vectorization of the trace term,
since the same graph operator \(\KGs\) acts independently on each one
of the \( |\Iact| \) adaptive weight fields.

Substituting the null-space parametrization~\ref{eq:nullp} into the quadratic expression \refeq{eq:quadar} yields a reduced quadratic problem in the
unknown reduced coordinates \(\yred\). Imposing stationarity with
respect to \(\yred\) in this problem  leads to the linear system
\begin{equation}
\label{eq:solveY}
\left(\Nmat^T\Hquad\Nmat\right)\yred
=
\Nmat^T(\gquad-\Hquad\zpart),
\end{equation}
which corresponds precisely to the equation solved in
Line~\ref{alg1:prune_reducedsolve}.

After solving this system,   \(\znew\) is reconstructed from
the reduced coordinates using Eq. \refeq{eq:nullp} and reshaped into   matrix \(\Wnew\). The
algorithm then inspects  \(\Wnew\).  If
no negative free entry is found the pruning step is accepted, and the final support is set to
\(\Irem=\Iact\setminus\{p\}\). The matrix \(\Wnew\) is then restricted
to the locally retained rows
\(\{1,\ldots,|\Iact|\}\setminus\{\ploc\}\), thereby returning the
redistributed rule on the pruned support
(Line~\ref{alg1:sss3333}).  Besides the redistributed weights themselves, the algorithm also returns the graph Dirichlet
energy of the admissible redistribution
(Line~\ref{alg1:energy}). This quantity is subsequently used in the
outer graph-based pruning stage to rank the feasible candidate
eliminations.

If, on the contrary, negative free entries are detected, the
corresponding rows are added to the appropriate local active sets
\(\Dj{j}\) and clamped to zero
(Line~\ref{alg1:prune_active}). The active-set counter is then updated
(Line~\ref{alg1:counterU}), and the local SVD parametrizations are
recomputed with the enlarged clamped sets. This process is repeated
until either no new negative entries are detected or the candidate
becomes infeasible because the local rank condition fails.

\begin{remark}
\label{remark:positive}
 Numerical experiments shown in Sections \ref{sec:infl_fw_ini} and \ref{sec:iniii} indicate that most  of weight fields
  are removed without activating the graph-coupled
positivity-enforcement stage, which typically becomes necessary only
during the final pruning steps, close to the final sparse solution. This favorable property may be explained by the presence of
the volume-preservation constraint within the local equality systems
(see Eq.~\ref{eq:invariants1}). Indeed, after removing one integration point and its associated weight
field, the remaining weights fields tend naturally to increase, on average,  in order to
compensate for the missing contribution, thereby favoring positivity. This compensation
mechanism may be better grasped through the analogy of a statically
redundant cable-supported structure of total self-weight \(T\). The  weights play the role of cable
tensions, the latent states represent different loading configurations,
and the equality constraints enforce equilibrium for each configuration.
Thus,
removing one integration point and its associated weight field is
analogous to cutting one cable: the remaining cable tensions must then
increase on average over all configurations  in order to equilibrate the total self-weight $T$. As redundancy
decreases, however, some cables may no longer be able to participate in
the equilibrium rearrangement for certain loading configurations and
simply go slack. In the present setting, this corresponds to weights
attempting to become negative and therefore being clamped to zero (in Line~\ref{alg1:prune_active} of Algorithm \ref{alg:prune_step_optionB}).
\end{remark}

%
%

\subsection{Candidate sweep and branch selection}
\label{sec:candidatesweepbranch}

\begin{algorithm}[!ht]
\footnotesize
\SetCommentSty{tcpstyleFOOT}
\DontPrintSemicolon

\KwData{
Current active support \(\Iact\); current weights \(\Wold\);
local equality systems restricted to the current active support,
\(\{(\Uactj{j},\bj{j})\}_{j=1}^{\nM}\) (where
\(\Uactj{j} =  (\Uj{j})_{\Iact}\in\mathbb{R}^{|\Iact|\times\ncondQ{j}}\) and
\(\bj{j}\in\mathbb{R}^{\ncondQ{j}}\));
graph operator \(\KGs\); regularization parameter \(\alphaG\ge0\);
maximum number of feasible pruning trials \(\ntry\).
}

\KwResult{Acceptance flag \(\SUCCESS\), selected candidate \(p^\star\),
remaining support \(\Irem\) and  redistributed weights \(\Wnew\) }

\(\SUCCESS\leftarrow\FALSE\)

\For{\(\ploc=1,\ldots,|\Iact|\)}{
\(
\bar\omega_{\ploc}
\leftarrow
\frac{1}{\nM}
\sum_{j=1}^{\nM}
(\Wold)_{\ploc,j}
\)
\tcp{Average weight carried by local candidate \(\ploc\) over the latent manifold
\label{alg2:average}}
}

Sort \(\Iact\) in increasing order of \(\bar\omega_p\), and denote the resulting ordered list by \(\Iord\)

Initialize empty lists \(\Jcand,\Jrem,\Qw,\Rsdir\); $\;\;$ \(k\leftarrow1\)\;

\While{\(k\le|\Iord|\) \(\AND\) \(|\Jcand|<\ntry\)}{

\(p\leftarrow(\Iord)_k\)\;

\(({\SUCCESS}_p,\Irem^{p},\Wnew^{p},\Sdir^{p})
\leftarrow
\textsc{PruneStep}
(\Iact,p,\Wold,\{(\Uactj{j},\bj{j})\}_{j=1}^{\nM},\KGs,\alphaG)\) \tcp{Algorithm \ref{alg:prune_step_optionB}\label{alg2:invoke}}
\If{\({\SUCCESS}_p=\TRUE\)}{
    Append \(p\), \(\Irem^{p}\), \(\Wnew^{p}\), and \(\Sdir^{p}\) to
    \(\Jcand,\Jrem,\Qw,\Rsdir\), respectively   \label{alg2:success} \;

    \If{\(\Sdir^{p}=0\)}{
        \textbf{break}
        \tcp{No graph redistribution was needed; no smoother feasible trial is possible \label{alg2:noneedrank}}
    }
}

\(k\leftarrow k+1\)  \label{alg2:updatek} \;
}

\eIf{\(\Jcand=\emptyset\)}{
\(\SUCCESS\leftarrow\FALSE\) \tcp{No feasible pruning candidate was found \label{alg2:noluck}}
\(\Irem\leftarrow\Iact\), \(\Wnew\leftarrow\Wold\) \tcp{Preserve current support and weights}
}{
\(a^\star\leftarrow\arg\min_a (\Rsdir)_a\) \tcp{Choose feasible pruning trial with smallest graph Dirichlet energy \label{alg2:rank}}

\(p^\star\leftarrow(\Jcand)_{a^\star}\);  $\;\;$ \(\Irem\leftarrow(\Jrem)_{a^\star}\);  $\;\;$   \(\Wnew\leftarrow(\Qw)_{a^\star}\)
\(\SUCCESS\leftarrow\TRUE\) \label{alg2:update}
}

\caption{\textsc{PruneSweep}: exploration and ranking of feasible
pruning candidates through repeated calls to
Algorithm~\ref{alg:prune_step_optionB}.  This
candidate-sweep procedure is in turn repeatedly invoked within
Algorithm~\ref{alg:MAW_global_pruning} during the global
pruning process.}
\label{alg:graph_candidate_sweep}
\end{algorithm}

Algorithm~\ref{alg:prune_step_optionB}  is inherently
path-dependent: once a candidate elimination is accepted, the subsequent
pruning sequence is conditioned by that choice, while the alternative
supports associated with other feasible eliminations are no longer pursued.
To mitigate this dependence on a single greedy decision, we refine the
search through Algorithm~\ref{alg:graph_candidate_sweep}, which repeatedly
invokes Algorithm~\ref{alg:prune_step_optionB}
(Line~\ref{alg2:invoke}) for different candidate points, generating up to
a user-prescribed number \(\ntry\) of feasible pruning branches before
committing to one elimination.

The sweep itself is driven by two successive indicators. First, the
candidate points are prioritized according to the average value of their
weights over the sampled latent manifold (Line~\ref{alg2:average}). Candidates carrying smaller
average weight are tested first,  since they are expected to play a less significant role in
satisfying the local exactness constraints.

The second indicator driving the sweep is the graph Dirichlet energy
computed in Line~\ref{alg1:energy} of
Algorithm~\ref{alg:prune_step_optionB}. Among all feasible pruning trials
stored in Line~\ref{alg2:success}, the sweep selects the candidate with
smallest  energy (Line~\ref{alg2:rank}).  It should be noted that this secondary ranking
becomes relevant only when the positivity-enforcement stage of
Algorithm~\ref{alg:prune_step_optionB} has actually modified the
redistribution, since otherwise the Dirichlet energy contribution does not
participate in the objective function of
Problem~\eqref{eq:quadraticGENERAL}. In such cases,
Algorithm~\ref{alg:prune_step_optionB} simply returns zero Dirichlet
energy, and consequently Algorithm~\ref{alg:graph_candidate_sweep} explores no further
branches (Line~\ref{alg2:noneedrank}).

\subsection{Global pruning loop}
\label{sec:globalpruning}
\begin{algorithm}[!ht]
\footnotesize
\SetCommentSty{tcpstyleFOOT}
\DontPrintSemicolon

\KwData{
Initial fixed ECM support
\(\ZD=\{z_1,\ldots,z_{\mInit}\}\subset\{1,\ldots,\ngp\}\);
initial fixed ECM weights
\(\omegaInit\in\mathbb{R}^{\mInit}\), defined on \(\ZD\);
local equality systems
\(\{(\Uj{j},\bj{j})\}_{j=1}^{\nM}\), with
\(\Uj{j}\in\mathbb{R}^{\mInit\times\ncondQ{j}}\) defined on \(\ZD\)
and \(\bj{j}\in\mathbb{R}^{\ncondQ{j}}\);
graph operator \(\KGs\in\mathbb{R}^{\nM\times\nM}\);
regularization parameter \(\alphaG\ge0\);
maximum number of feasible pruning trials \(\ntry\).
}

\KwResult{
Subset of retained Gauss points
\(\Zomega\subseteq\ZD\) and associated   weights
\(\Wadc\in\mathbb{R}^{|\Zomega|\times\nM}\) (the \(i\)-th row of \(\Wadc\)
contains samples over the latent space of the adaptive weight field
associated with the \(i\)-th retained Gauss point).
}

\(\Iact\leftarrow\{1,\ldots,\mInit\}\)
\tcp{Active local indices with respect to the fixed support \(\ZD\) \label{alg3:ini1}}

\(\Wold\leftarrow[\omegaInit,\ldots,\omegaInit]\in\mathbb{R}^{\mInit\times\nM}\)
\tcp{Initial adaptive weights, still represented on the full support \(\ZD\)  \label{alg3:ini2}}

\(\SUCCESS\leftarrow\TRUE\)\;

\While{\(\SUCCESS=\TRUE\)}{

\(\Uactj{j}\leftarrow(\Uj{j})_{\Iact}\in\mathbb{R}^{|\Iact|\times\ncondQ{j}}\),
\(j=1,\ldots,\nM\)
\tcp{Restrict local systems to the current active support \label{alg3:U}}

\((\SUCCESS,\Irem,\Wnew,\Sdir)
\leftarrow
\textsc{PruneSweep}
(\Iact,\Wold,\{(\Uactj{j},\bj{j})\}_{j=1}^{\nM},\KGs,\alphaG,\ntry)\)
\tcp{Algorithm~\ref{alg:graph_candidate_sweep}   \label{alg3:alg2}}

\eIf{\(\SUCCESS=\TRUE\)}{
\(\Iact\leftarrow\Irem\) \tcp{Update the current active support after successful pruning \label{alg3:up}}
\(\Wold\leftarrow\Wnew\)
\tcp{\(\Wold\in\mathbb{R}^{|\Iact|\times\nM}\), now supported on \(\ZD(\Iact)\) \label{alg3:up2}}
}{
\textbf{stop}
\tcp{No further feasible pruning step was found \label{alg3:finish}}
}
}

\(\Zomega\leftarrow\{z_i\}_{i\in\Iact}\)\;
\(\Wadc\leftarrow\Wold\)\;

\caption{\textsc{MAWPrune}: global feasible pruning of the initial
fixed ECM rule through repeated calls to
Algorithm~\ref{alg:graph_candidate_sweep}. The algorithm progressively
reduces the common support of the adaptive cubature rule while
preserving the local exactness constraints over the sampled latent
manifold.}\label{alg:MAW_global_pruning}
\end{algorithm}

Finally, Algorithm~\ref{alg:MAW_global_pruning} describes how the initial
fixed ECM rule introduced in
Section~\ref{sec:MAW_local_constraints} (see, in particular,
Eq.~\eqref{eq:initialW}) is progressively transformed into a more compact
 adaptive cubature rule through repeated calls
(Line~\ref{alg3:alg2}) to the candidate-sweep procedure described above. The process terminates when no further feasible
pruning step can be found (Line~\ref{alg3:finish}). The resulting matrix
\(
\Wad\in\mathbb{R}^{|\Zomega|\times\nM}
\)
contains the sampled representation of the desired manifold-adaptive weight
fields supported on the retained set of Gauss points
\(\Zomega=\{g_1,\ldots,g_{\nomega}\}\subseteq\ZD\). More precisely, the
\(i\)-th row of \(\Wad\) stores the samples over the latent space of the
 weight field associated\footnote{It should be noted that, with a mild
abuse of notation, the symbol \(\Wad\) is used both for the sparse
solution of the original global sparsification
problem~\eqref{eq:MAW_global_l0} (represented on the initial support
\(\ZD\) by inserting zeros at the eliminated points, hence
\(\Wad\in\mathbb{R}^{\mInit\times\nM}\)) and for its compressed dense
representation restricted to the retained support \(\Zomega\), for which
\(\Wad\in\mathbb{R}^{|\Zomega|\times\nM}\).}  with the retained Gauss point \(g_i\in\Zomega\).
%

\subsection{Regression of the adaptive weight fields}
\label{sec:regression}
The matrix \(\Wad\) is not itself the final outcome of the procedure; it
only provides the sampled data required for constructing the continuous
manifold-adaptive weight fields entering
Eq.~\eqref{eq:approxGEN2}. More specifically, for each retained Gauss
point \(g_i\in\Zomega\), \(i=1,\ldots,\nomega\), one seeks to construct a nonlinear map
\begin{equation}
\label{eq:3dssdddd}
    \qM
    \mapsto
    \omegag{g_i}(\qM),
\end{equation}
from the pairs
\begin{equation}
\label{eq:pairs}
  \left\{
    \big(\qMj{j},(\Wad)_{i,j}\big)
    \right\}_{j=1}^{\nM}
\end{equation}.

Generic regression procedures for Eq.~\ref{eq:3dssdddd} do not guarantee
positivity or exact preservation of the total volume, particularly outside the convex hull of the training data. Possible strategies for
mitigating this issue include, on the one hand, enriching the regression
dataset with additional low-fidelity samples outside the FE training
region\footnote{The additional samples could be generated
  using the manifold HROM with fixed weights. Incidentally, the
same multi-fidelity philosophy could be applied to the decoder by augmenting its training
set with additional low-fidelity snapshots generated, for instance, from a
standard fixed-weight HROM.}, or, on the other hand,
employing regression procedures explicitly designed to preserve positivity, volume,
or related geometric constraints~\cite{pennec2006riemannian}.  The study of such
constraint-preserving regression procedures lies beyond the scope of the
present work; in the numerical examples discussed in the sequel, we simply
employ the same unconstrained regression techniques used for the nonlinear
decoder.

\subsection{Summary}
\label{sec:summary}

\begin{BOX}
\caption{Offline construction of a MAW--ECM rule.}
\label{box:MAW_ECM_summary}

\footnotesize

\begin{enumerate}
\setlength{\itemsep}{2pt}
\setlength{\parsep}{0pt}

\item Generate the FE training trajectories over \(\PmuTRAIN\) and
collect the independent displacement snapshots, together with the
Gauss-point and history data required to reconstruct the internal-force
and output integrands; see
Section~\ref{sec:FOM}  and
Eq.~\eqref{eq:manifold_snapshot_matrix}.

\item Construct the decoder \(\Decod\) and encoder \(\Encod\), identify
the latent coordinates of the training states, and evaluate the tangent
basis \(\PhiD(\qM)\). For the factorized decoder adopted here, this
amounts to determining \(\PhiROM\), \(\TAU(\qM)\), and its Jacobian
\(\JTAU(\qM)\); see
Eqs.~\eqref{eq:dredNON}, \eqref{eq:inversemapping},
\eqref{eq:decoder_tau_generic}, and~\eqref{eq:JTAU}.

\item At each training state, assemble the generic manifold-integrand
block \(\Agen_j\) associated with projected equilibrium and/or output
evaluation, and concatenate the resulting blocks into the
global matrix \(\Agen\); cf.
Eqs.~\eqref{eq:rDec_reduced_operator},
\eqref{eq:Agen_j_definition}, and~\eqref{eq:AgenDEF}.

\item Compress the column space of \(\Agen\) to obtain \(\Ucub\), and
compute the initial fixed-weight ECM rule from
Problem~\eqref{eq:ECM_problem_manifold_fixed}. This provides the
candidate support
\(\ZD=\supp(\wD)\), containing \(\mInit\) Gauss points, together with
the corresponding positive weights; see
Eqs.~\eqref{eq:Ucub_definition}, \eqref{eq:ZD_support},
and~\eqref{eq:mInit_fixed_ECM}.

\item Select the representative parameter subset \(\PmuM\) and its
latent states according to Eq.~\eqref{eq:PmuM_definition}. Initialize
the sampled adaptive weights by restricting and replicating the
fixed-weight rule, as in Eq.~\eqref{eq:initialW}, and construct the
local exactness systems
\[
\Uj{j}^{T}\omegaj{j}=\bj{j},
\qquad j=1,\ldots,\nM,
\]
by projecting the local integrand blocks, adjoining the prescribed
invariant conditions, and removing linearly dependent constraints;
cf. Eqs.~\eqref{eq:MAW_projected_local_block}--%
\eqref{eq:constraintsLOCAL}.

\item Construct the latent graph \(\Glat\), its graph operator
\(\KGs\), and choose the regularization parameter \(\alphaG\).

 \item Apply Algorithm~\ref{alg:MAW_global_pruning}
(\textsc{MAWPrune}) to the initial fixed-weight ECM rule. The algorithm
returns the retained support
\(\Zomega\subseteq\ZD\) and the matrix
\(\Wad\in\mathbb{R}^{|\Zomega|\times\nM}\) containing the sampled
adaptive weights; see Section~\ref{sec:globalpruning}.

\item For each retained Gauss point \(g_i\in\Zomega\), regress the
\(i\)-th row of \(\Wad\) over the sampled latent states to construct
the continuous weight field
\[
\qM\mapsto\omegag{g_i}(\qM),
\qquad i=1,\ldots,\nomega,
\]
from the data pairs in Eq.~\eqref{eq:pairs}; see
Section~\ref{sec:regression}. The resulting online cubature rule has a
fixed common support \(\Zomega\), while its weights vary over the
latent manifold.

\end{enumerate}

\end{BOX}

By way of conclusion, Box~\ref{box:MAW_ECM_summary}
summarizes the complete offline construction of the MAW--ECM rule.

\section{Benchmark 1: homogenization of a  metamaterial unit  cell}
\label{sec:metamaterial}


\subsection{Problem setting}
 \label{sec:meta_problem_setting}

\begin{figure}[!ht]
    \centering
    \includegraphics[width=0.8\textwidth]{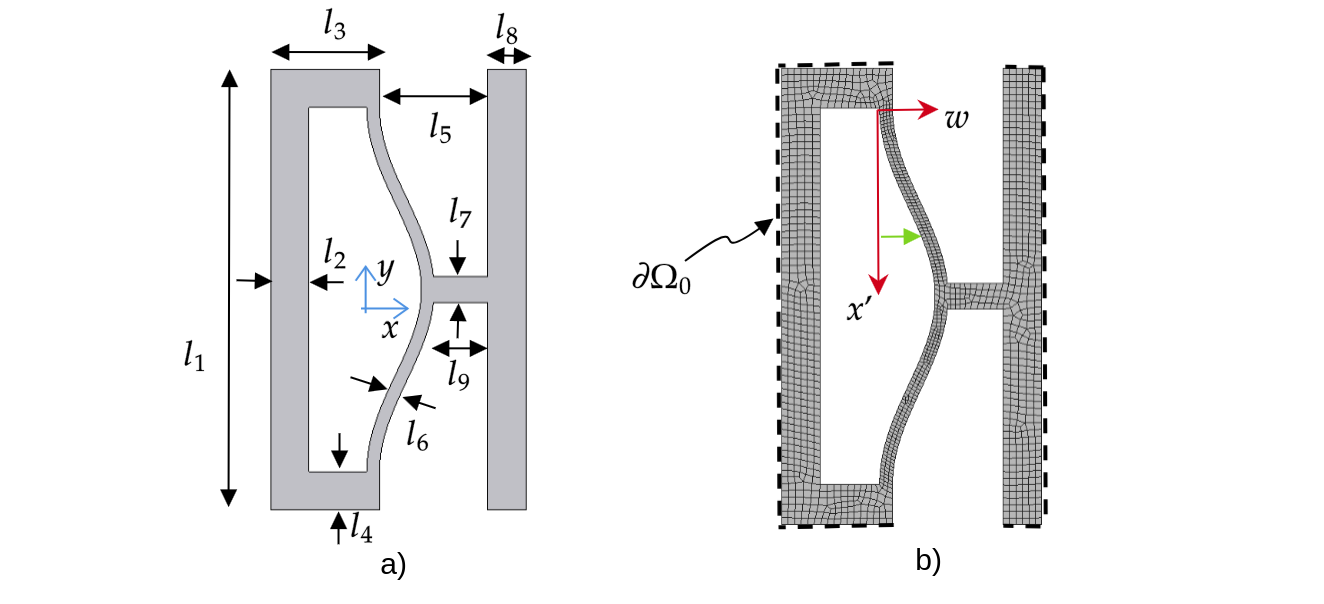}
    \caption{Geometry and finite element discretization of the
    metamaterial unit cell (the geometry is based on the design proposed in Ref.~\cite{chen2025negative}).
    (a) Reference domain, with
    \(l_1=60.40~\mathrm{mm}\),
    \(l_2=l_4=l_8=5.20~\mathrm{mm}\),
    \(l_3=14.50~\mathrm{mm}\),
    \(l_5=14.90~\mathrm{mm}\),
    \(l_6=1.76~\mathrm{mm}\),
    \(l_7=3.50~\mathrm{mm}\), and
    \(l_9=7.49~\mathrm{mm}\).
    The centerline of the lower
    curved ligament is described in the local coordinates
    \((x',w)\) by
    \(w(x')=\frac{h}{2}
    [1-\cos(\pi x'/L)]\), \(0\leq x'\leq 2L\), where
    \(2L=l_1-2l_4=50.00~\mathrm{mm}\) is its projected span and
    \(h=l_5-l_9-l_6=5.65~\mathrm{mm}\) is its maximum transverse rise.
    (b) Finite element mesh; the dashed markers identify the external
    boundary \(\partial\Omega_0\), whose opposite sides are paired for
    the imposition of periodic boundary conditions.}
    \label{fig:cell_geometry}
\end{figure}

The first numerical benchmark considers the two-dimensional periodic
``metamaterial'' unit cell\footnote{A \emph{metamaterial} is an architected material formed by the
repetition of unit cells whose effective macroscopic response differs
markedly from that of the constituent solid alone. In the present case,
snap-through of the curved ligaments produces negative incremental
stiffness under compression and a pronounced asymmetry between the
tensile and compressive responses, making the problem a challenging
benchmark for model reduction. Further details on mechanical
metamaterials may be found in
Refs.~\cite{qiu2004curved,correa2015mechanical,chen2025negative}.} displayed in  Fig.~\ref{fig:cell_geometry}.
The  material is modeled as an isotropic compressible
Neo--Hookean solid under plane-strain conditions, with Young's modulus
\(E=1628~\mathrm{MPa}\) and Poisson's ratio \(\nu=0.4\). The domain is
discretized using \(\nel = 1315\) bilinear quadrilateral elements (average size $0.1~mm$) and \(1602\)
nodes; $\ngpel = 4$  Gauss points  are employed per element, giving a total of
\(\ngp=\nel \ngpel =  5260\) integration points.

The cell is deformed    exclusively through its external boundary $\partial \Omega_0$ via imposed displacements (thus
\(\Fext=\zero\) in Eq.~\refeq{eq:Rfe_definition}). Furthermore, the boundary term $\dBOUND$ in Eq.~\refeq{eq:dvecreconstruction} is assumed to be linear with the input parameter $\muvec$, i.e., $\dBOUND=\Umacro\muvec$,
  the columns of $\Umacro$ defining the prescribed displacement patterns for unit variations of each component of $\muvec$. To give this   parametrization a concrete mechanical
interpretation, \(\Umacro\) is constructed from the periodic boundary
conditions of first-order computational homogenization\footnote{First-order homogenization is employed here solely as a practical setting for illustrating the proposed strategy, as it provides a conveniently low-dimensional input space in terms of macroscopic deformation measures. We are aware that its macroscopic kinematics are too restrictive to represent the full complexity of bistable unit cells, for which richer nonlinear multiscale formulations involving additional microstructural fields, such as micromorphic variables, are more appropriate; see, e.g., Ref.~\cite{sperling2024comparative}.}
  at finite
strains; see, e.g., Ref.~\cite{miehe2003computational}. 
  More specifically, the input parameters are identified with two   components (thus $\nmu = 2$) of the macroscopic displacement gradient $\Gmacro \in \RRn{2}{2}$:  the normal
component along the \(x\)-direction and the shear component, namely,
\begin{equation}
\label{eq:meta_macro_gradient}
    \muvec
    =
    \begin{bmatrix}
        \Gmacrox\\
        \Gmacroxy
    \end{bmatrix}.
\end{equation}
The transverse normal component is set to zero, and \(\Gmacro\) is assumed
symmetric, so that no macroscopic rigid-body rotation is prescribed; thus:
\begin{equation}
    \Gmacro(\muvec)
    =
    \begin{bmatrix}
        \Gmacrox  & \Gmacroxy\\
        \Gmacroxy & 0
    \end{bmatrix} =\begin{bmatrix}
        \mu_1  & \mu_2\\
        \mu_2 & 0
    \end{bmatrix}.
\end{equation}
Nonzero entries in $\Umacro$ appear only for slave boundary nodes   and   corner
nodes; such entries are determined by considering that,  for a node $\bm{X}=(X,Y)^T$, the   macroscopic contribution is defined   as $\bar{\bm{u}}(\bm{X})=\Gmacro \bm{X}$, which in the present
case yields $\bar{\bm{u}}(\bm{X})=( \mu_1 X+ \mu_2 Y,\,
\mu_2 X)^T$.  Consistently with this interpretation, \(\dred \in \Rn{\nL}\)  in
Eq.~\refeq{eq:dvecreconstruction} collects the degrees of freedom of
the interior nodes and the independent boundary nodes on the
designated master sides for imposing periodic conditions (here $\nL = 3018$). On the other hand, the lifting operator \(\Tlift\) in Eq.~\refeq{eq:dvecreconstruction} transfers the
periodic fluctuation displacements from the master sides to their
paired slave sides.

As customary in first-order computational homogenization, the output of
interest is the homogenized first Piola--Kirchhoff stress,
\begin{equation}
\label{eq:meta_homogenized_stress}
    \PKmacro
    =
    \frac{1}{|\Omega_0|}
    \int_{\Omega_0}\PKone\,d\Omega_0,
\end{equation}
which is a volumetric output \(\yOUT\) of the form introduced in
Eq.~\refeq{eq:generic_output_FE_quadrature}; hence,  $\rOUT = \PKone/|\Omega_0| $ and  \(\ncond=4\).


\subsection{Training data}

\begin{figure}[!ht]
    \centering
    \includegraphics[width=0.8\textwidth]{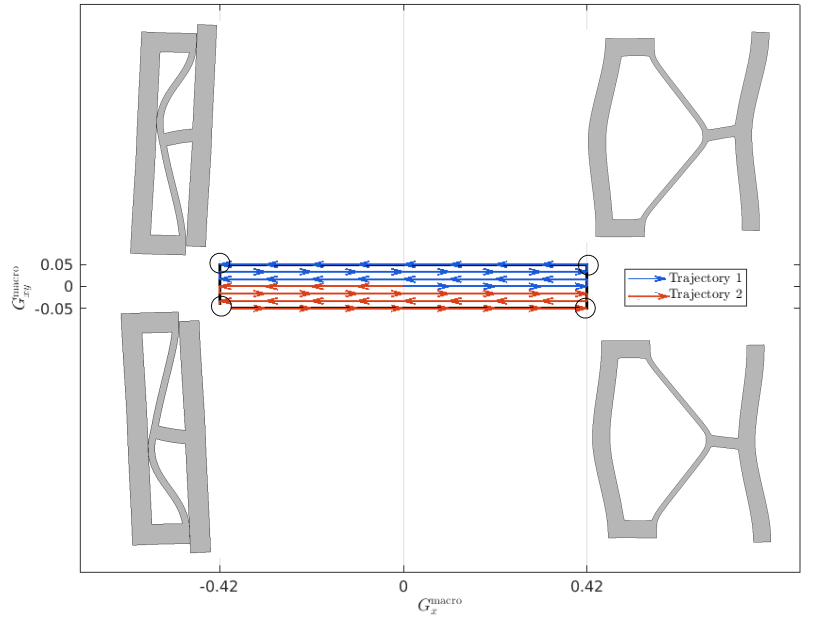}
    \caption{Training parameter domain and zig-zag trajectories used to
generate the FE displacement snapshots. The corner configurations
correspond to the extreme parameter values.}
    \label{fig:param_space}
\end{figure}

The training set \(\PmuTRAIN\) (see Eq.~\ref{eq:pmutrain})  is   restricted to a
two-dimensional rectangle, in order to keep the regression problem simple and to
allow a direct graphical inspection of the adaptive weights. The compression/stretching
component along the soft direction of the cell is varied in the interval
$ \Gmacrox \in [-0.42,\,0.42]$, whereas the shear component is varied in
$\Gmacroxy \in [-0.05,\,0.05]$. As illustrated in Fig.~\ref{fig:param_space}, the lower bound
\((\Gmacrox)_{\textrm{min}}=-0.42\) corresponds to the most compressed configuration before self-contact occurs\footnote{Our model does not include any
representation of self-contact; hence, extrapolation beyond this
compression level is physically meaningless.}.    Although these cells are primarily intended to operate under compression,
the parameter domain is  kept symmetric, with
\((\Gmacrox)_{\textrm{max}} = 0.42\). The purpose is to test the reduced model in
a prediction regime in which input parameters of the same magnitude but
opposite sign generate outputs of markedly different magnitude (as shown
later, in Figure~\ref{fig:standard_HROM_stress}, the corresponding stress levels in the $x$ direction  may differ by up to a factor of
\(25\), since buckling of the curved ligaments leads to a substantially softer
response under compression).

We sample the rectangular parameter domain   using a
Cartesian grid with \(230\) and \(20\) uniform intervals in the
\(\Gmacrox\)- and \(\Gmacroxy\)-directions, respectively. This gives
\(\nsnap=231\cdot21=4851\) distinct parameter states.  The FE simulations corresponding to each of these parameters are obtained by   traversing this grid along two continuous zig-zag
loading trajectories, covering the nonnegative and nonpositive shear
regions, respectively (see Fig.~\ref{fig:param_space}).

\subsection{Standard HROM  }
 \label{sec:meta_standard_hrom}

We begin by assessing the
efficiency of the standard HROM. This
model is recovered as a particular case of the manifold HROM introduced
in Section~\ref{sec:decoderENCODER} by identifying, in the decoder
factorization of Eq.~\refeq{eq:decoder_tau_generic},  \(\qM\) with the standard reduced coordinates
\(\qrom\), and by taking \(\TAU(\qrom)=\qrom\). It then follows that
\(\JTAU=\ident\) in Eq.~\refeq{eq:JTAU} and
\(\PhiD=\PhiROM\) in Eq.~\refeq{eq:PhiD_tau_generic}.

\begin{figure}[!ht]
    \centering
    \begin{subfigure}[t]{0.2\textwidth}
        \centering
        \includegraphics[width=\textwidth]{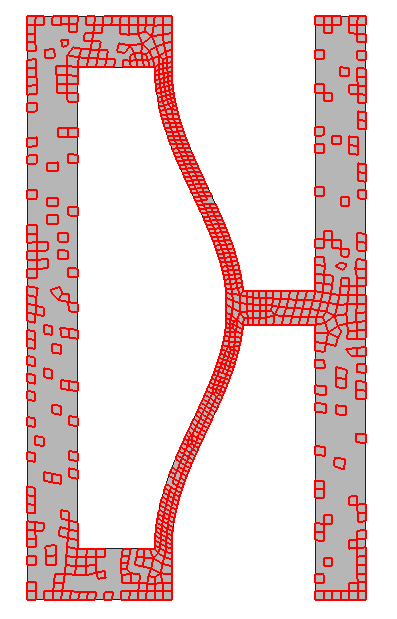}
        \caption{Elements containing the selected Gauss points when using standard HROM.}
        \label{fig:standard_HROM_support}
    \end{subfigure}
    \hfill
    \begin{subfigure}[t]{0.7\textwidth}
        \centering
        \includegraphics[width=\textwidth]{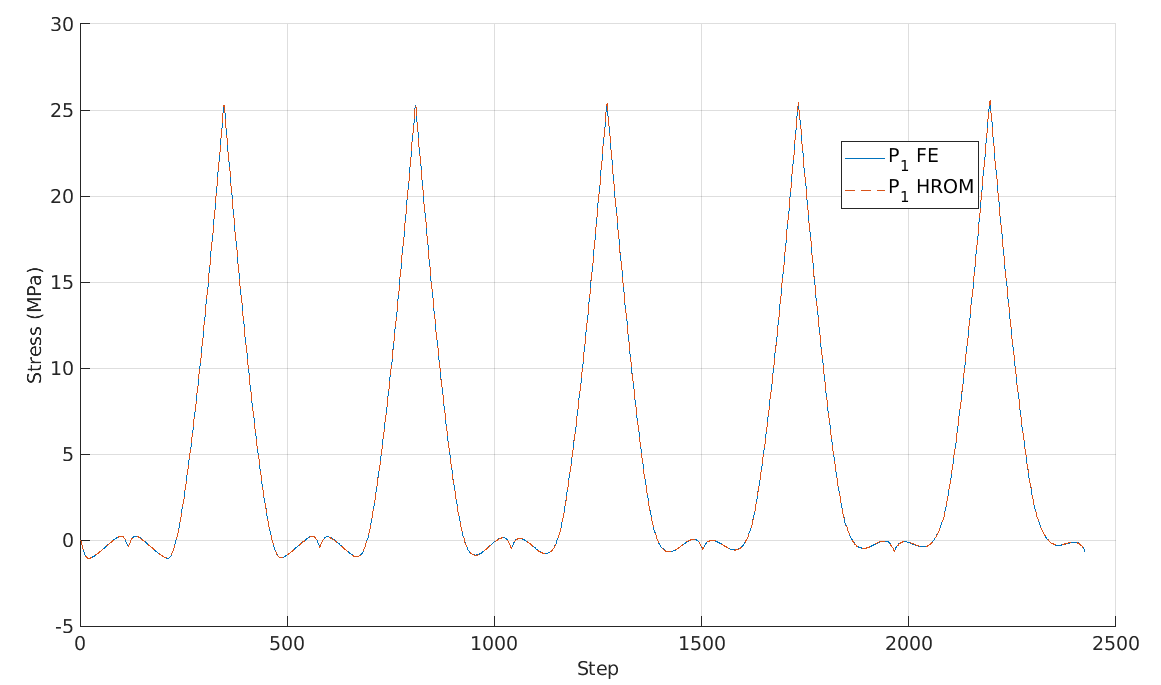}
        \caption{Homogenized \(P_x\) response.}
        \label{fig:standard_HROM_stress}
    \end{subfigure}
    \caption{Standard-HROM results for
    \(\tolSVDd=10^{-6}\) and \(\tolSVD=10^{-4}\).
    (a) The ECM rule comprises \(1340\) Gauss points distributed over
    \(671\) finite elements. (b) FE and standard-HROM homogenized PK1
    stress histories along the second trajectory in
    Fig.~\ref{fig:param_space}; the HROM uses \(\rROM=52\) displacement
    modes.}
    \label{fig:standard_HROM_results}
\end{figure}

 The accuracy and complexity of this model are governed solely by the
truncation tolerances $\tolSVDd$, for determining the displacement basis  \(\PhiROM\), see Eq.~\ref{eq:PhiROM_truncation}, and   \(\tolSVD\), for the integrand matrix \(\Ucub\), see Eq.~\refeq{eq:Ucub_definition}.  Both bases are computed using the
Sequential Randomized SVD   of
Ref.~\cite{hernandez2024cecm}, briefly introduced in
Section~\ref{sec:fixed_weight_cubature}.\footnote{For
\(\rROM=52\), the integrand matrix \(\Agen\) has \(5260\) rows, one per
candidate Gauss point, and
\((\rROM+\nOUT)\nsnap=(52+4)\times4851=271656\) columns. Its explicit
storage in double precision requires approximately \(11.4\,\mathrm{GB}\),
excluding factorization costs. The SRSVD is therefore particularly useful
under memory limitations, as it processes the state-wise blocks
\(\Agen_j\) sequentially without assembling the complete matrix. \label{footnote:size}}
We then performed a systematic tolerance study and selected the least
expensive configuration for which the normalized maximum error in the
homogenized stress remained below \(0.1\%\) over the complete training
set. The study revealed that this criterion requires \(\tolSVDd=10^{-6}\), yielding a basis matrix $\PhiROM$ with
\(\rROM=52\) modes, and  \(\tolSVD=10^{-4}\), which results in an
ECM rule with \(1340\) Gauss points. For illustration,
Fig.~\ref{fig:standard_HROM_results}(a) shows the \(671\) finite elements
containing these points. The tight displacement truncation tolerance
\(\tolSVDd=10^{-6}\) required to meet the accuracy criterion for the stresses reflects the previously mentioned
markedly different scales of the tensile and compressive responses shown
in Fig.~\ref{fig:standard_HROM_results}(b). Indeed, modes associated with small
singular values contribute little to the snapshot energy, yet are necessary to  capture local details related to
 internal buckling mechanisms governing the low-stress compressive
branch.

In summary, even though the problem is governed by  \(\nmu=2\) parameters,    $\rROM=52$ linear displacement modes are required
to reproduce the nonlinear response with the desired accuracy. This  illustrates the poor
linear compressibility of the problem under study. The
limitation is even more apparent at the hyperreduction level, where the
fixed ECM rule retains  1340 of the original $\ngp=5260$ integration points (approximately 25\%).

\subsection{Manifold HROM}
\subsubsection{Identification of latent coordinates}
\label{sec:meta_manifold_representation}

 The poor linear compressibility observed above motivates replacing the
linear reduced representation with the nonlinear-manifold approximation
of Eq.~\refeq{eq:decoder_master_slave}, introduced in
Section~\ref{sec:fact}. We first address the identification of the
latent coordinates. As established in Eq.~\refeq{eq:qMdefq}, these are
obtained from the modal coefficients through
\(\qM=\Tm\qrom\); following the input-informed strategy discussed in
Section~\ref{sec:inputinformed}, the task here is therefore to determine
the transformation matrix \(\Tm\) by exploiting the particular
parametrization of the present benchmark.

The  problem under consideration is static (no inertial effects) and the material response is hyperelastic (no path dependence). It can be argued that, under such circumstances,   the displacement
state $\dred$ admits a single-valued parametrization in terms of prescribed
  input \(\muvec\). This implies that the dimension of the solution manifold is   equal to the dimension of the input parameter space ($\rD=\nmu=2$), and, furthermore, that   there is a one-to-one mapping between input parameters and the desired latent coordinates:
  \begin{equation}
   \qM = \varphiM(\muvec).
  \end{equation}
  For the present benchmark, \(\Tm\) is determined by requiring that the above map     approximate the identity
map over the training set. In view of Eq.~\refeq{eq:qMdefq}, this
amounts to making $
    \qM(\muvec_j)-\muvec_j
    =
    \Tm\qrom(\muvec_j)-\muvec_j
$ as small as possible for every training input \(\muvec_j\) ($j=1,2 \ldots \nsnap$). Collecting
the modal coefficients and the corresponding input parameters into
\begin{equation}
\label{eq:meta_coordinate_matrices}
    \Qrom
    \defeq
    \PhiROM^T\Dsnap,
    \qquad
    \Mtrain
    \defeq
    \begin{bmatrix}
        \muvec_1 & \muvec_2 & \cdots & \muvec_{\nsnap}
    \end{bmatrix},
\end{equation}
(here, \(\Dsnap\) is the displacement snapshot matrix introduced in
Eq.~\refeq{eq:manifold_snapshot_matrix}), the matrix \(\Tm\) is accordingly obtained from the least-squares problem
\begin{equation}
\label{eq:least_squares_master_coordinates}
    \Tm
    =
    \arg\min_{\widetilde{\Tm}}
    \left\|
        \widetilde{\Tm}\Qrom-\Mtrain
    \right\|_F^2.
\end{equation}
  Assuming $\Qrom$ has full rank, the  solution of Problem \refeq{eq:least_squares_master_coordinates} is unique and  given by
\begin{equation}
\label{eq:least_squares_master_coordinates_solution}
    \Tm
    =
    \Mtrain\Qrom^\dagger
\end{equation}
where \(\Qrom^\dagger\) denotes the Moore--Penrose pseudoinverse.

The relative value of the least-squares residual at the optimum for the problem under consideration is
\begin{equation}
\label{eq:latent_parameter_residual}
    \frac{
        \left\|\Tm\Qrom-\Mtrain\right\|_F
    }{
        \left\|\Mtrain\right\|_F
    }
    =
    3.09\cdot10^{-7}.
\end{equation}
Thus, to numerical accuracy, the latent coordinates associated with the
training states can be identified with the   input parameters:
\begin{equation}
\label{eq:input_aligned_coordinates}
    \qM\approx\muvec.
\end{equation}

\begin{remark}
\label{remark:noiterations}
This identification has an important consequence for the online stage
of the present benchmark. Because \(\qM\) may be obtained directly from the
prescribed input \(\muvec\), the reduced equilibrium equations
\refeq{eq:MAW_ROM_equilibrium} need not be solved to determine the
latent state. The decoder can instead be evaluated directly at
\(\qM\approx\muvec\). Hyperreduction through the proposed adaptive
cubature rule is therefore required only for evaluating the homogenized
stresses in Eq.~\refeq{eq:meta_homogenized_stress} (i.e., the equilibrium-residual integrands need not be included in
the cubature training problem).
\end{remark}

\subsubsection{Nonlinear closure}
\label{sec:nonlinearclosure}
\begin{figure}[!ht]
    \centering

    \begin{subfigure}{0.21\textwidth}
        \centering
        \includegraphics[width=\textwidth]{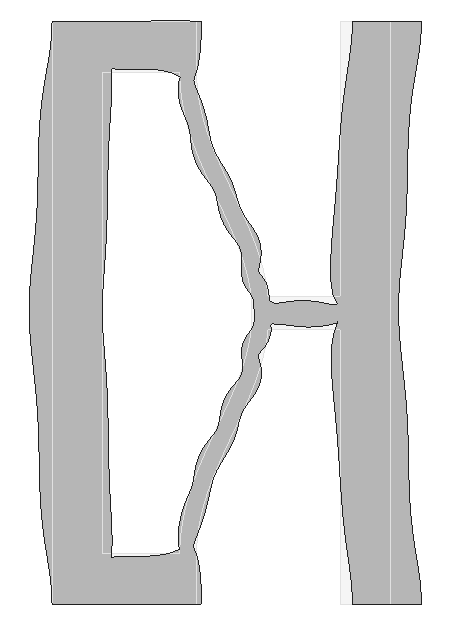}
        \caption{Master mode \({\PhiM}_1\).}
                \label{fig:master_mode_1}
    \end{subfigure}
    \hspace{0.05\textwidth}
    \begin{subfigure}{0.21\textwidth}
        \centering
        \includegraphics[width=\textwidth]{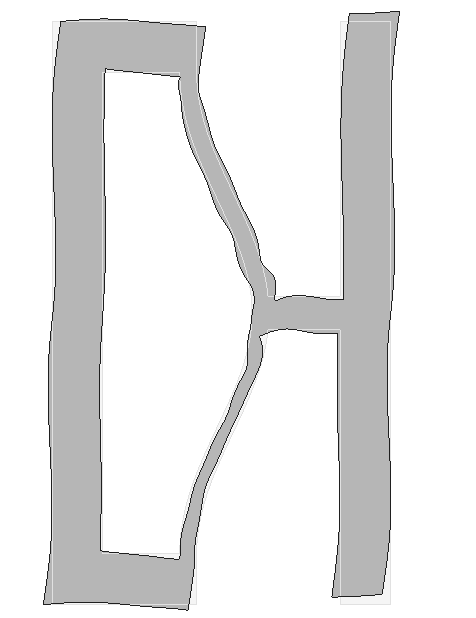}
        \caption{Master mode \({\PhiM}_2\).}
                \label{fig:master_mode_2}
    \end{subfigure}

    \vspace{0.5em}

    \begin{subfigure}{0.21\textwidth}
        \centering
        \includegraphics[width=\textwidth]{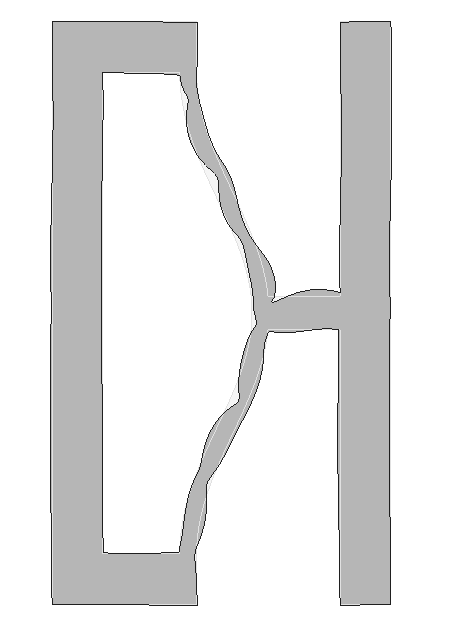}
        \caption{Slave mode \({\PhiS}_1\).}
          \label{fig:slave_mode_1}
    \end{subfigure}
    \begin{subfigure}{0.21\textwidth}
        \centering
        \includegraphics[width=\textwidth]{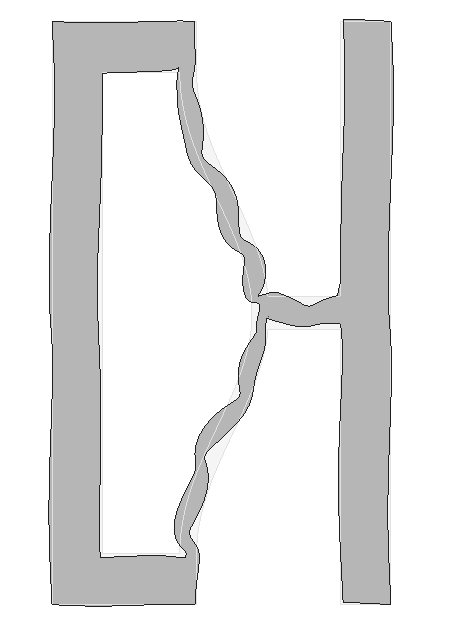}
        \caption{Slave mode \({\PhiS}_6\).}
         \label{fig:slave_mode_2}
    \end{subfigure}
    \begin{subfigure}{0.21\textwidth}
        \centering
        \includegraphics[width=\textwidth]{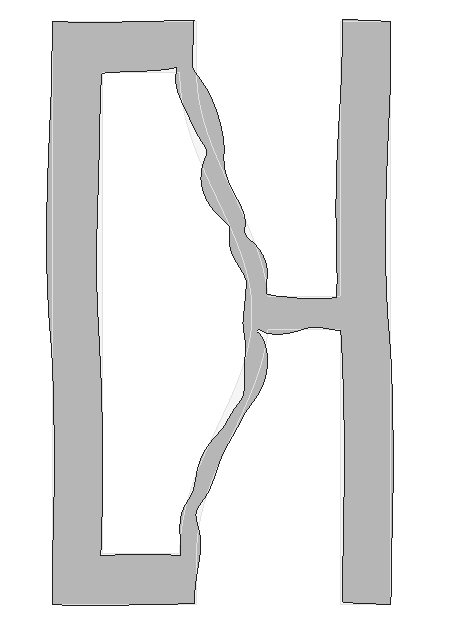}
        \caption{Slave mode \({\PhiS}_{13}\).}
        \label{fig:slave_mode_3}
    \end{subfigure}
    \begin{subfigure}{0.18\textwidth}
        \centering
        \includegraphics[width=\textwidth]{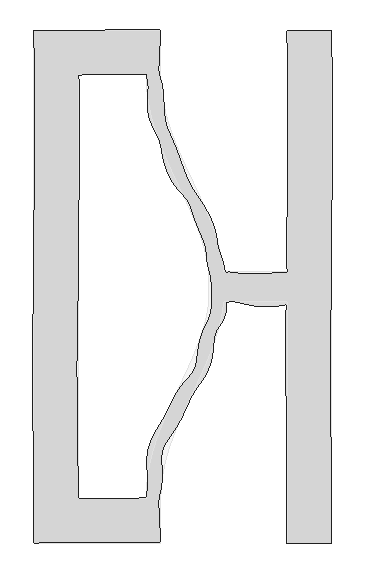}
        \caption{Slave mode \({\PhiS}_{44}\).}
        \label{fig:slave_mode_4}
    \end{subfigure}

    \caption{Deformed shapes of the $\rD = 2$ master modes $\PhiM$ and four of the  $\rROM-\rD = 50$    slave modes $\PhiS$.  The modal vectors,
    originally defined on the independent degrees of freedom, are extended
    to the complete nodal space through the lifting operator \(\Tlift\)
    for visualization (see Eq.~\refeq{eq:dvecreconstruction}).}
    \label{fig:master_slave_modes_manifold}
\end{figure}

Once \(\Tm\) has been determined, the master basis \(\PhiM\) and the
scaling matrix \(\Am\) follow from the SVD of
\(\PhiROM\Tm^T\) introduced in Section~\ref{sec:fact}, together with
Eq.~\refeq{eq:Am_definition}. The slave basis \(\PhiS\) is subsequently
obtained from the orthogonality conditions in Eq.~\refeq{eq:orto}. For
the present problem, this decomposition yields the $\rD = 2$ master modes whose deformed shapes are depicted in Figures \ref{fig:master_mode_1} and \ref{fig:master_mode_2}, and
\(\rROM-\rD=50\) slave modes (the deformed shapes of some of these modes are in turn shown in Figures \ref{fig:slave_mode_1} to \ref{fig:slave_mode_4} ).

%

  It remains to identify the nonlinear closure mapping the
\(\rD=2\) latent coordinates onto the
\(\rROM-\rD=50\) slave coordinates,
\(\qS\approx\Nslave(\qM)\), from the training pairs
\((\qM,\qS)\) obtained from the displacement snapshots. To this end, we employ anisotropic
Gaussian radial basis functions (RBFs) augmented with polynomial
terms (quadratic). The regression coefficients are determined through
standard-form Tikhonov-regularized least squares, using the identity
regularization matrix ( see
Refs.~\cite{poggio1990networks,wendland2005scattered}). Only a subset of the available snapshot pairs is retained as RBF
centers. A structured subsampling frequency of \([3,1]\) is used:
every third sample is retained along the first latent direction
(compression/extension), whereas all samples are retained along the
second direction (shear). This results in \(1540\) centers out of $\nsnap =4851$ data points. The
remaining states are employed to assess the accuracy of the closure away
from the retained centers. The characteristic lengths of the anisotropic Gaussian kernels are set to
\(\ell_1=3.3013\cdot10^{-2}\) and
\(\ell_2=1.5789\cdot10^{-2}\), and the Tikhonov regularization
parameter  to \(\lambda_{\mathrm{RBF}}=10^{-10}\). With these
settings, the relative reconstruction error is
\(4.30\cdot10^{-6}\) at the retained centers and
\(1.03\cdot10^{-4}\) over the complete dataset.

\begin{figure}[!ht]
    \centering

    \begin{subfigure}{0.45\textwidth}
        \centering
        \includegraphics[width=\textwidth]{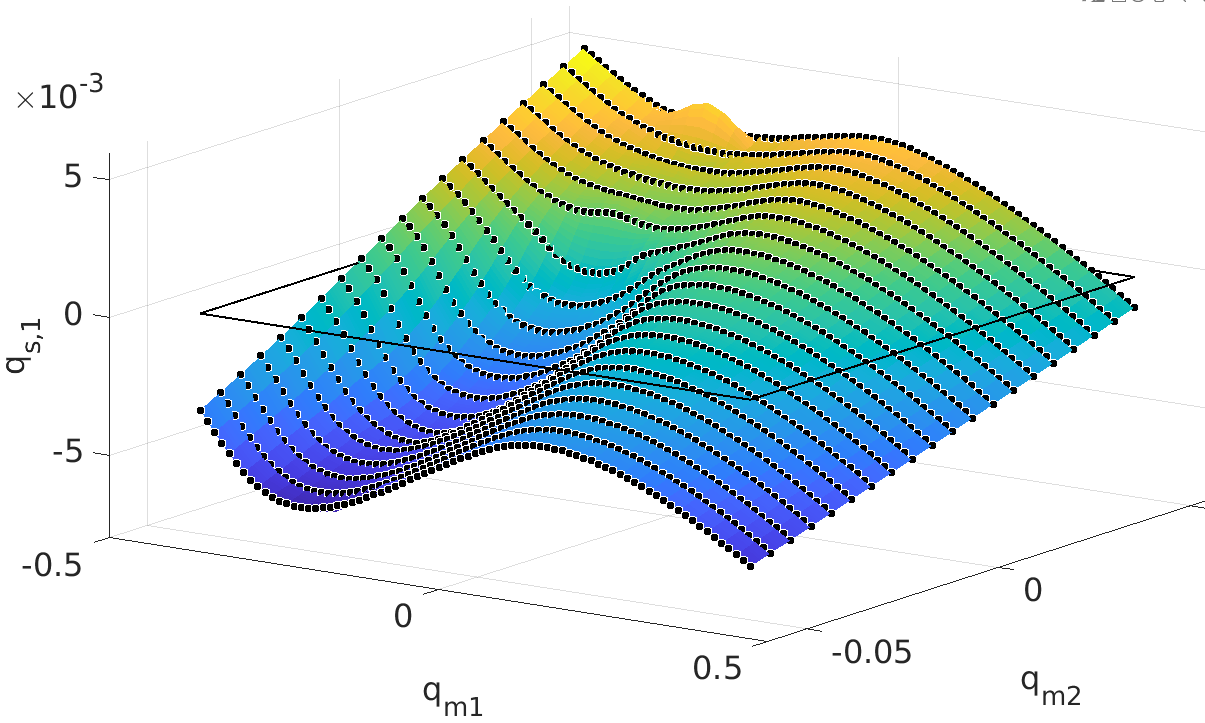}
        \caption{Slave amplitude \((\qS)_1\).}
        \label{fig:RBF_slave_surface_1}
    \end{subfigure}
    \hfill
    \begin{subfigure}{0.45\textwidth}
        \centering
        \includegraphics[width=\textwidth]{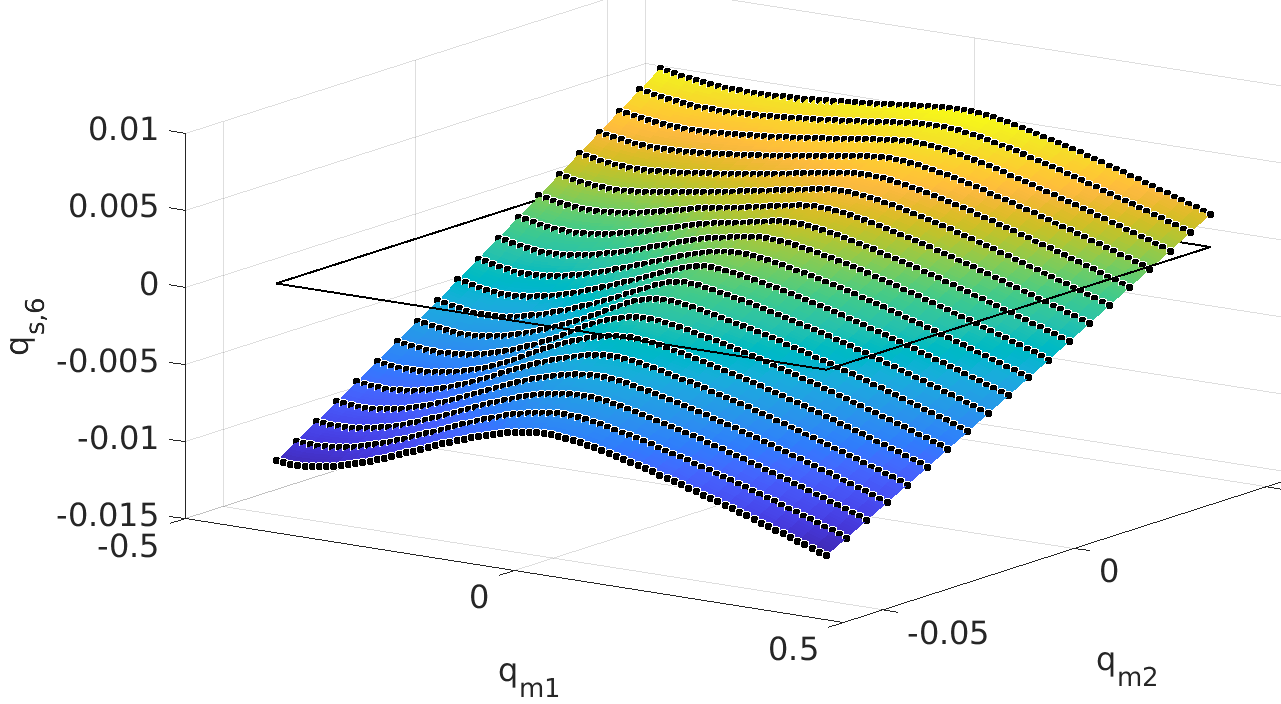}
        \caption{Slave amplitude \((\qS)_6\).}
        \label{fig:RBF_slave_surface_6}
    \end{subfigure}

    \vspace{0.5em}

    \begin{subfigure}{0.45\textwidth}
        \centering
        \includegraphics[width=\textwidth]{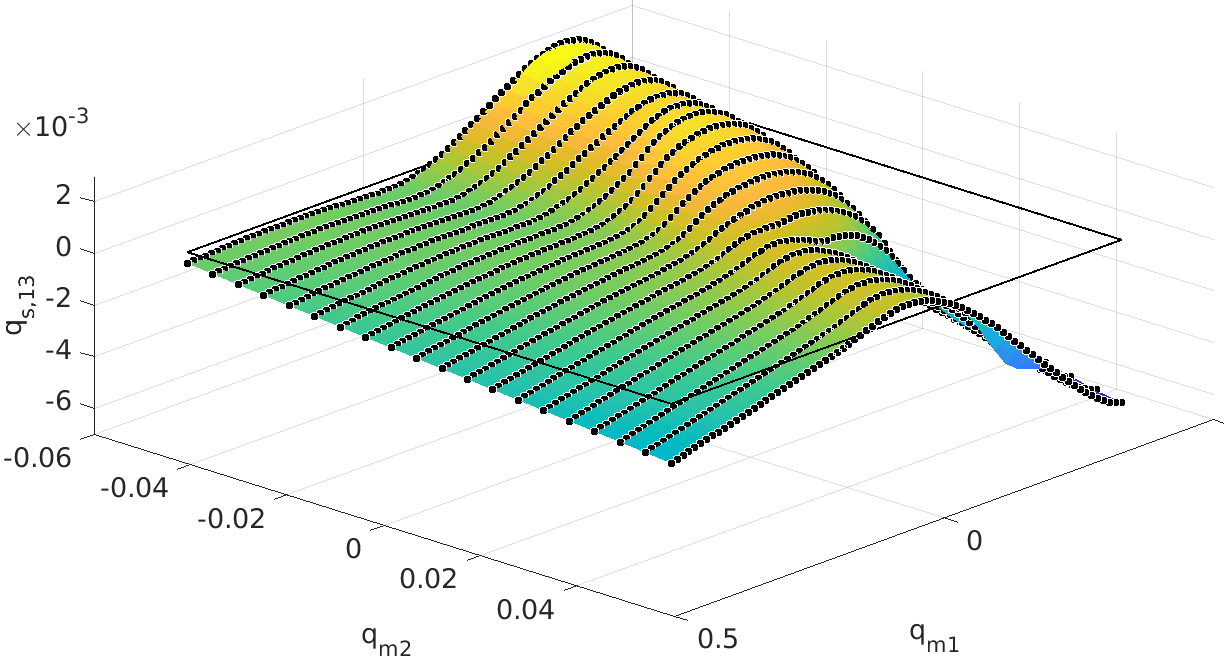}
        \caption{Slave amplitude \((\qS)_{13}\).}
        \label{fig:RBF_slave_surface_13}
    \end{subfigure}
    \hfill
    \begin{subfigure}{0.45\textwidth}
        \centering
        \includegraphics[width=\textwidth]{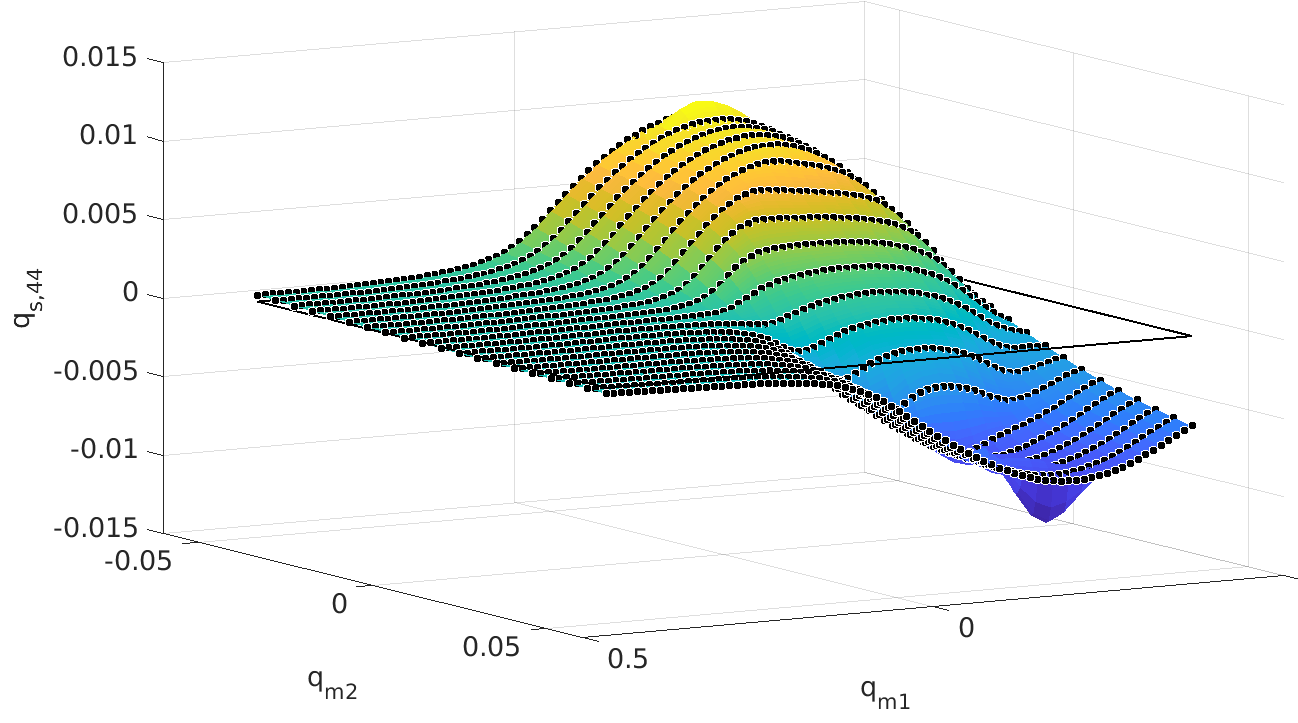}
        \caption{Slave amplitude \((\qS)_{44}\).}
        \label{fig:RBF_slave_surface_44}
    \end{subfigure}

    \caption{RBF reconstructions of four representative slave amplitudes
    over the two-dimensional latent space. The black markers denote the
    retained RBF centers.}
    \label{fig:meta_slave_closure}
\end{figure}

Figure~\ref{fig:meta_slave_closure} shows the resulting mappings for
four representative slave coordinates. The reconstructed surfaces are
smooth and single-valued over the sampled latent domain (the black markers indicate the retained RBF centers).

\subsubsection{Fixed-weight ECM}
\label{sec:fixedECM1}
With the manifold coordinates and the slave closure in place, we next
construct the fixed-weight ECM rule, discussed in Section~\ref{sec:fixedweight_cubature}. As noted in Remark~\ref{remark:noiterations},  the latent state is prescribed directly by the macroscopic
input, and no reduced equilibrium problem needs
to be solved for the evaluation of the homogenized response. We
therefore adopt the output-only case described after
Eq.~\refeq{eq:approxGEN} in Section~\ref{sec:fixedweight_cubature}  and construct the cubature rule exclusively
for the homogenized stress in Eq.~\refeq{eq:meta_homogenized_stress}, that is,  $\rGENg{g}= {\PKone_g}/{|\Omega_0|}$ ($g=1,2 \ldots \ngp$), which implies that there are  $\ncond = \nstress = 4$ integration conditions per state. The same \(\nM=1540\) retained latent states used for the RBF closure described previously in Section \ref{sec:nonlinearclosure} are employed to assemble the output-integrand matrix $\Agen$ defined in Eq.~\ref{eq:AgenDEF}; hence\footnote{This matrix is markedly smaller than its standard-HROM
counterpart, which contains
\((\rROM+\nOUT)\nsnap=271656\) columns; see
Footnote~\ref{footnote:size} in
Section~\ref{sec:meta_standard_hrom}. The number of scalar fields
processed by the SRSVD is therefore reduced by a factor of
approximately \(44\), substantially alleviating the memory demands
associated with the construction of the ECM basis. }  $\Agen \in
    \mathbb{R}^{\ngp\times(\nstress\nM)}     =     \mathbb{R}^{5260\times6160}$. Using the SRSVD with \(\tolSVD=10^{-4}\) to compress \(\Agen\), and
augmenting the resulting basis with the volume constraint in
Eq.~\refeq{eq:ECM_volume_constraint}, the ECM yields a fixed rule with
\(\mInit=54\) Gauss points distributed over \(48\) finite elements.

Table~\ref{tab:online_fixed_manifold} compares  the accuracy of the fixed-weight manifold
approximation with those of the standard ROM and HROM. Here,
\(\errdisp\) denotes the relative Frobenius error in the independent
displacement histories, whereas \(\errPK\) is the corresponding
relative Frobenius error in the homogenized PK1 stress histories, in
both cases using the FE solution as reference. The errors are computed
separately over each complete loading trajectory and then averaged over
the two trajectories.

\begin{table}[!ht]
\centering
\caption{Online complexity (in terms of number of master/slave coordinates, and number of integration points  ) and accuracy  (relative Frobenius norms in displacements, $\errdisp$ and stresses $\errPK$)  of the fixed-weight manifold
HROM, together with the standard ROM and HROM.   }
\label{tab:online_fixed_manifold}
\begin{tabular}{lccccc}
\toprule
Model
& Master coords.
& Slave coords.
& Integ. points
& \(\errdisp\)
& \(\errPK\) \\
\midrule
ROM, full Gauss rule
& \(52\) & \(0\) & \(5260\)
& \(2.52\cdot10^{-6}\)
& \(5.05\cdot10^{-5}\) \\
HROM, ECM
& \(52\) & \(0\) & \(1340\)
& \(1.85\cdot10^{-5}\)
& \(7.13\cdot10^{-5}\) \\
M-HROM, fixed ECM
& \(2\) & \(50\) & \(54\)
& \(7.87\cdot10^{-5}\)
& \(2.98\cdot10^{-4}\) \\
\bottomrule
\end{tabular}
\end{table}

 Table~\ref{tab:online_fixed_manifold} shows that the manifold HROM
increases the homogenized-stress error by a factor of approximately four
relative to the standard HROM. This loss of accuracy is due to the    subsampling and subsequent regression (recall that  the
error is evaluated over all \(\nsnap=4851\) states, whereas the closure
is fitted using only \(\nM=1540\) retained centers). The resulting
accuracy--compression trade-off is nevertheless strongly favorable: the number of
independent coordinates decreases from \(\rROM = 52\) to \(\rD=2\), and the
number of integration points from \(1340\) to \(\mInit = 54\). Interestingly,
the corresponding compression ratios are almost identical, namely
\(26\) and approximately \(25\), respectively.

\subsection{Manifold-adaptive weight ECM }

The preceding fixed-weight manifold HROM completes the kinematic
compression, but its cubature rule remains global: the same
\(\mInit=54\) weights must represent the homogenized-stress integrands
over the entire latent domain. The final stage of the compression
process is therefore to exploit the manifold structure also at the
cubature level, by allowing the weights to vary with the latent state.
Following the formulation of Section~\ref{sec:MAWECM}, the fixed-weight
ECM rule provides the initial candidate support \(\ZD\) and the feasible
initial weight vector \(\omegaInit\). The greedy strategy of
Section~\ref{sec:Greedy} then progressively removes integration points
and redistributes the remaining weights over the sampled manifold while
preserving the local exactness and positivity conditions.

In the present output-only setting, the state-dependent integrands are
the $\nstress =4$ components of the   PK1 stress entering the
homogenized response. At each retained latent state, these four
integration conditions are supplemented with the invariant volume
constraint. This gives    a lower bound of $\mLower = 4+1 = 5$ points (see Eq.~\ref{eq:ncondmin}) for the desired
  adaptive rule .

The same \(\nM=1540\) states retained as RBF centers for the decoder are used as the
nodes of the latent graph \(\Glat\). Here, the numerical identification
\(\qM\approx\muvec\), together with the Cartesian sampling of the input
parameter domain (see Figure \ref{fig:standard_HROM_results}), makes the graph construction immediate: the retained
latent states inherit the connectivity of the structured parameter
grid. Adjacent states are connected to define an auxiliary
quadrilateral mesh in latent space, and the graph operator \(\KGs \in \RRn{1540}{1540}\) is
assembled from the standard finite element discretization of an
isotropic Laplacian using bilinear shape functions, consistently with
the graph-regularized redistribution described in
Section~\ref{sec:Greedy}.

\subsubsection{Influence of the fixed-weight initialization and pruning
efficiency}
\label{sec:infl_fw_ini}

We first assess the reduction achieved by the MAW--ECM pruning strategy
of Algorithm~\ref{alg:MAW_global_pruning}, and its sensitivity to the
fixed-weight initialization. The regularity of the resulting adaptive
weight fields, their regression over the latent domain, and the
associated online accuracy are examined subsequently.

Four initial cubature rules are considered, obtained with $\tolSVDfixed\in
\left\{10^{-4},10^{-5},10^{-6},10^{-7}\right\}$
in the truncation of the training integrand matrix \(\Agen\). The case
\(\tolSVDfixed=10^{-4}\) is precisely the \(\mInit=54\)-point rule
assessed in Table~\ref{tab:online_fixed_manifold}. The other three rules
are introduced here solely to determine whether the final adaptive
support and the pruning efficiency depend   on the size and
accuracy of the fixed-weight initialization.  The graph-regularization parameter and the number of candidate trials
are held fixed throughout the study at
\(\alphaG=10^4\), appearing in
Algorithm~\ref{alg:prune_step_optionB}, and
\(\ntry=5\), appearing in
Algorithm~\ref{alg:graph_candidate_sweep}. These values were selected
after preliminary experimentation as a compromise between regularity of
the adaptive weight fields and computational efficiency during pruning.

\begin{table}[ht!]
\centering
\caption{Influence of the fixed-weight initialization on the MAW--ECM pruning process for different values of the   SVD truncation tolerance $\tolSVDfixed$ for the integrand matrix $\Agen$.  Here, $\mInit$ denotes the number of points of the resulting initial fixed-weight ECM rule, and $\nomega$ the final number of points after pruning. \emph{First-stage limit} denotes the minimum number of points reached before positivity enforcement and graph regularization are activated in Algorithm~\ref{alg:prune_step_optionB}. \emph{Total time} corresponds to the wall-clock time of Algorithm~\ref{alg:MAW_global_pruning}, whereas \emph{Full-regularized time} denotes the wall-clock time obtained when regularization is applied during all pruning steps.    }
\label{tab:maw_pruning_summary}
\begin{tabular}{cccccc}
\hline
$\tolSVDfixed$ & $\mInit$ & First-stage limit &   $\nomega$ & Total time  (s) & Full-regularized time (s) \\
\hline
$10^{-4}$ & $54$  & $20$ & $10$ & $50.0$ & $  1172.0$ \\
$10^{-5}$ & $87$  & $18$ & $9$  & $44.0$ & $  7682.0$  \\
$10^{-6}$ & $135$ & $18$ & $10$ & $58.0$ & --- \\
$10^{-7}$ & $187$ & $19$ & $10$ & $89.0$ & --- \\
\hline
\end{tabular}
\end{table}

 Table~\ref{tab:maw_pruning_summary} summarizes the results obtained
with the four initial cubature rules. The most striking observation is
the weak dependence of the final adaptive rule on the initialization.
As the fixed-weight tolerance becomes more stringent, the size of the
initial rule increases markedly, from \(\mInit=54\) to
\(\mInit=187\). Nevertheless, the pruning algorithm consistently
reduces these rules to either \(\nomega=9\) or \(\nomega=10\)
integration points. Thus, the final support appears to be nearly
independent of the initial one. Moreover, the resulting adaptive rules
contain only about twice the theoretical lower bound,
\(\mLower=5\).

 An equally important observation concerns the behavior of the first pruning stage, namely the stage in which positivity enforcement is not explicitly activated.  As can be gleaned from the third column of Table~\ref{tab:maw_pruning_summary}, independently of the initial number of points, this stage systematically reduces the number of points  down to approximately $18$--$20$   before positivity violations begin to appear (Line~\ref{alg1:enforceREQUIRE} of Algorithm~\ref{alg:prune_step_optionB}). Thus, not only the final number of cubature points, but also the number of points at which positivity violations appear, seems to be essentially independent of the initial number of points.

 This behavior has a direct impact on computational efficiency (see last two columns of Table~\ref{tab:maw_pruning_summary}). For instance,   increasing the initial support from $54$ to $87$ points produces   no increase in wall-clock time (which remains below one minute in both cases\footnote{All computations
were performed in MATLAB on a Linux laptop,
equipped with a 13th-generation Intel Core i9-13900H processor and
16~GiB of RAM.}). This confirms the rationale behind Algorithm~\ref{alg:prune_step_optionB} and the proposed two-stage strategy: the computational cost is governed primarily not by the initial number of candidate points, but by the number of pruning iterations requiring activation of the graph-regularized positivity-enforcement stage.

In order to quantify the computational gains introduced by the proposed two-stage strategy, we compare the performance of the algorithm against a modified version in which the first pruning stage is artificially suppressed. This is achieved by temporarily disabling, for testing purposes, the conditional skip in Line~\ref{alg1:success3} of Algorithm~\ref{alg:prune_step_optionB} (namely, the line stating ``Successful pruning step, no need for regularization''). Under this modification, graph regularization is enforced throughout the entire pruning process, leading to a dramatic increase in computational cost.  In particular, the case $\mInit = 54$ increases from   $50$ seconds to approximately $20$ minutes, whereas the case $\mInit = 87$ increases from   $44$ seconds to more than two hours.

  \PROMPT{\input{PromptMeta1}}

 %


\subsubsection{Adaptive cubature rule and weight-field regression}
\label{subsec:num_char_adaptive_weights}

\begin{figure}[!ht]
\centering

\subfloat[Initial fixed-weight cubature rule ($\mInit=54$).]{
\includegraphics[width=0.48\textwidth]{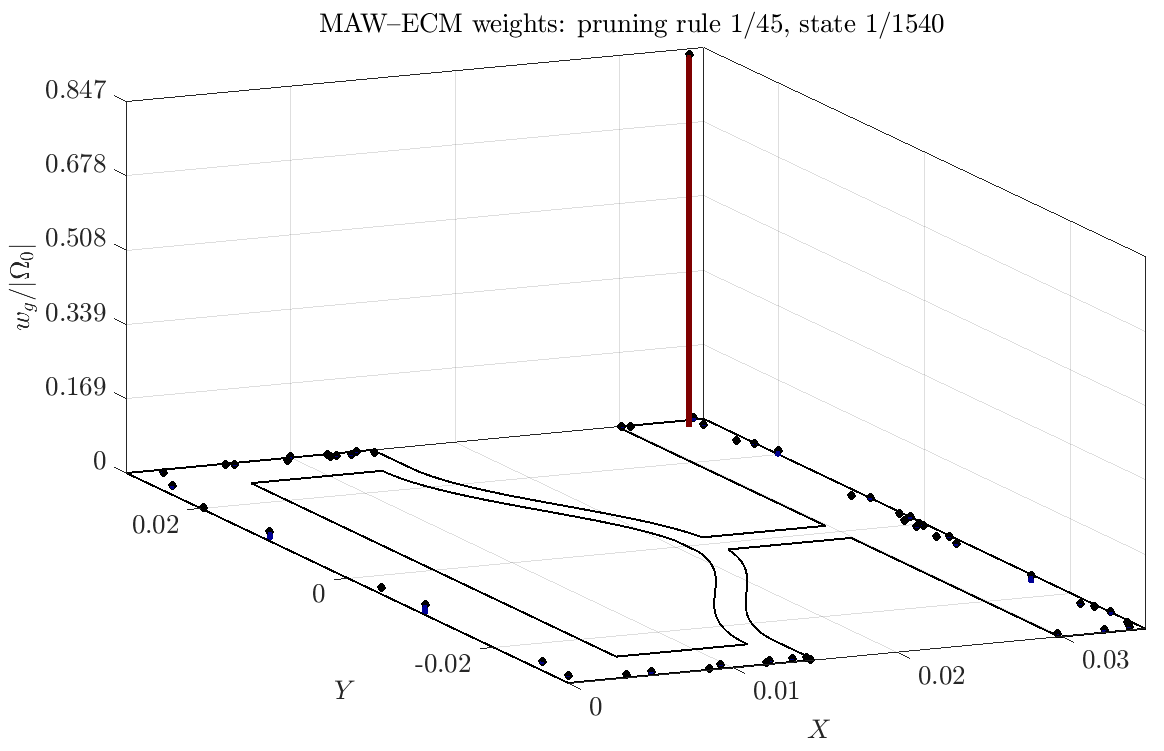} \label{fig:ecm54a}}
\hfill
\subfloat[Final adaptive cubature rule for $\qMone = \qMtwo = 0$ ($\nomega=10$).]{
\includegraphics[width=0.48\textwidth]{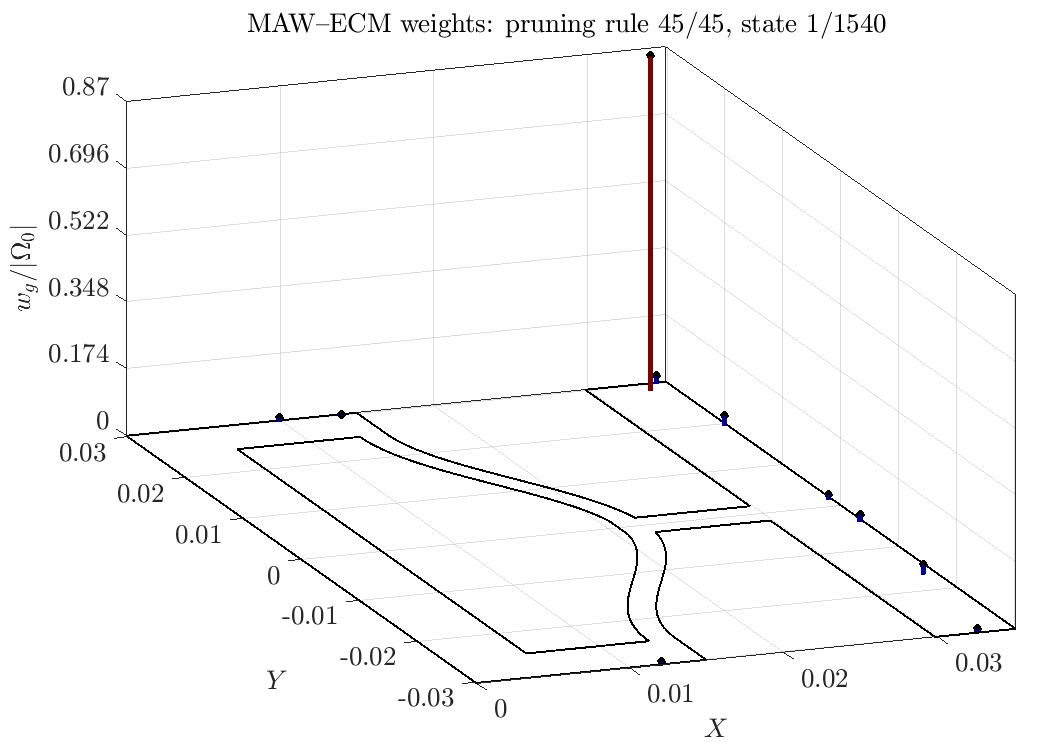}  \label{fig:ecm10a}}

\vspace{0.3cm}

\subfloat[Elements containing the $\mInit = 54$ integration points of the initial fixed-weight rule.]{
\includegraphics[width=0.4\textwidth]{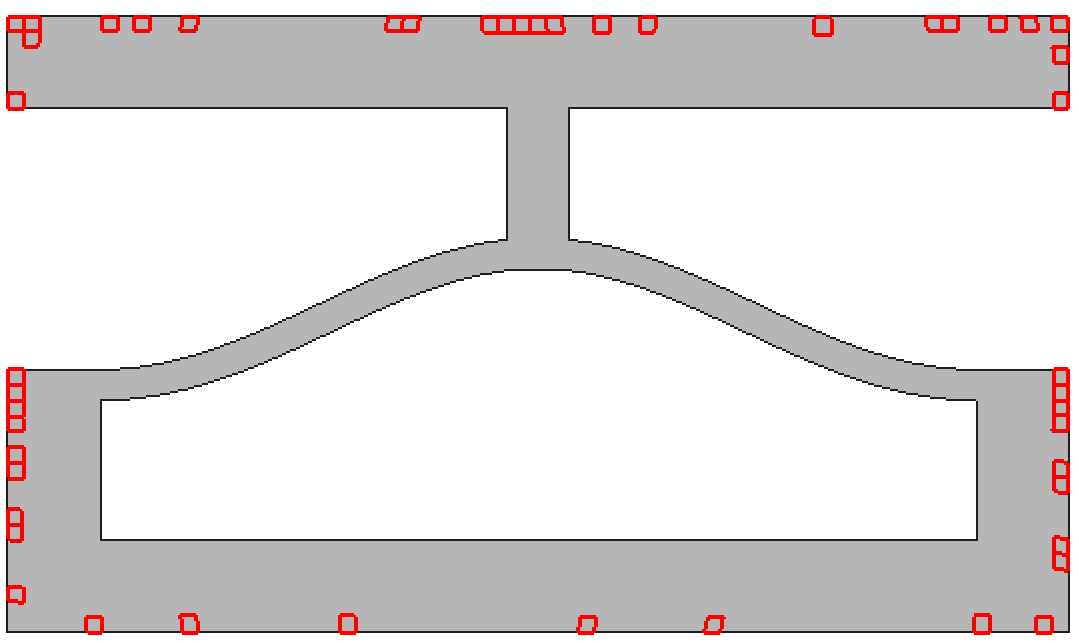} \label{fig:ecm54b}}
\hfill
\subfloat[Subset of elements containing the $\nomega  = 10$ integration points of the final adaptive rule.]{
\includegraphics[width=0.4\textwidth]{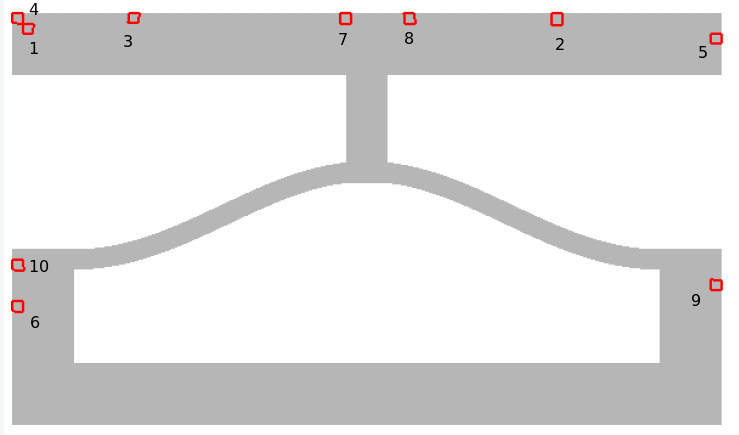}  \label{fig:ecm10b}}

\caption{
Comparison between the initial fixed-weight manifold ECM rule and the final adaptive cubature rule obtained after pruning. The upper row shows the associated cubature weights, whereas the lower row identifies the elements containing the selected Gauss points.  The adaptive rule shown in the right column corresponds to latent state $\qMone=\qMtwo=0$. Heights in the upper row are proportional to the corresponding cubature weights.
}
\label{fig:ecm54_vs_maw10}
\end{figure}

Having quantified the level of reduction afforded by MAW--ECM, we next examine the structure of the resulting adaptive rule in greater detail. In particular, we focus on the case $\mInit=54$, and analyze both the spatial distribution of the retained integration points and the associated distribution of weights.
Figure~\ref{fig:ecm54_vs_maw10} compares the initial fixed-weight rule   with the final adaptive rule containing $\nomega=10$ points (  for the particular state $\qMone=0$ and $\qMtwo=0$). The upper row displays the spatial location of the selected integration points and their associated cubature weights, whereas the lower row identifies the finite elements containing such points. Notice that the initial rule in Fig.~\ref{fig:ecm54a} is already highly non-uniform: one of the $54$ integration points carries approximately $85\%$ of the total volume. As discussed in Ref.~\cite{hernandez2024cecm}, this lack of uniformity is a characteristic feature of ECM-type cubature rules, reflecting the greedy, sparsity-promoting nature of the algorithm (which prioritizes error reduction over weight uniformity). It can be seen in Fig.~\ref{fig:ecm10a} that the adaptive pruning process preserves this dominant point (labelled as point~1), with its relative weight   increasing from $85\%$ to $87\%$ in the final rule.

\begin{figure}[!ht]
\centering
\includegraphics[width=\textwidth]{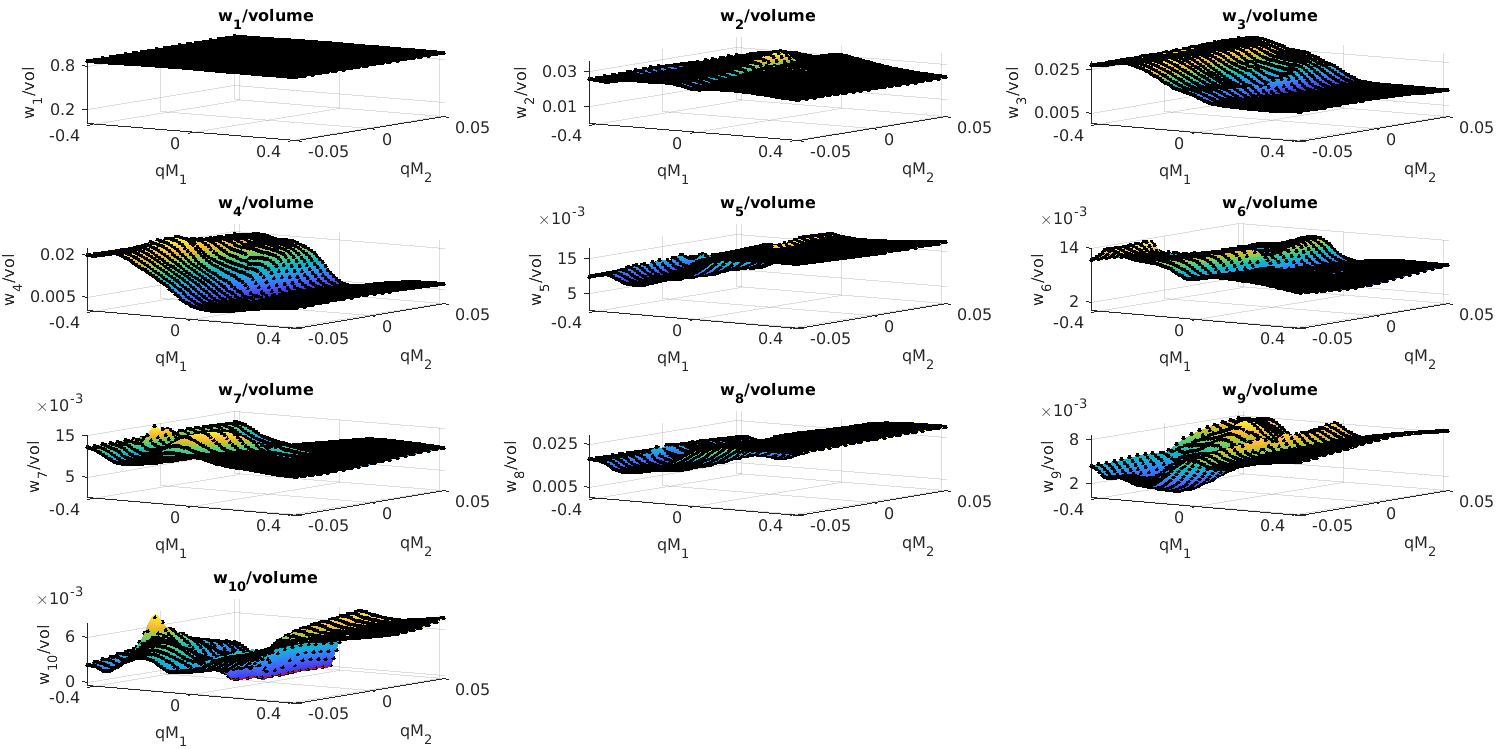}
\caption{RBF reconstruction of the $\nomega = 10$ adaptive weight fields $\omegag{g}=\omegag{g}(\qM)$, $g=1,\ldots,\nomega$ (normalized by the volume $|\Omega_0|$). Black dots denote the sampled manifold states used by the adaptive cubature procedure. The spatial location of the integration point associated with each weight field is shown in Fig.~\ref{fig:ecm10b}.}
\label{fig:weights_alpha1e4_rbf}
\end{figure}

 Once the support of the adaptive cubature rule has been identified, the next step consists of constructing continuous approximations of the associated adaptive weight fields, $\omegag{g}=\omegag{g}(\qM)$, $g=1,\ldots,\nomega$. Here, we use the same  anisotropic RBFs employed for the decoder, but with constant polynomial correction,   Tikhonov regularization parameter $\lambda=10^{-10}$ and characteristic length scales $\ellRBF_1 = 2.2\times10^{-2}$ and $\ellRBF_2 = 1.053\times10^{-2}$. Figure~\ref{fig:weights_alpha1e4_rbf} shows the resulting weight fields (normalized by the volume $|\Omega_0|$), together with the sampled manifold states (black dots).

The first weight field, $\omegag{1}$  remains nearly constant throughout the latent manifold and accounts for approximately $87\%$ of the total volume. This observation further reinforces the interpretation of this point as primarily responsible for enforcing the volume-preservation constraint. The nontrivial variations are therefore concentrated in the remaining nine weight fields. Their average magnitude ranges from approximately $3.0\%$ of the total volume for $\omegag{2}$ down to only $0.05\%$ for $\omegag{10}$. Furthermore, the amplitude of the oscillations tends to increase as the average value of the corresponding weight decreases. At the same time, the spatial structure of the corresponding fields becomes progressively richer as their average magnitude decreases.  Nevertheless, all weight fields remain smooth functions of the latent coordinates, with no evidence of abrupt localized variations.  The RBF reconstruction yields a relative interpolation error of $3.89\cdot10^{-9}$ over the complete dataset.

\begin{figure}[!ht]
\centering
\subfloat[$\eta = 0.25 $]{
\includegraphics[width=0.24\textwidth]{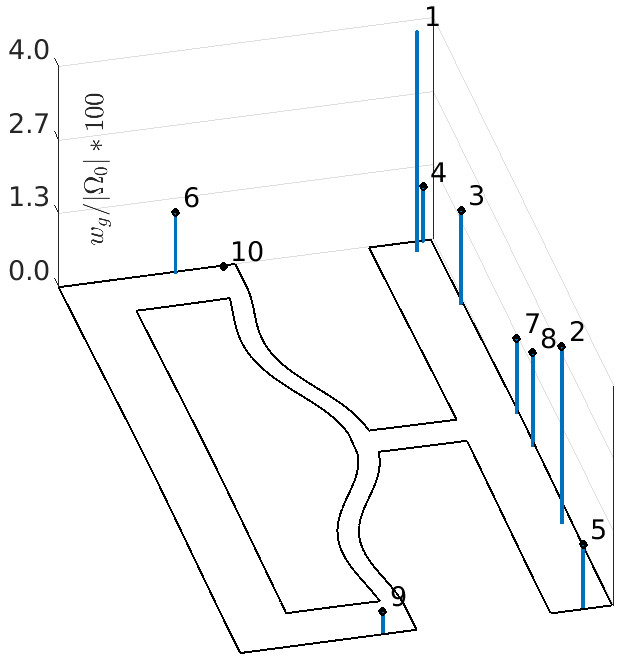}}
\hfill
\subfloat[$\eta = 0.5 $]{
\includegraphics[width=0.21\textwidth]{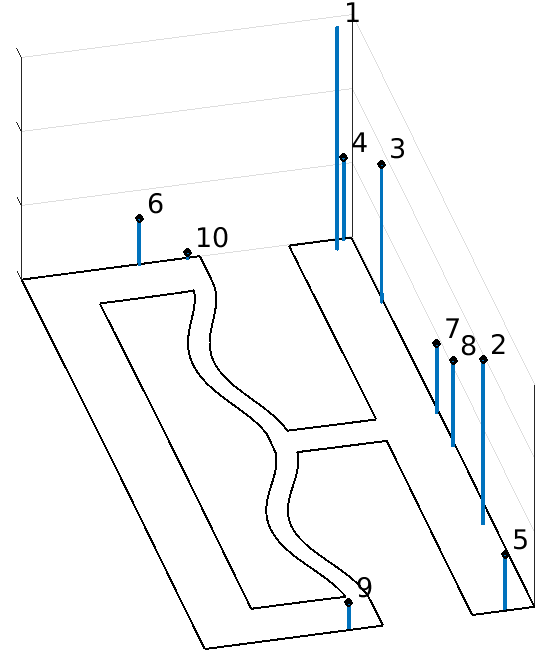}}
\hfill
\subfloat[$\eta = 0.75 $]{
\includegraphics[width=0.21\textwidth]{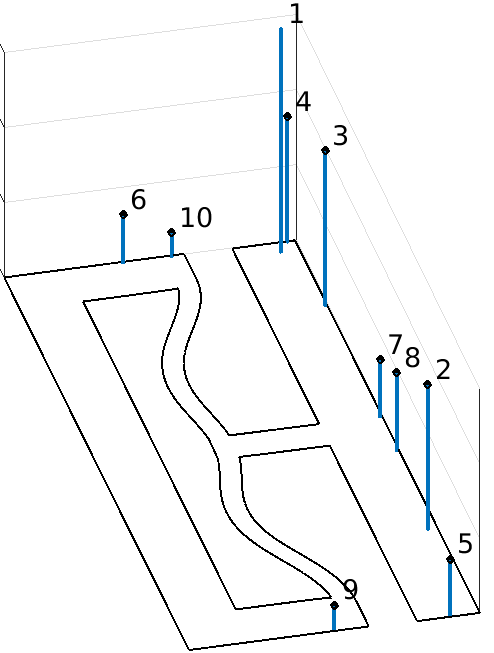}}
\hfill
\subfloat[$\eta = 1.0 $]{
\includegraphics[width=0.21\textwidth]{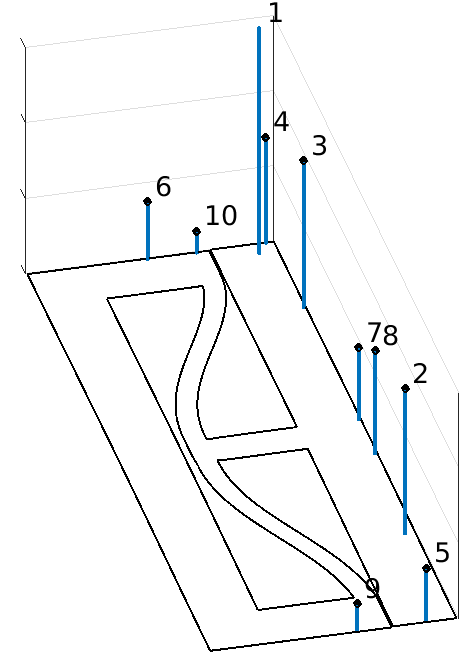}}
\caption{Evolution of the adaptive cubature weights along a purely compressive loading path ($\qMtwo=\Gmacroxy=0$, $\qMone=\Gmacrox<0$). The quantity $\eta=\eta(\qMone)$ denotes the compression level ($\eta = 0$ for the undeformed configuration, and $\eta = 1$ for maximum compression). Each subfigure shows the deformed configuration of the unit cell together with the location of the $\nomega=10$ integration points and their associated weights, which evolve according to the regressed fields shown in Fig.~\ref{fig:weights_alpha1e4_rbf}. The height of each bar represents the percentage of the total volume carried by the corresponding integration point.  To make the redistribution among the smaller weights visible, the plotted height is truncated at $4\%$ of the total volume.  Consequently, the dominant weight $\omegag{1}$ is displayed only partially; its actual value is approximately $88.37\%$, $87.83\%$, $87.21\%$, and $86.99\%$ for $\eta=0.25$, $0.5$, $0.75$, and $1.0$, respectively.}
\label{fig:adaptive_weights_compression_path}
\end{figure}

An alternative way of visualizing the   weight variability is to exploit the fact that the latent coordinates are approximately equal to the input macro-gradients ($\qMone\approx\Gmacrox$ and $\qMtwo\approx\Gmacroxy$). Thus, the adaptive weights can be interpreted directly as deformation-dependent quadrature weights.   Figure~\ref{fig:adaptive_weights_compression_path}
shows this interpretation along a purely compressive path, with
$\Gmacroxy=0$ and
$\Gmacrox<0$.  The deformed configuration of the unit cell is displayed together with the corresponding adaptive cubature weights.

\subsubsection{Effect of graph regularization}

\begin{figure}[!ht]
\centering
\subfloat[$\alphaG=0$.]{
\includegraphics[width=0.48\textwidth]{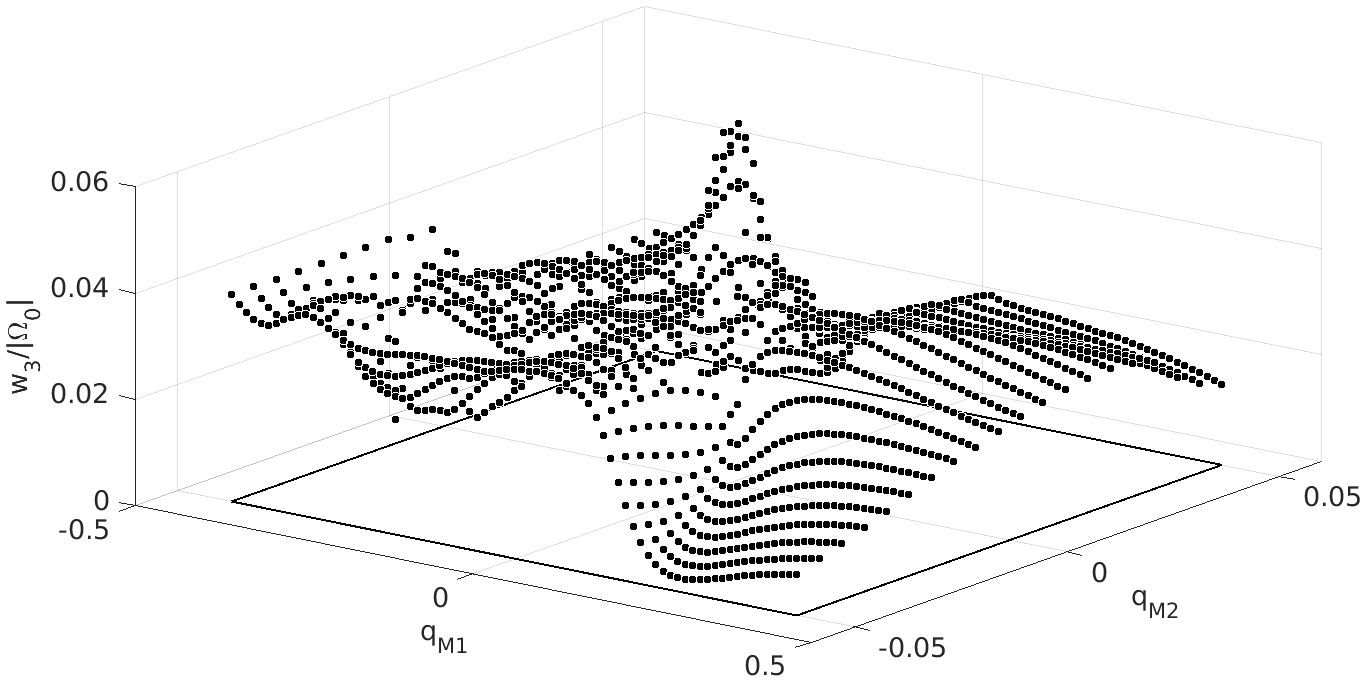} \label{fig:nor}}
\hfill
\subfloat[$\alphaG=10^{4}$.]{
\includegraphics[width=0.48\textwidth]{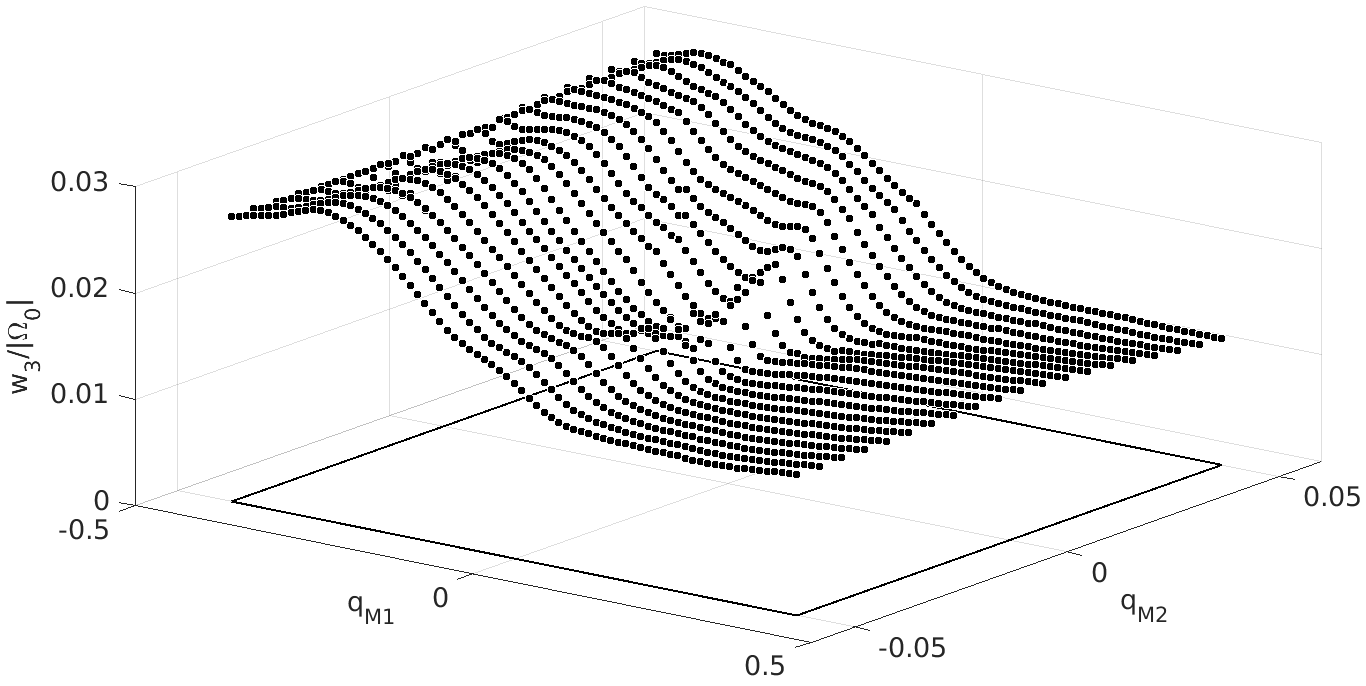} \label{fig:r}}
\caption{Sampled values of the third adaptive weight field obtained without graph regularization ($\alphaG=0$, left) and with graph regularization ($\alphaG=10^{4}$, right).}
\label{fig:w3_alpha_comparison}
\end{figure}

Let us now qualitatively examine the effect of the graph regularization parameter $\alphaG$ on the smoothness of the resulting weight fields.  To this end, Fig.~\ref{fig:w3_alpha_comparison} compares the sampled values of the third weight field for $\alphaG=10^{4}$ (the same field previously shown in Fig.~\ref{fig:weights_alpha1e4_rbf}) and for $\alphaG=0$, i.e., in the absence of regularization. The benefits of the regularization are  evident: for $\alphaG=0$, the sampled values clearly exhibit stronger local fluctuations and a less coherent evolution over the latent manifold.

\begin{table}[!ht]
\centering
\caption{Influence of graph regularization on the accuracy of the adaptive cubature rule. Errors correspond to averages over the two training trajectories.}
\label{tab:MAW_ECM_accuracy}
\begin{tabular}{lccc}
\toprule
Model & $\nECM$ & $\errdisp$ & $\errPK$ \\
\midrule
M-HROM, fixed-weight ECM & 54
& $7.87\cdot10^{-5}$
& $2.98\cdot10^{-4}$ \\

MAW--ECM ($\alphaG=10^{4}$) & 10
& $7.87\cdot10^{-5}$
& $7.21\cdot10^{-4}$ \\

MAW--ECM ($\alphaG=0$) & 10
& $7.87\cdot10^{-5}$
& $5.21\cdot10^{-3}$ \\
\bottomrule
\end{tabular}
\end{table}

 A more quantitative assessment of the role of graph regularization is obtained by comparing the online errors shown in Table~\ref{tab:MAW_ECM_accuracy}.  The first row  reproduces the results previously reported for the fixed-weight M-HROM in Table~\ref{tab:online_fixed_manifold}, and is included here solely to facilitate a direct comparison with the adaptive cubature rules.  The displacement error $\errdisp$ remains unchanged across all three cases because the displacement   is reconstructed directly from the decoder (since $\qM \approx \muvec$),  and, therefore does not depend on the particular cubature rule employed. As for the error $\errPK$ in predicting homogenized PK1 stresses, we see that  reducing the support from the original $\mInit=54$ integration points to only $\nomega=10$   points   leads to an error increase with respect to the original fixed-weight rule of  $7.21/2.98\approx 2.4$ when the regularization is active. In the absence of regularization, this error level rises up to $52.1/2.98\approx 18$.

\begin{figure}[!ht]
\centering
\includegraphics[width=0.7\textwidth]{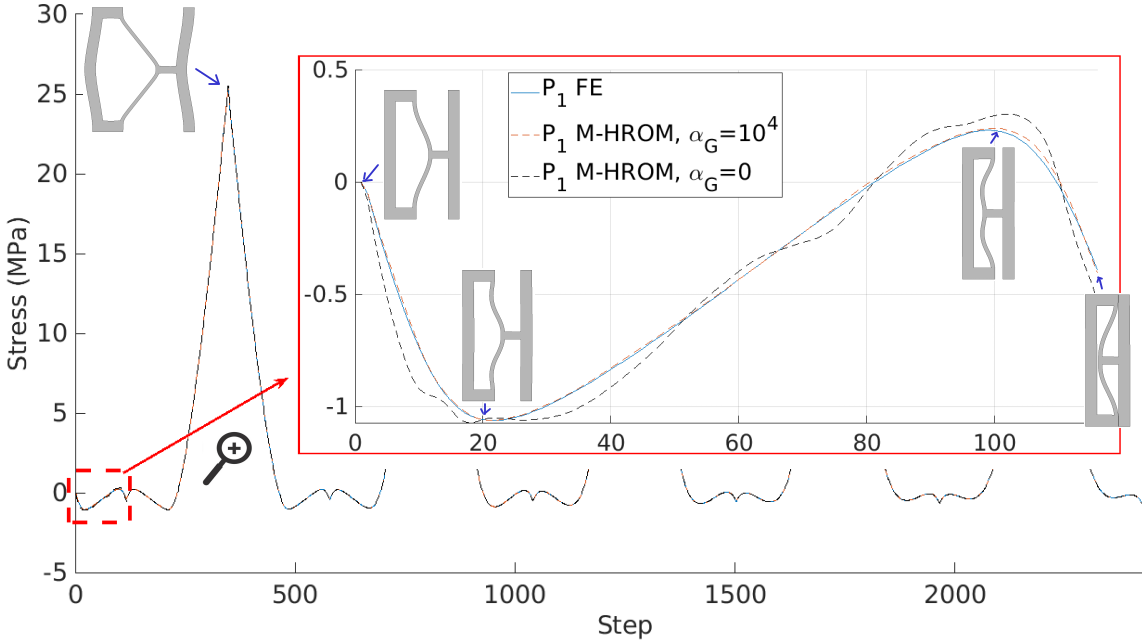}
\caption{Comparison of the first component of the homogenized PK1 stress predicted by the FE model, the regularized MAW--ECM rule ($\alphaG=10^{4}$), and the unregularized MAW--ECM rule ($\alphaG=0$) along the second zig-zag trajectory used for training. The inset shows a magnified view of the first portion of the trajectory (steps $1$--$116$), corresponding to the compression path previously shown in Fig.~\ref{fig:adaptive_weights_compression_path}. The deformed configurations illustrate the progressive buckling of the two internal beams and the associated negative-stiffness response of the unit cell.}
\label{fig:stress_regularization_comparison}
\end{figure}

To better understand the origin of the accuracy degradation caused by the absence of regularization, we compare in  Fig.~\ref{fig:stress_regularization_comparison}   the evolution of the first component of the homogenized PK1 stress along the second zig-zag trajectory used for training. The inset provides a magnified view of the first $116$ loading steps, corresponding to the purely compressive path previously shown in Fig.~\ref{fig:adaptive_weights_compression_path}. The accompanying deformed configurations illustrate the progressive buckling of the two internal beams responsible for the apparent negative-stiffness response of the cell.

It can be readily seen that the deterioration associated with the unregularized rule is not uniformly distributed. Rather, it is concentrated in the compressive regime\footnote{It is worth noting that the global Frobenius-type error measure employed in Table~\ref{tab:MAW_ECM_accuracy} tends to attenuate the impact of these localized discrepancies because the unit cell is approximately twenty-five times stiffer in tension than in compression. Interestingly, this asymmetry is also reflected in Fig.~\ref{fig:w3_alpha_comparison}, where the strongest departures from smooth behavior occur precisely in the compressive regime ($\qMone\equiv\Gmacrox<0$).
}, as can be appreciated  in the enlarged view. Although the unregularized M-HROM prediction ($\alphaG=0$) captures the overall trend of the response, it  exhibits noticeable local deviations from the FE solution. By contrast, the regularized prediction remains in close agreement with the reference response throughout the entire compression path.

In summary, the preceding results provide evidence that graph regularization does improve the smoothness of the adaptive weight fields and, in doing so, also the quality of the resulting stress predictions.


\section{Benchmark 2: continuum damage problem}
\label{sec:DamageProblem}
\subsection{Problem setting}
\label{sec:probsetting_2}
\begin{figure}[htbp]
\centering
\includegraphics[width=0.95\textwidth]{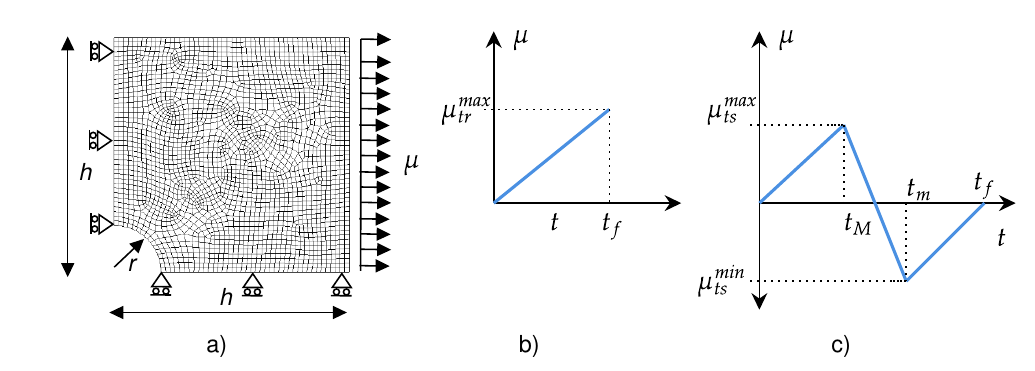}
\caption{
Damage benchmark and loading trajectories.
(a) Quarter model of a plate with a circular hole
(\(h=80\,\mathrm{mm}\) and \(r=20\,\mathrm{mm}\)).  Symmetry conditions are imposed along the left
and lower boundaries, while a uniform traction of amplitude \(\fload\)
is applied on the right boundary.
(b) Training trajectory used for the construction of the HROMs.
(c) Testing trajectory employed for model validation, involving loading,
unloading, and reloading stages.
}
\label{fig:damage_problem}
\end{figure}

Having examined the performance of the proposed methodology in a
problem involving geometric nonlinearities and elastic material
behavior, we next consider a problem in which the nonlinear response
originates from inelastic constitutive effects. The benchmark selected
for this purpose is the plane-strain deformation of a plate containing
a circular hole under small-strain conditions, governed by a simple
isotropic damage model with linear hardening.  Owing to symmetry, only one quarter of the structure is
modeled, as illustrated in Figure~\ref{fig:damage_problem}(a).   A uniformly distributed traction is applied along the right boundary;  the amplitude of such a  traction, denoted by \(\fload\), is selected as
the input parameter in the problem. Displacement boundary conditions are homogeneous, and thus $\dBOUND = \zero$ in Eq.~\eqref{eq:dvecreconstruction}; likewise $\Tlift$ in the same equation reduces to  a Boolean mapping matrix from the
unconstrained degrees of freedom to the full nodal space.

The output of interest is not a derived quantity, but the state variable itself, namely, the vector $\dred$ of
unconstrained nodal displacements. Unlike the Neo-Hookean-based benchmark considered in the previous
section, the constitutive memory embodied in
Eq.~\eqref{eq:damage_history_variable} prevents the solution from being
parameterized by the current value of the load factor \(\fload\)
alone. Consequently, the structural response must be recovered through
the solution of the reduced-order equilibrium equations.

Although the employed damage model
is classical, see e.g. Ref.
\cite{simo1987strain}, it is convenient to briefly recall its main
features, since they will later motivate several aspects of the
proposed input-informed decoder. In an isotropic damage model, the constitutive response at each Gauss point is driven by a single scalar
internal variable, \(\rDAMAGE\), whose evolution equation can be
integrated explicitly as follows:
\begin{equation}
\label{eq:damage_history_variable}
\rDAMAGE(t)
=
\max\left\{
\rDAMAGE(0),
\max_{0\leq s\leq t}
\sqrt{\strain(s)^{T}\Celas\strain(s)}
\right\},
\end{equation}
where  \(\strain\) and \(\Celas\) denote the infinitesimal strain vector
and the elasticity matrix, respectively.   The    damage variable  $\damage\in [0,1]$ is then recovered as $\damage =
1- {\qDAMAGE}/{\rDAMAGE}$, where
$ \dot{\qDAMAGE}
=
\hardeningD\,\dot{\rDAMAGE}$,
 \(\hardeningD >0\) being the hardening modulus.   The corresponding stress
vector is finally obtained as
\begin{equation}
\label{eq:damage_stress}
\stress
=
(1-\damage)\Celas\strain.
\end{equation}

The finite element mesh consists
of \(\nel = 2450\) quadratic quadrilateral elements (average size $2 \,\mathrm{mm}$) and \(9997\) nodes. The material  parameters are as follows: Young's modulus
\(\youngD=70\,\mathrm{GPa}\), Poisson ratio \(\poissonD=0.3\), damage
threshold \(\strengthD=70\,\mathrm{MPa}\), and hardening parameter
\(\hardeningD=0.01\).

Concerning training and generalization, we adopt a deliberately
parsimonious strategy: the  reduced-order models
considered in this work---namely the standard HROM, the manifold HROM,
and the proposed manifold HROM equipped with adaptive weights---are trained exclusively on a monotonic tensile loading path
(Figure~\ref{fig:damage_problem}(b)). No unloading, reloading, or
load-reversal information is supplied as data\footnote{ By stating that no unloading,
reloading, or load-reversal information is supplied as data, we mean
that such trajectories are not included among the training snapshots.
Naturally, this deliberate choice relies on prior knowledge of the
problem; for instance, the assumptions of small strains and isotropic damage imply that the model is symmetric with respect to changes in the sign of \(\fload\), so there is no need to include compression data.}. The
reduced-order models will then be challenged with a deformation history
that does contain unloading and reloading stages
(Figure~\ref{fig:damage_problem}(c)).

For the training stage, the loading parameter increases monotonically
from zero to \(\floadmax=70\,\mathrm{MPa}\) over a pseudo-time interval
\(t_f=1\) $s$, discretized using \(\nsnap=1000\)   increments.
The generalization test is performed using the non-monotonic loading
history shown in Figure~\ref{fig:damage_problem}(c), with
\(\mu^{\max}_{ts}=0.9\floadmax\),
\(\mu^{\min}_{ts}=-0.95\floadmax\), \(t_f=1\,\mathrm{s}\),
\(t_M=0.45t_f\), \(t_m=0.9t_f\), and \(1500\) increments.

\subsection{Standard HROM}



%



We begin the assessment by  constructing the standard   HROM (as pointed out in Section \ref{sec:meta_standard_hrom}, this model is recovered by setting   $\qM = \qrom$ and $\TAU = \qrom$
in the manifold HROM of  Section~\ref{sec:decoderENCODER}). As customary, the unconstrained displacement snapshots are
stored in the matrix
\(\Dsnap\in\mathbb{R}^{\nL\times\nsnap}\), where, in this case,
\(\nL=19832\). The reduced basis
\(\PhiROM\) is then sought as a linear combination of the columns of
\(\Dsnap\).  In contrast to the metamaterial example, however, \(\PhiROM\) is not
obtained directly from the singular value decomposition of
\(\Dsnap\). Instead, we follow the procedure introduced by the first author in
Ref.~\cite{hernandez2014high} for elastoplastic problems, and herein
adapted to the present damage setting, whereby the displacement vector
is decomposed into a linear reference contribution and an orthogonal
nonlinear correction:
\begin{equation}
\label{eq:decomNON}
 \dred = \dredLIN + \dredNON, \qquad \textrm{where}   \qquad  {\dredLIN}^T\dredNON = 0.
\end{equation}
The component \(\dredLIN\) represents the projection of the solution
onto the undamaged subspace, whereas \(\dredNON\) contains the
remaining information associated with damage. Separate reduced bases are then constructed for each contribution. This
requires distinguishing between the undamaged and damaged portions of
the training response. Accordingly, the displacement snapshots are partitioned
into two submatrices:
\begin{equation}
\label{eq:dsnapDECOMPOSITION}
 \Dsnap = \rowdos{\DsnapLIN}{\DsnapNON}
\end{equation}
  corresponding to
the undamaged regime, \(\DsnapLIN\), and the subsequent
nonlinear regime, \(\DsnapNON \in \RRn{\nL}{\nsnapNON}\) (for this case, \(\nsnapNON = 678\)). The snapshots in \(\DsnapLIN\) are compressed without truncation,
yielding the orthonormal basis \(\PhiROMlin \), which in the present case consists
of a single mode, since the undamaged response is linear and the applied load is parameterized by a single scalar. The remaining snapshots \(\DsnapNON\) are then
projected onto the orthogonal complement of \(\PhiROMlin\), producing
the residual matrix \(\DsnapNONorth\). A truncated singular value
decomposition is subsequently applied to \(\DsnapNONorth\), and the
corresponding left singular vectors define the basis
\(\PhiROMnon\). Hence,
\begin{equation}
\label{eq:bothdecom2s}
\dredLIN
=
\PhiROMlin \qMlin,
\qquad
\dredNON
\approx
\PhiROMnon \qNONall.
\end{equation}
where $\qMlin \in \Rn{}$ and $\qNONall \in \Rn{\rROM -1}$ are the corresponding linear and nonlinear amplitudes.

\begin{figure}[htbp]
\centering
\includegraphics[width=1\textwidth]{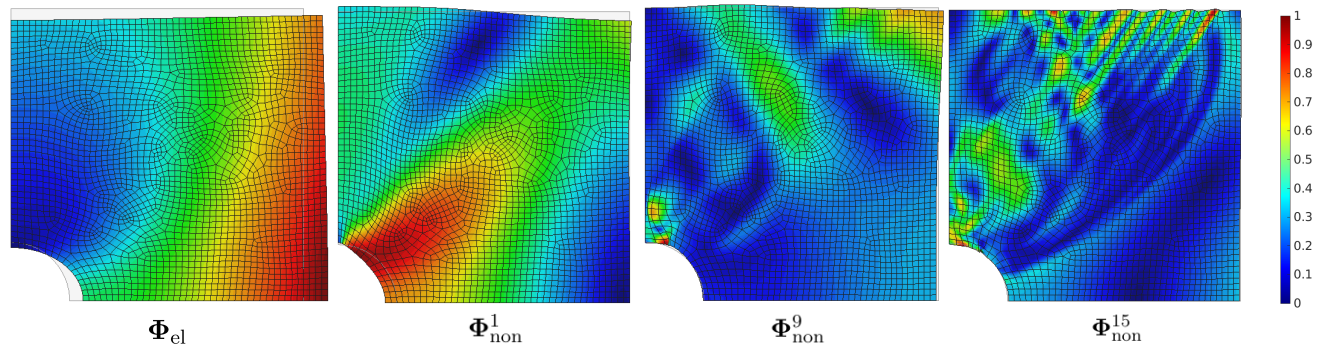}
\caption{
``Undamaged'' mode \(\PhiROMlin \in \Rn{\nL}\) and three representative nonlinear correction modes
from \(\PhiROMnon \in \RRn{\nL}{\rROM-1}\) ($\nL = 19832$, $\rROM = 16$), obtained from the displacement snapshots   generated
along the monotonic training trajectory shown in
Figure~\ref{fig:damage_problem}(b). The deformed configurations are displayed together with contour plots
of the Euclidean norm, normalized by its maximum value.
}
\label{fig:damage_modes}
\end{figure}

Using a relative truncation tolerance \(\tolSVDd=10^{-4}\), this
procedure yields 15 damage-correction modes. The resulting reduced
basis and associated generalized coordinates $\qrom$ are therefore
\begin{equation}
\label{eq:hromdamage}
\PhiROM
=
\left[
\PhiROMlin
\;
\PhiROMnon
\right], \qquad   \qrom = \coldos{\qMlin}{\qNONall},  \qquad \qrom = \PhiROM^T \dred,
\end{equation}
with a total dimension of \(\rROM=16\). Note that, by construction,
\begin{equation}
\label{eq:orthcond23}
{\PhiROMlin}^T\PhiROMnon=\bm{0}, \hspace{1cm}  {\PhiROMlin}^T\PhiROMlin=\ident,  \hspace{1cm} {\PhiROMnon}^T\PhiROMnon=\ident.
\end{equation}
By way of illustration, we show in Figure~\ref{fig:damage_modes}  the linear mode
\(\PhiROMlin\) together with three representative modes from
\(\PhiROMnon\).

\begin{figure}[!ht]
\centering
\includegraphics[width=0.95\textwidth]{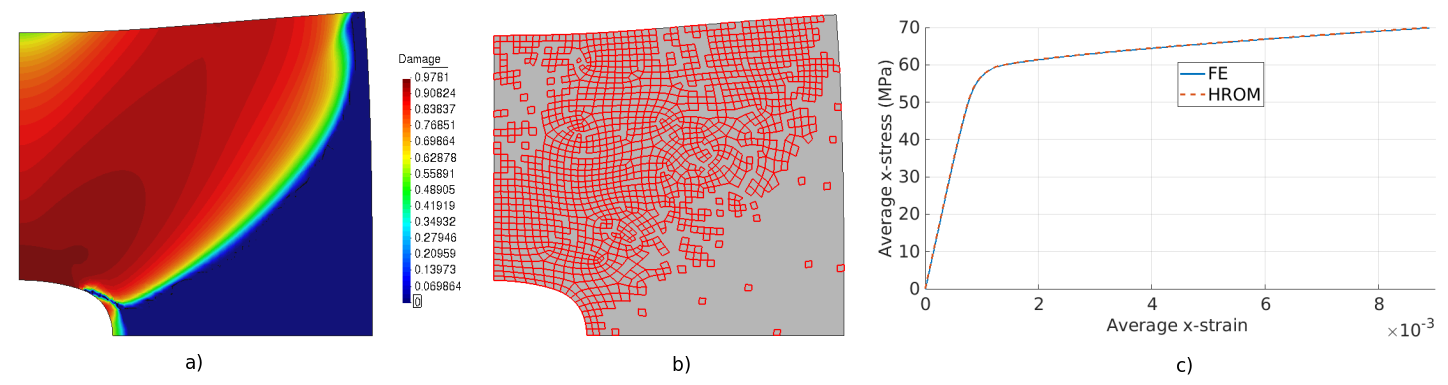}
\caption{Standard HROM for the damage benchmark.
(a) Damage field at the maximum training load.
(b) Elements containing the \(2370\) Gauss points selected by the
fixed ECM rule.
(c) FE and HROM average stress--strain responses along the monotonic
training trajectory.}
\label{fig:damage_hrom_std}
\end{figure}

We next construct the standard fixed-weight ECM rule following
Sections~\ref{sec:fixedweight_cubature}
and~\ref{sec:fixed_weight_cubature}. The only quantity entering the cubature training matrix is the set of projected internal-force densities; in this case, the matrix has dimensions $\ngp \times \rROM \nsnap$, where $\ngp = 2450\cdot9=22050$ and $\rROM \nsnap = 16\cdot 1000 = 16000$.
Using the SRSVD with \(\tolSVD=10^{-5}\), we get \(2369\)
left singular vectors. Augmenting this basis with the volume constraint
in Eq.~\refeq{eq:ECM_volume_constraint}, the ECM selects \(2370\)
Gauss points distributed over \(1140\) finite elements; see
Fig.~\ref{fig:damage_hrom_std}(b). Thus, we see that, despite the apparent simplicity of the problem, which is driven by the history of a single
input parameter, the resulting rule remains relatively dense, retaining 2370 of the $\ngp = 22050$ integration points, namely around $10.7\%$ of the original set. The explanation behind this limited sparsity may lie in the fact that the growth of damage produces
an evolving nonlinear region whose spatial support changes throughout
the loading process. From the standpoint of linear approximation, this
is analogous to the difficulties encountered in transport- and
propagation-dominated problems, where the movement or evolution of localized
features is known to result in a slow decay of singular values (see e.g. Refs.~\cite{barnett2022quadratic,barnett2023neural}). The contour plot of the damage variable \(\damage\) at the end of the
loading process and the spatial distribution of the resulting ECM
points, shown in
Figures~\ref{fig:damage_hrom_std}(a)--(b), corroborate this
interpretation, as most selected Gauss points cluster in the region
where damage is most severe.

The accuracy of the resulting HROM is assessed  using the relative
Frobenius displacement error \(\errdisp\), defined as in
Section~\ref{sec:fixedECM1}. On the training trajectory,
the standard HROM gives \(\errdisp=1.38\cdot10^{-4}\). The agreement
between the FE and HROM responses can be further appreciated  in Figure~\ref{fig:damage_hrom_std}(c), where we plot  the FE and HROM
average stress--strain curves in the loading direction. The average stress is
computed by dividing the applied force $\fload$ by the length of the loaded
boundary $h$, while the average strain is obtained from the predicted mean horizontal
displacement of the right boundary.    This graph also helps to identify the main stages of the response: an initial elastic branch is followed by a relatively sharp transition
associated with damage initiation and, subsequently, by a nonlinear
hardening regime.

\begin{figure}[htbp]
\centering
\includegraphics[width=\textwidth]{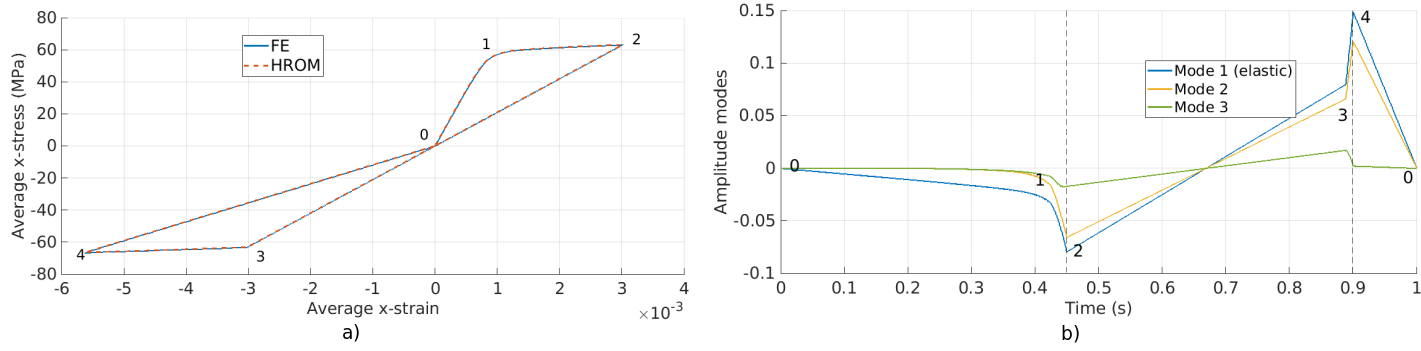}
\caption{Generalization test under cyclic loading. (a) Average stress--strain response obtained during the loading cycle.  Comparison between FE and standard HROM. The labels indicate the successive stages of the loading history: primary tensile loading (0--2), unloading from the tensile state (2--3), reverse loading into compression (3--4), and unloading from the compressive state (4--0).   (b) Evolution of the amplitudes  corresponding to  the first three modes $\PhiROM_1 = \PhiROMlin, \PhiROM_2$ and $\PhiROM_3$. }
\label{fig:damage_cycle}
\end{figure}

As for the  generalization test of
Figure~\ref{fig:damage_problem}(c), the relative displacement error remains of the same
order, namely \(\errdisp=1.28\cdot10^{-4}\). This result corroborates our initial assumption that, for the present
benchmark, training on the monotonic tensile branch is sufficient to
capture the subsequent unloading and reloading behavior,  including the (unseen) compression regime. The distinct stages of the cycle are illustrated in the average stress--strain graph of Figure~\ref{fig:damage_cycle}(a) and in the evolution of the amplitudes of the first three modes in Figure~\ref{fig:damage_cycle}(b).

\subsection{Manifold HROM}
 \label{sec:meta_manifold_representation2}



\subsubsection{Identification of latent variables}
\label{sec:id2}

We now turn to the manifold HROM. We begin with the construction of the input-informed
decoder (as discussed in
Section~\ref{sec:inputinformed}), for which we adopt the same master--slave structure introduced
in Eq.~\ref{eq:decoder_master_slave}, but adapted to the linear/nonlinear
decomposition of Eq.~\ref{eq:decomNON}. The latent coordinate
\(\qMlin\) and the corresponding basis matrix \(\PhiROMlin\) for the
linear component \(\dredLIN\) have already been identified. It remains to determine how many
additional latent coordinates are required to represent the nonlinear
component \(\dredNON\), how to identify them by exploiting the available
loading information, and, finally, how to construct the map from the
resulting master coordinates to the slave coordinates,
\(\qS=\Nslave(\qM)\).

The constitutive response at each Gauss point is governed by a single
scalar internal variable, \(\rDAMAGE\), defined in
Eq.~\eqref{eq:damage_history_variable}. This feature, together with the
fact that the evolution of these local variables along the training
trajectory is ultimately induced by the single scalar input
\(\fload\), suggests that one additional latent coordinate may suffice
to parameterize the nonlinear component \(\dredNON\). We denote this additional
coordinate by \(\qMnon\). In accordance with the hypothesis embodied in Eq.~\ref{eq:qMdefq}, we seek $\qMnon$ as a linear combination of the
nonlinear amplitudes $\qNONall = {\PhiROMnon}^T \dred$, that is:
\begin{equation}
\label{eq:TmNONdef}
 \qMnon = \TmNON \qNONall.
\end{equation}
Thus, the problem of identifying $\qMnon$  amounts to determining the row matrix $\TmNON \in \RRn{1}{(\rROM-1)}$. It is worth noting that, with the above choice, the transformation
matrix introduced in Eq.~\ref{eq:qMdefq} adopts the form
\begin{equation}
\Tm = \matcdos{1}{\zero}{0}{\TmNON}.
\end{equation}

The conditions used to determine \(\TmNON\) will not be imposed
directly on Eq.~\refeq{eq:TmNONdef}, but on its normalized counterpart,
obtained by dividing both sides of that equation by \(\qMlin\):
\begin{equation}
\label{eq:TmNONdef2}
\qMnonrel = \TmNON \qNONallrel,
\end{equation}
where
\begin{equation}
\label{eq:ratiosdEF}
\qMnonrel \defeq \dfrac{\qMnon}{\qMlin},
\qquad
\qNONallrel \defeq \dfrac{\qNONall}{\qMlin},
\qquad
\qMlin \neq 0.
\end{equation}

The rationale behind this normalization  is that  these normalized variables
remain unchanged during purely elastic loading or unloading and evolve
only when damage develops at one or more Gauss points. The
normalization therefore filters out the reversible dependence on the
loading amplitude and isolates changes associated with the evolution
of damage. This property follows directly from the nature of the employed damage model. Whenever the
internal variable \(\rDAMAGE\) in Eq.~\eqref{eq:damage_history_variable}, and consequently the damage variable
\(\damage\), remain frozen at all Gauss points, the structure responds
elastically with the fixed degraded stiffness implied by
Eq.~\ref{eq:damage_stress}. The reduced amplitudes are then proportional
to the input force $\fload$ and vanish at the same zero-load
state (as can be appreciated in Fig.~\ref{fig:damage_cycle}(b), between states 2 and 3); hence,  the ratios defined in Eq.~\ref{eq:ratiosdEF} remain
constant.

 Having clarified the role of the normalization, we now turn to the
determination of \(\TmNON\). We restrict attention to the portion of
the training trajectory over which damage evolves, represented by the
nonlinear snapshot block \(\DsnapNON\), see Eq.~\ref{eq:dsnapDECOMPOSITION}, whose columns correspond to the ordered load
levels
\begin{equation}
\label{eq:ordered}
\floadNON{1}
<
\floadNON{2}
<
\cdots
<
\floadNON{\nsnapNON},
\end{equation}
where \(\floadNON{1}\) marks the onset of damage, and
\(\floadNON{\nsnapNON} = \floadmax\) is the maximum load. The desired latent coordinate
\(\qMnonrel=\TmNON\qNONallrel\) must parameterize this one-dimensional
trajectory so that every component of \(\qNONallrel\) can be expressed
as a single-valued function of \(\qMnonrel\).  This requires \(\qMnonrel\) to be injective with respect to the applied
load \(\fload\), which in this case amounts to requiring strict monotonicity of
\(\qMnonrel\) with respect to \(\fload\).

  Monotonicity is the desired outcome, but not the quantity that we
enforce directly. Instead, we seek \(\TmNON\) as the coefficients of
the smoothest linear combination, on the grounds that oscillatory coordinates are incompatible with the
monotonicity requirement. The   resulting optimization problem leads to  the generalized eigenvalue
problem
\begin{equation}
\label{eq:damage_graph_eig}
(\QNONrel\KGsM\QNONrel^{T})\bm{v}_i
=
\lambda_i\,
(\QNONrel\QNONrel^{T})\bm{v}_i,
\qquad
i=1,\ldots,\rROM-1.
\end{equation}
Here, $\KGsM \in \RRn{\nsnapNON}{\nsnapNON}$ denotes the graph-Laplacian operator  associated with the sequence of loads
\(\{\floadNON{i}\}_{i=1}^{\nsnapNON}\), whereas $\QNONrel\in\RRn{(\rROM-1)}{\nsnapNON}$ is obtained by dividing each column of $\Qnon = {\PhiROMnon}^T \DsnapNON$ by the corresponding entries of $\Qlin = {\PhiROMlin}^T \DsnapNON$.
 The desired coefficients are finally chosen as
\begin{equation}
\label{eq:aplha}
\TmNON
=
\dfrac{\bm{v}_1^T}
{\normd{\PhiROMnon \bm{v}_1}},
\end{equation}
where \(\bm{v}_1\) is the eigenvector associated with the smallest
eigenvalue. The normalization by
\(\normd{\PhiROMnon \bm{v}_1}\) ensures that
\(\Am=\ident\) in Eq.~\ref{eq:decoder_master_slave}.

 \begin{remark}
The eigenvalue problem~\refeq{eq:damage_graph_eig} is a standard
construction in spectral graph theory; see, e.g.,
\cite{spielman2019spectral}. In the present problem, it arises from interpreting the ordered
load samples \(\{\floadNON{i}\}_{i=1}^{\nsnapNON}\) as the nodes of a
one-dimensional graph, and the rows of \(\QNONrel\) as scalar fields
defined on those nodes. Within this interpretation, the quadratic form
appearing in the left-hand side of Eq.~\refeq{eq:damage_graph_eig}  represents the
graph-Dirichlet energy (the same energy introduced when dealing with the regularization of the MAW--ECM in Section~\ref{sec:regoptprobl}), whereas the
quadratic form  in the right-hand side provides a
measure of its magnitude. The objective is therefore to find the
linear combination with minimum graph-Dirichlet energy among all
linear combinations of unit magnitude. This constrained minimization
problem takes the form of a Rayleigh quotient \cite{spielman2019spectral}, whose optimality
conditions are precisely Eq.~\refeq{eq:damage_graph_eig}. In the present one-dimensional setting, the graph-Laplacian operator $\KGsM$ can be interpreted as the stiffness (Laplacian) matrix  of a one-dimensional finite-element mesh whose nodes coincide with the ordered load samples.
\end{remark}

\begin{figure}[htbp]
\centering
 \begin{subfigure}[b]{0.48\textwidth}
\centering
\includegraphics[width=\textwidth]{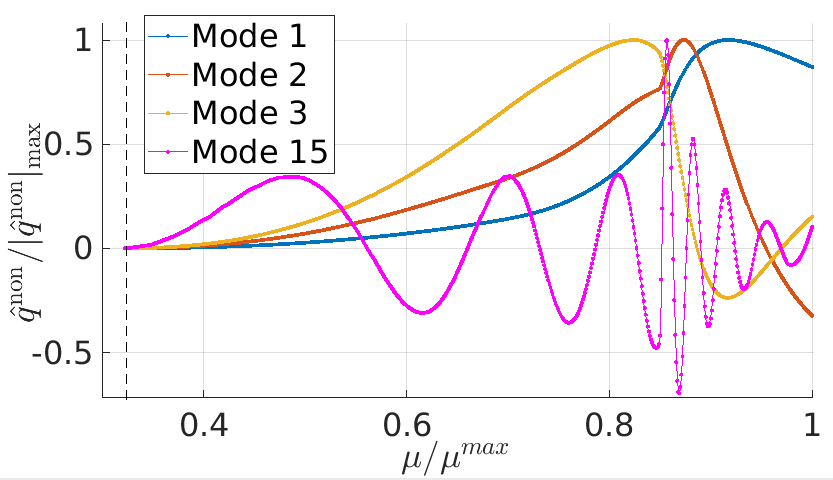}
\caption{}
\label{fig:damage_individual_coordinates}
\end{subfigure}
\hfill
\begin{subfigure}[b]{0.48\textwidth}
\centering
\includegraphics[width=\textwidth]{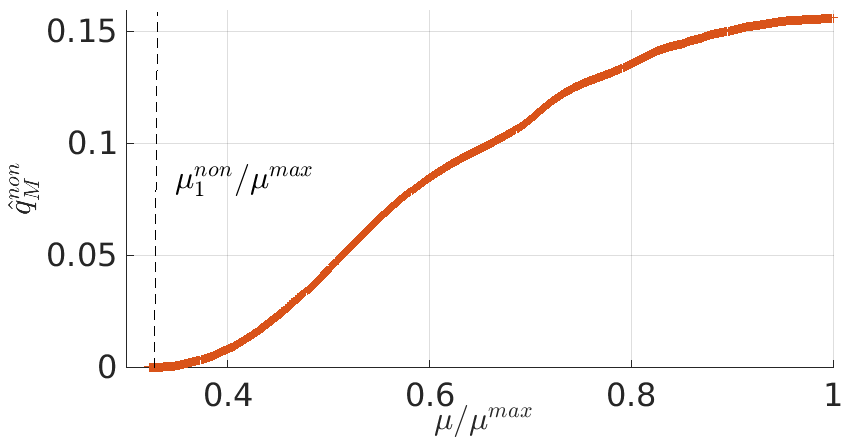}
\caption{}
\label{fig:damage_selected_coordinate}
\end{subfigure}
 \caption{Identification of the nonlinear latent coordinate for the
damage benchmark. The horizontal axis represents the normalized load
parameter \(\mu/\mu^{\max}\) for the nonlinear portion of the training
tensile test (see expression~\ref{eq:ordered}). (a) Components \(1\), \(2\), \(3\), and \(15\) of the
normalized nonlinear coordinate vector
\(\qNONallrel\in\Rn{\rROM-1}\), where \(\rROM=16\); see
Eq.~\refeq{eq:ratiosdEF}. (b) Normalized latent variable $\qMnonrel$ obtained as a linear combination of $\qNONallrel$ ($\qMnonrel = \TmNON \qNONallrel$). The matrix of coefficients in this linear combination, $\TmNON$,  is determined by solving    the generalized
eigenvalue problem \eqref{eq:damage_graph_eig}.  }
 \label{fig:damage_latent_coordinate}
\end{figure}

 Figure~\ref{fig:damage_latent_coordinate} provides an a posteriori
illustration of the identification. The   components of $\qNONallrel$ shown
in Fig.~\ref{fig:damage_individual_coordinates} fail to satisfy the
monotonicity requirement and are therefore unsuitable as scalar latent
coordinates (although not shown, the remaining components  also fail to satisfy the monotonicity requirement). By contrast, the linear combination associated with the
smallest generalized eigenvalue, shown in
Fig.~\ref{fig:damage_selected_coordinate}, exhibits the desired
monotone evolution. Thus, for the present benchmark, the adopted
smoothness criterion successfully yields an injective parametrization
of the normalized nonlinear trajectory.

\subsubsection{Nonlinear closure}
\label{sec:nonlinearclosure2}
Consistently with the normalization introduced above, we express the
closure map \(\Nslave=\Nslave(\qMlin,\qMnon)\) appearing in the decoder
of Eq.~\refeq{eq:decoder_master_slave} in the multiplicative form
\begin{equation}
\label{eq:nslave_map}
\Nslave(\qMlin,\qMnon)
=
\qMlin\,\NslaveONE(\qMnonrel).
\end{equation}
The map \(\NslaveONE=\NslaveONE(\qMnonrel)\) is constructed from the
input--output pairs  $\{\qMnonrel(\floadNON{j}), \gSnon(\floadNON{j})\}_{j=1}^{\nsnapNON}$, where $\gSnon(\floadNON{j})$ contains the normalized slave amplitudes:
\begin{equation}
\label{eq:gSnon_snap}
\gSnon(\floadNON{j})
=
\dfrac{{\PhiS}^{T}\dred(\floadNON{j})}
{\qMlin(\floadNON{j})}.
\end{equation}

The condition \(\qMlin=0\) occurs at zero applied load. Apart from the
initial state, it is encountered when the load passes through zero during
unloading. Since \(\qMnonrel\) remains constant throughout unloading, we
resolve the resulting indeterminacy by retaining its value from the
previous converged increment:
\begin{equation}
\label{eq:qMnonrel_zero}
\qMnonrel(t_{n+1})
=
\qMnonrel(t_n),
\qquad
\text{if}
\qquad
\qMlin(t_{n+1})=0,
\end{equation}
with \(\qMnonrel(0)=0\), where \(t_n\) denotes the pseudo-time of the
\(n\)-th load increment.

With the multiplicative closure in Eq.~\refeq{eq:nslave_map}, the
generic coefficient map in Eq.~\refeq{eq:tau_master_slave} specializes to
\begin{equation}
\label{eq:PhiTAU_damage}
 \TAU(\qM)
=
\begin{bmatrix}
\qMlin &
\qMnon &
\qMlin\,\NslaveONE(\qMnonrel)^T
\end{bmatrix}^T,
\end{equation}
 where \(\qMnonrel=\qMnon/\qMlin\) for \(\qMlin\neq0\), with its value
at \(\qMlin=0\) defined by Eq.~\refeq{eq:qMnonrel_zero}. Its Jacobian, see Eq.~\refeq{eq:JTAU}, follows directly from the
chain rule:
\begin{equation}
\label{eq:JTAU_damage}
\JTAU
=
\begin{bmatrix}
1 & 0\\[1mm]
0 & 1\\[1mm]
\NslaveONE(\qMnonrel)
-\qMnonrel\,\NslaveONEder(\qMnonrel)
&
\NslaveONEder(\qMnonrel)
\end{bmatrix}.
\end{equation}
Interestingly, the dependence on the master coordinates enters only
through the normalized coordinate \(\qMnonrel\). Hence, the tangent
basis \(\PhiD=\PhiROM\JTAU\), defined in
Eq.~\refeq{eq:PhiD_tau_generic} and used to project the equilibrium
residual, see Eq.~\refeq{eq:rDec_definition}, remains unchanged during elastic unloading and reloading at
frozen damage, since we have argued that \(\qMnonrel\) remains constant along
these branches.

\begin{figure}[htbp]
\centering
\includegraphics[width=0.7\textwidth]{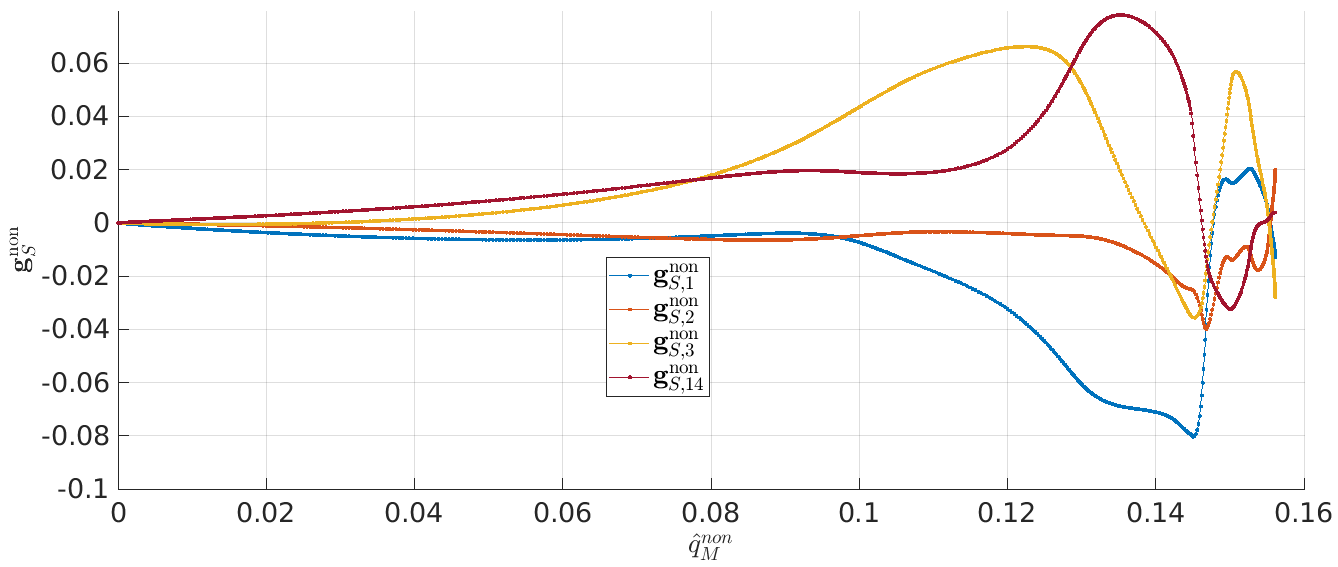}
\caption{Normalized slave amplitudes \(\gSnon\), defined in
Eq.~\refeq{eq:gSnon_snap}, as functions of the normalized latent
coordinate \(\qMnonrel\). This dataset is used for regressing the mapping component \(\NslaveONE = \NslaveONE(\qMnonrel)\) in Eq.~\ref{eq:nslave_map}.}
\label{fig:damage_decoder_data}
\end{figure}

 Figure~\ref{fig:damage_decoder_data} displays four representative
normalized slave amplitudes, \(\ggSnon{i}\) (\(i=1,2,3,14\)), as
functions of \(\qMnonrel\). They vary smoothly over most of the training
range, with localized oscillations and kinks appearing only near the
final, strongly damaged states. To prevent these irregularities from
producing oscillatory decoder derivatives in
Eqs.~\refeq{eq:JTAU_damage}, we fit \(\NslaveONE\) using \(100\)
uniformly distributed samples from the \(\nsnapNON=678\) snapshots and
cubic least-squares B-splines with \(90\) knots. The resulting smoothing
error is assessed in the following subsection, in the discussion of Table~\ref{tab:online_fixed_manifold_damage}. Outside the
training interval, the spline is extended by a quadratic Taylor
expansion at the corresponding endpoint.

\subsubsection{Fixed-weight ECM}
\label{sec:fixedECM2}

\begin{figure}[htbp]
\centering
\begin{subfigure}[t]{0.28\textwidth}
\centering
\includegraphics[width=\textwidth]{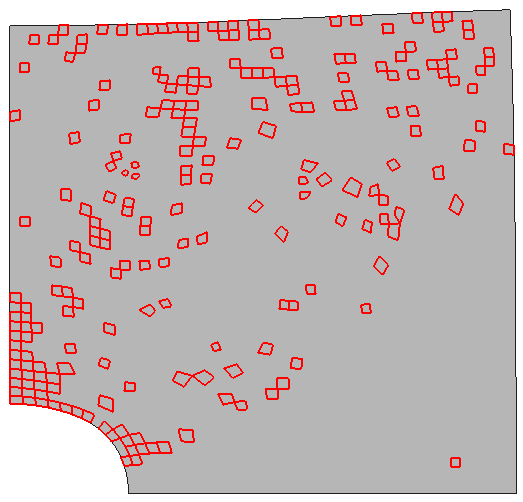}
\caption{Elements containing the selected Gauss points.}
\label{fig:ECM_FIXED2}
\end{subfigure}
\hfill
\begin{subfigure}[t]{0.62\textwidth}
\centering
\includegraphics[width=\textwidth]{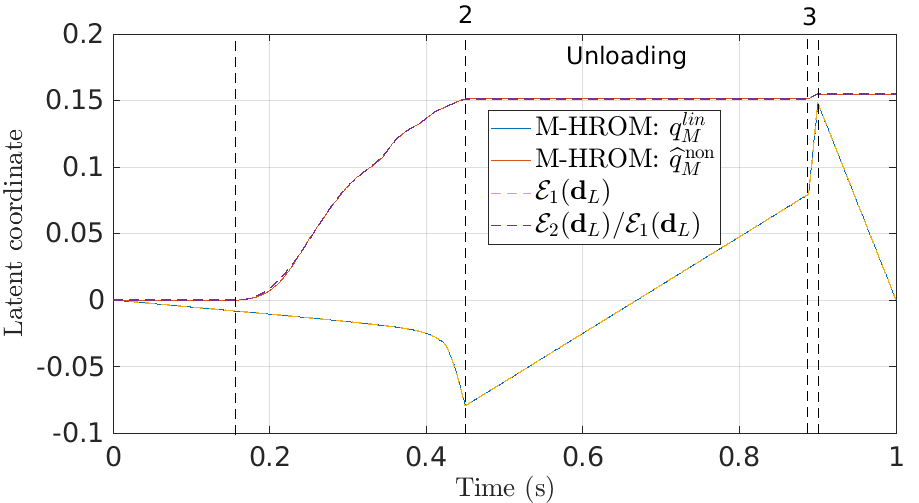}
\caption{Latent coordinates along the generalization trajectory.}
\label{fig:vinterna3}
\end{subfigure}
\caption{Fixed-weight manifold HROM for the damage benchmark.
(a) Elements containing the \(277\) Gauss points selected by the ECM.
(b) Evolution of \(\qMlin\) and \(\qMnonrel\) along the generalization
trajectory of Fig.~\ref{fig:damage_problem}(c), compared with the
coordinates obtained by applying the encoder in
Eq.~\refeq{eq:encoder_explicit} to the FE solutions.}
\label{fig:damage_mhrom_fixed_ecm}
\end{figure}

 Having defined the latent coordinates and the associated slave closure,
we now construct the fixed-weight ECM rule of
Section~\ref{sec:fixedweight_cubature}.  Since the present benchmark involves no volumetric output of interest,
we adopt the equilibrium-only case of Eq.~\refeq{eq:approxGEN}, for
which \(\rGENg{g}=\rDec{g}\) and \(\ncond=\rD=2\). Accordingly, the
training matrix \(\Agen\) of Eq.~\refeq{eq:AgenDEF} contains
only the projected internal-force densities and has dimensions
\(\Agen\in\mathbb{R}^{22050\times2000}\).  Using the SRSVD with
\(\tolSVD=10^{-5}\), and augmenting the resulting basis with the volume
constraint in Eq.~\refeq{eq:ECM_volume_constraint}, the ECM yields a
fixed rule with \(\mInit=277\) Gauss points distributed over \(213\)
finite elements, see Fig.~\ref{fig:damage_mhrom_fixed_ecm}(a).

Table~\ref{tab:online_fixed_manifold_damage} compares the online
complexity and accuracy of the resulting manifold HROM with those of
the standard HROM. The errors are the relative Frobenius displacement
errors over the monotonic training trajectory and the cyclic
generalization trajectory, respectively.

\begin{table}[!ht]
\centering
\caption{Online complexity and displacement accuracy of the fixed-weight
manifold HROM for the damage benchmark, together with the standard
HROM.}
\label{tab:online_fixed_manifold_damage}
\begin{tabular}{lccccc}
\toprule
Model
& Master coords.
& Slave coords.
& Integ. points
& \((\errdisp)^{\mathrm{train}}\)
& \((\errdisp)^{\mathrm{test}}\) \\
\midrule
HROM, ECM
& \(16\) & \(0\) & \(2370\)
& \(1.38\cdot10^{-4}\)
& \(1.28\cdot10^{-4}\) \\
M-HROM, fixed ECM
& \(2\) & \(14\) & \(277\)
& \(1.78\cdot10^{-3}\)
& \(1.81\cdot10^{-3}\) \\
\bottomrule
\end{tabular}
\end{table}

The manifold HROM reduces the number of independent coordinates from
\(\rROM=16\) to \(\rD=2\), and the number of integration points from
\(2370\) to \(277\). As in the first benchmark, the corresponding
compression factors are nearly identical, namely \(8\) and
approximately \(8.6\), respectively.

On the other hand, according to Table~\refeq{tab:online_fixed_manifold_damage}, the displacement errors increase by approximately one order of
magnitude relative to the standard HROM, but remain nearly identical
over the training and generalization trajectories. These results
indicate that the dominant additional error arises from the smoothing
introduced in the closure regression discussed in Section~\ref{sec:nonlinearclosure2}, rather
than from a loss of generalization over the unseen cyclic trajectory. Moreover, quadratic
Newton convergence was observed at every load increment.

Figure~\ref{fig:damage_mhrom_fixed_ecm}(b) further compares the latent
coordinates predicted by the manifold HROM along the generalization
trajectory with those obtained by applying the encoder in
Eq.~\refeq{eq:encoder_explicit} to the FE solutions. The linear
coordinates are virtually indistinguishable, whereas
\(\qMnonrel\) exhibits only minor discrepancies near the onset of
damage.


\begin{remark}
\label{remark:interpretation}
Although introduced solely to parameterize the nonlinear closure in Eq.~\refeq{eq:nslave_map},
\(\qMnonrel\) exhibits in
Figure~\ref{fig:damage_mhrom_fixed_ecm}(b) a behavior that admits a
compelling physical interpretation: as anticipated, it remains constant
during elastic unloading and reloading and, furthermore, increases
steadily as damage progresses. This is the same behavior exhibited by
the variable \(\rDAMAGE\) at the Gauss-point level; see
Eq.~\refeq{eq:damage_history_variable}. Thus, in seeking a maximally
compressed parametrization of the displacement manifold, we have
identified the structural-scale counterpart of the local internal
variable \(\rDAMAGE\).
\end{remark}

\subsection{Manifold-adaptive weight ECM }
\PROMPT{\input{Prompt_MAWD1}}
%
%

We now apply the MAW--ECM methodology summarized in
Box~\ref{box:MAW_ECM_summary} to the manifold HROM of the present
damage benchmark. Consistently with the normalized manifold
parametrization, we choose the adaptive weights to depend only on
\(\qMnonrel\):
\begin{equation}
\label{eq:normalizedTRAININGw}
\omegag{g}=\omegag{g}(\qMnonrel),
\qquad
g\in\Zomega.
\end{equation}


This choice assigns the same weight vector to all states of the virgin
elastic branch, for which \(\qMnonrel=0\) while \(\qMlin\) varies. The
adaptive rule must therefore reproduce exactly the projected internal
forces over this branch. To enforce this requirement, we partition the
training matrix into its elastic and nonlinear blocks,
$\Agen=\rowdos{\AgenLIN}{\AgenNON}$,
analogously to the displacement-snapshot decomposition in
Eq.~\refeq{eq:dsnapDECOMPOSITION}. Let \(\Uel\) denote a basis for the
column space of \(\AgenLIN\), which, owing to the linearity of the
response in this regime, has dimension \(\rD=2\). We include this basis
among the invariant constraints of Eq.~\refeq{eq:invariants1}:
\begin{equation}
\label{eq:Uinv_damage}
\Uinv=
\left[
\Uel,\,
\onevec
\right].
\end{equation}

The nonlinear block
\(\AgenNON\in\RRn{\ngp}{\nsnapNON\rD}\) supplies the state-dependent
conditions along the damaged portion of the training trajectory. In
summary, at each sampled value \(\qMnonrel(\floadNON{j})\), the local
system comprises the \(\ncond=\rD=2\) conditions associated with the
projected internal forces, the \(\rD=2\) elastic invariants in \(\Uel\),
and the volume constraint. Hence, the lower bound defined in
Eq.~\refeq{eq:ncondmin} becomes
\begin{equation}
\label{eq:ncond_damage}
\mLower = \ncond+\ninv = \rD+\rD+1 = 5.
\end{equation}
In other words, the desired adaptive rule satisfying the local
constraints must contain at least \(\mLower=5\) integration points.

\begin{remark}
A key consequence of the proposed parametrization for the adaptive weights  is that no additional
constraints are required to describe unloading. Indeed, since \(\qMnonrel\)
remains constant during unloading,   the Jacobian matrix
\(\JTAU\) of Eq.~\eqref{eq:JTAU_damage} remains constant as well. Furthermore, according to Eq.~\eqref{eq:PhiTAU_damage},
the decoder is linear with respect to \(\qMlin\) when
\(\qMnonrel\) is fixed. Hence, all states belonging to the same
unloading branch differ only through a scalar scaling factor and,  as a result, enforcing the local integration constraints at a given damage state
of the training set automatically guarantees the same level of
integration accuracy along any unloading trajectory emanating from
that state.
\end{remark}

\subsubsection{Influence of the fixed-weight initialization and pruning efficiency}
\label{sec:iniii}
We now proceed with the  assessment of the efficiency  of
the pruning process of
Algorithm~\ref{alg:MAW_global_pruning} for different initializations, namely three initial cubature
rules containing \(\mInit=97\), \(\mInit=277\), and
\(\mInit=534\) integration points, respectively.  These rules are obtained by varying the   truncation tolerance
\(\tolSVDfixed\) in Eq.~\refeq{eq:Ucub_definition}, applied to the
internal-force training matrix \(\Agen\). The initial rule displayed
previously in Figure~\ref{fig:ECM_FIXED2} corresponds to the
intermediate case, \(\tolSVDfixed=10^{-5}\), yielding \(\mInit=277\)
integration points.
The ordered samples
\(\{\qMnonrel(\floadNON{j})\}_{j=1}^{\nM}\) ($\nM = \nsnapNON$) define the nodes of
a one-dimensional latent mesh. Accordingly, the graph operator
\(\KGs\) entering the pruning algorithm is assembled from the standard
finite element discretization of the one-dimensional Laplacian using
linear shape functions, analogously to the construction of \(\KGsM\)
in Eq.~\refeq{eq:damage_graph_eig}.  The graph-regularization parameter is fixed at \(\alphaG=0.1\), and the
pruning algorithm performs \(\ntry=5\) trials per branch.

\begin{table}[ht!]
\centering
\caption{Influence of the fixed-weight initialization on the MAW--ECM
pruning process for different values of the SVD truncation tolerance
\(\tolSVDfixed\) employed in the construction of the
internal-force matrix \(\Agen\). Here, \(\mInit\) denotes the number
of points of the initial fixed-weight ECM rule, whereas
\emph{Unregularized pruning} reports the percentage of integration
points eliminated without activating graph regularization or explicit
positivity enforcement; \(\nomega\) denotes the final number of
integration points, and \emph{Total time} corresponds to the
wall-clock time of the complete MAW--ECM pruning algorithm.}
\label{tab:maw_pruning_summary_damage}
\begin{tabular}{ccccc}
\hline
\(\tolSVDfixed\) &
\(\mInit\) &
Unregularized pruning (\%) &
\(\nomega\) &
Total time (s) \\
\hline
$10^{-4}$ & $97$  & $80.4$ & $8$ & $46.1$ \\
$10^{-5}$ & $277$ & $87.0$ & $7$ & $131.0$ \\
$10^{-6}$ & $534$ & $96.4$ & $7$ & $63.9$ \\
\hline
\end{tabular}
\end{table}

Table~\ref{tab:maw_pruning_summary_damage} summarizes the results of
the three pruning tests. Despite the large variability in the size of
the initial fixed-weight rules, all three pruning processes converge
to adaptive cubature rules containing only
\(\nomega=7\)--\(8\) integration points, i.e., just two or three
points above the theoretical lower bound \(\mLower=5\) established in
Eq.~\eqref{eq:ncond_damage}.  This is
essentially the same behavior observed in the metamaterial benchmark
(Table~\ref{tab:maw_pruning_summary}) and corroborates that the
number of integration points is dictated primarily by the intrinsic dimension
of the solution manifold, rather than by the size of the initial
fixed-weight cubature rule.

 Table~\ref{tab:maw_pruning_summary_damage} also confirms another
trend observed in the metamaterial benchmark discussed in Section \ref{sec:infl_fw_ini}: graph regularization and
explicit positivity enforcement are required only during a relatively
small fraction of the pruning iterations. Indeed, in the three
studied cases, between \(80\%\) and \(96\%\) of the integration
points are eliminated without activating the positivity-enforcement
stage of Algorithm~\ref{alg:prune_step_optionB}. It is worth noting that the   tests also revealed an
intermittent activation of the regularization stage, with the pruning
algorithm   switching between regularized and unregularized
iterations before reaching the final rule. Such behavior was not
observed in the metamaterial benchmark, where the first activation of
the positivity-enforcement stage invariably led to a permanent
transition to the regularized regime.
\begin{figure}[t]
\centering
\includegraphics[width=0.9\textwidth]{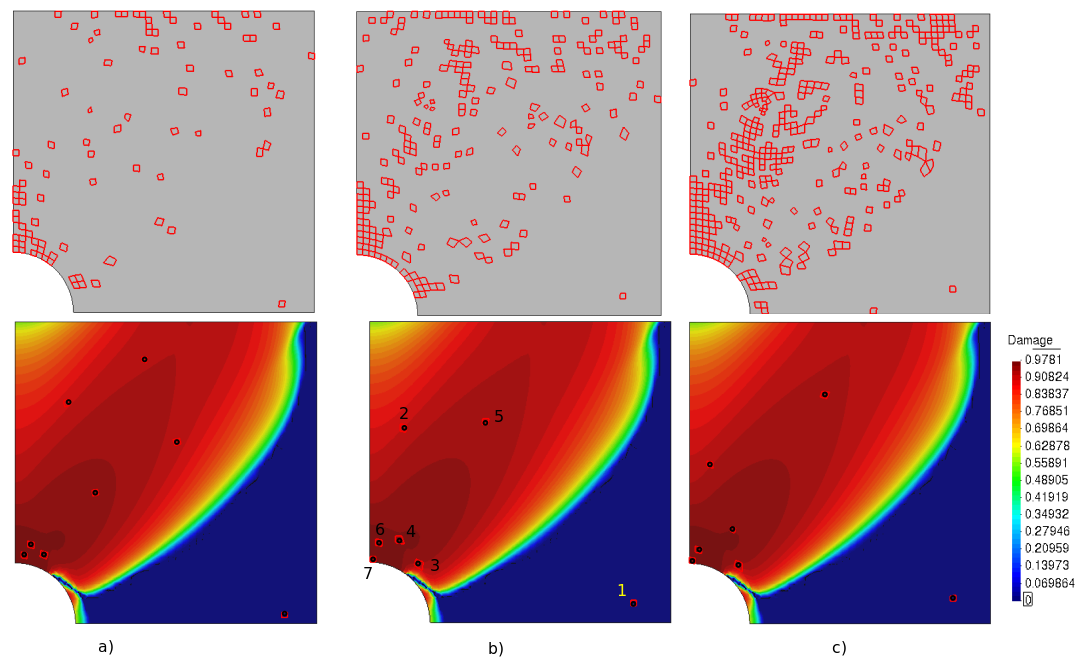}
\caption{Spatial distribution of the integration points associated with the
initial and final cubature rules corresponding to the three
fixed-weight initializations reported in
Table~\ref{tab:maw_pruning_summary_damage}. The upper row shows the
elements containing the integration points of the initial fixed-weight
rules. The lower row shows the retained integration points of the final
adaptive rules, superposed on the damage field \(\damage\) at the end
of the training trajectory. The three columns correspond to the
initializations \(\mInit=97\), \(\mInit=277\), and \(\mInit=534\),
respectively. The labels in panel (b) identify the seven retained
points used in the subsequent analysis.}
\label{fig:damage3cases}
\end{figure}

By way of illustration, Figure~\ref{fig:damage3cases} displays, in
the top row, the finite elements containing the integration points of
the initial fixed-weight rules considered in
Table~\ref{tab:maw_pruning_summary_damage};  the lower row, on the other hand, shows the
location of the integration points of the final sparsest rule,  superposed on the contour
plot of the  damage field  predicted by the FE model at the final step of the training trajectory.
It is interesting to note that the spatial distribution of the retained points exhibits, in all three
cases, a remarkably similar pattern: one point is invariably located
in the undamaged region (point \(1\) in
Fig.~\ref{fig:damage3cases}(b)); three or four points are placed in the vicinity of
the hole boundary, where damage nucleates and the stress
concentration is highest (points \(3\), \(4\), \(6\), and \(7\) in
Fig.~\ref{fig:damage3cases}(b)); and the other two or three points are located within the damaged region,
but away from the hole.

\begin{figure}[t]
\centering
\includegraphics[width=0.65\textwidth]{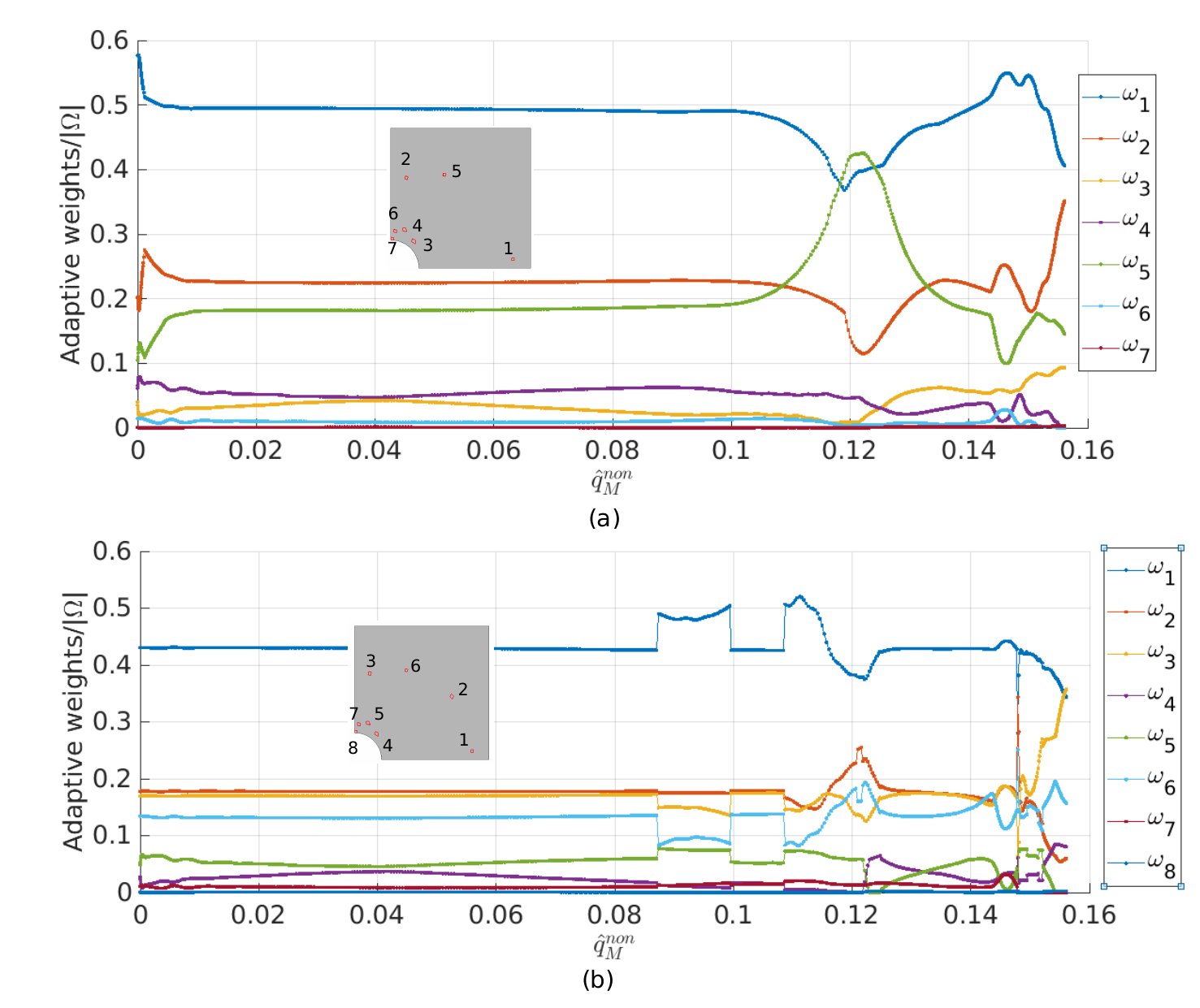}
\caption{Sampled adaptive weights obtained from the pruning algorithm
for the damage benchmark with \(\mInit=277\).  The insets show the retained supports and the
labels used to identify the individual weight samples.  (a) Graph-regularized pruning procedure.  (b) Pruning without regularization. }
\label{fig:damage_weights_samples}
\end{figure}

On the other hand, Figure~\ref{fig:damage_weights_samples}(a) displays the sampled
adaptive weights
\(\{(\Wad)_{g,j}\}_{j=1}^{\nM}\),
\(g=1,2,\ldots,\nomega\)
(see Eq.~\eqref{eq:pairs}) computed by the MAW--ECM for the case
\(\mInit=277\). The weights are represented as functions of the
corresponding latent states
\(\{\qMnonrel(\floadNON{j})\}_{j=1}^{\nsnapNON}\). For illustration purposes, the finite elements containing the integration points
of the corresponding adaptive rule, already shown in
Fig.~\ref{fig:damage3cases}(b), are reproduced as an inset.  Inspection
of these graphs  reveals a smooth dependence
on \(\qMnonrel\) throughout most of the sampled interval, with some irregularities appearing near \(\qMnonrel=0\) (which   may be
attributed to boundary effects inherent to the graph-regularization
procedure).   Nevertheless, such irregularities are relatively minor when compared
with the behavior of the graphs obtained when the  regularization stage of
Algorithm~\ref{alg:prune_step_optionB} is disabled, in
the same spirit as the comparison reported in
Fig.~\ref{fig:w3_alpha_comparison} for the metamaterial benchmark.
The corresponding results for the present problem are shown in
Fig.~\ref{fig:damage_weights_samples}(b).  We can see that  the sampled
weights exhibit abrupt variations  and a markedly less coherent
dependence on \(\qMnonrel\) than in the regularized case, confirming    that graph
regularization  does play a crucial role in generating weight distributions which are relatively smooth, and thus
amenable to subsequent regression.

\begin{figure}[t]
\centering
\includegraphics[width=0.8\textwidth]{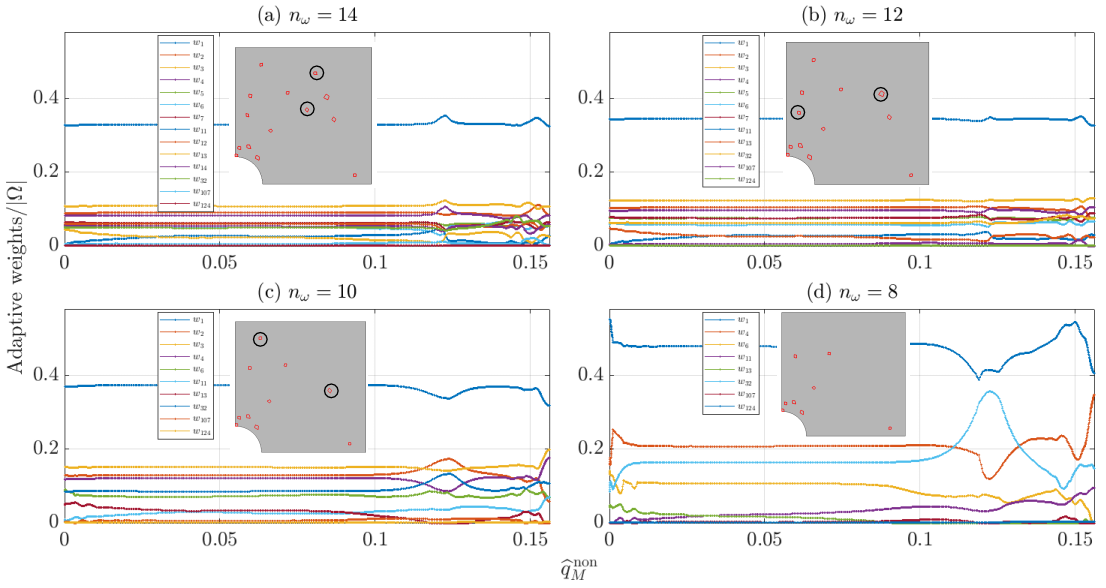}
\caption{Adaptive weights as a function of the normalized, nonlinear latent coordinate $\qMnonrel$, corresponding to four feasible
cubature rules extracted from different stages of the pruning process (containing
\(\nomega=14\), \(12\), \(10\), and \(8\) integration points).  The insets show the
elements containing the integration points of each rule. Integration
points enclosed by circles are those selected for removal in the
subsequent pruning stage.  }
 \label{fig:damage_pruning_path}
\end{figure}

The adaptive nature of the proposed cubature rules can be appreciated
more clearly in Fig.~\ref{fig:damage_pruning_path}, which displays the
sampled weights associated with four intermediate rules (with
\(\nomega=14\), \(12\), \(10\), and \(8\) integration points)
extracted along the pruning process that led to the \(\nomega=7\)
adaptive weight fields shown previously in
Fig.~\ref{fig:damage_weights_samples}(a).     To visualize the progressive evolution of the spatial distribution of
the integration points, each panel includes an inset showing the
corresponding integration rule, with circles marking the integration
points that will be removed in the subsequent pruning stage. We can
see that, as the rule is progressively pruned, the variability of the
weight curves increases.  Indeed, for \(\nomega=14\) and \(\nomega=12\), the weights exhibit
only mild variations along the latent manifold, whereas for
\(\nomega=10\) and, even more markedly, for \(\nomega=8\), they vary
more significantly with \(\qMnonrel\) to compensate for the discarded
integration points.

\begin{remark}
It should be recalled that all the rules shown in
Fig.~\ref{fig:damage_pruning_path} are feasible solutions of the
original optimization problem, i.e., they satisfy the local exactness
and positivity constraints. Consequently, if the optimal (i.e., the
sparsest) rule is deemed excessively irregular for subsequent
regression, one may deliberately step back along the pruning path and
select a denser, yet smoother, rule instead.
\end{remark}

%
%
%
%
%

\subsubsection{Online accuracy assessment}
\label{sec:online2}

We finally assess whether the substantial reduction in the number of integration points achieved by
MAW--ECM compromises the online accuracy of the manifold HROM. To this
end, we consider the adaptive rule shown in
Figure~\ref{fig:damage_weights_samples}(a), obtained by pruning the
\(\mInit=277\) fixed-weight rule down to only \(\nomega=7\) integration
points. The corresponding fixed-weight manifold HROM of
Section~\ref{sec:fixedECM2} (see Table~\ref{tab:damage_maw_accuracy}) is retained as reference, so that the
comparison isolates the effect of replacing the fixed cubature rule by
the adaptive one.

\begin{figure}[t]
\centering
\subfloat[Evolution of the latent coordinates.]{
\includegraphics[width=0.4\textwidth]
{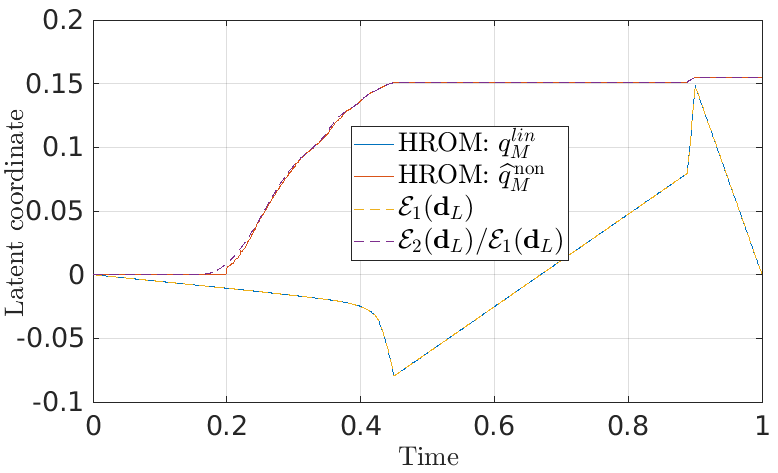}
\label{fig:damage_maw_latent}
}
\hfill
\subfloat[Average axial stress--strain response.]{
\includegraphics[width=0.35\textwidth]
{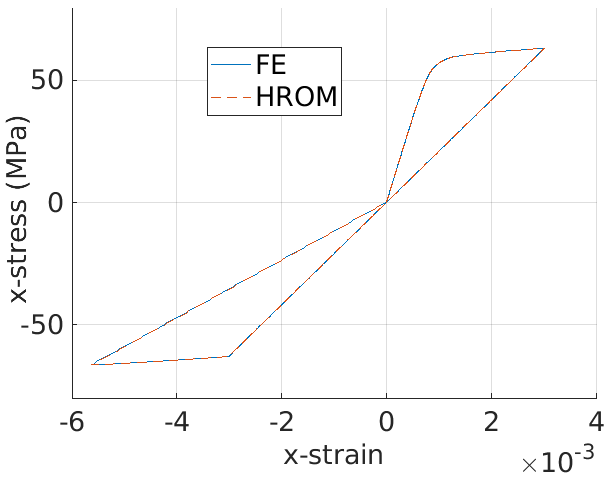}
\label{fig:damage_maw_stress_strain}
}
\caption{Online response of the manifold HROM with the adaptive rule of
\(\nomega=7\) integration points along the generalization trajectory of
Fig.~\ref{fig:damage_problem}(c).
(a) Evolution of \(\qMlin\) and \(\qMnonrel\), compared with the
coordinates obtained by applying the encoder in
Eq.~\refeq{eq:encoder_explicit} to the FE solutions (see Fig.~\ref{fig:damage_mhrom_fixed_ecm}(b) for the analogous
fixed-weight comparison).
(b) Average axial stress--strain response predicted by the FE and
manifold-HROM models (the same comparison carried out
in Fig.~\ref{fig:damage_cycle}(a) for the standard HROM).}
\label{fig:damage_maw_global}
\end{figure}
For online evaluation, the sampled adaptive weights
\(\omegag{g}=\omegag{g}(\qMnonrel)\) are represented using the same
spline framework employed for the decoder in
Section~\ref{sec:nonlinearclosure2}. The regressions use \(500\)
uniformly distributed samples selected from the \(\nsnapNON=678\)
available snapshots. Outside the sampled interval, each weight is held
constant at its nearest endpoint value, thereby preserving positivity
and the exact volume constraint inherited from the adaptive cubature
construction.

\begin{table}[t]
\centering
\caption{Relative displacement errors for the fixed- and
adaptive-weight manifold HROMs.}
\label{tab:damage_maw_accuracy}
\begin{tabular}{lccc}
\hline
Model
&
Number of points
&
\((\errdisp)^{\mathrm{train}}\)
&
\((\errdisp)^{\mathrm{test}}\)
\\
\hline
Fixed-weight manifold HROM
&
\(277\)
&
\(1.78\cdot10^{-3}\)
&
\(1.81\cdot10^{-3}\)
\\
Manifold HROM with MAW--ECM
&
\(7\)
&
\(1.8125\cdot10^{-3}\)
&
\(1.6538\cdot10^{-3}\)
\\
\hline
\end{tabular}
\end{table}

The resulting displacement errors over the training and generalization
trajectories are reported in Table~\ref{tab:damage_maw_accuracy}.
Inspection of this table shows that the adaptive rule essentially
preserves the accuracy of the fixed-weight manifold HROM: the small
variations observed in the training and test errors are negligible
relative to the approximation error already introduced by the decoder.

\begin{figure}[t]
\centering
\includegraphics[width=0.7\textwidth]{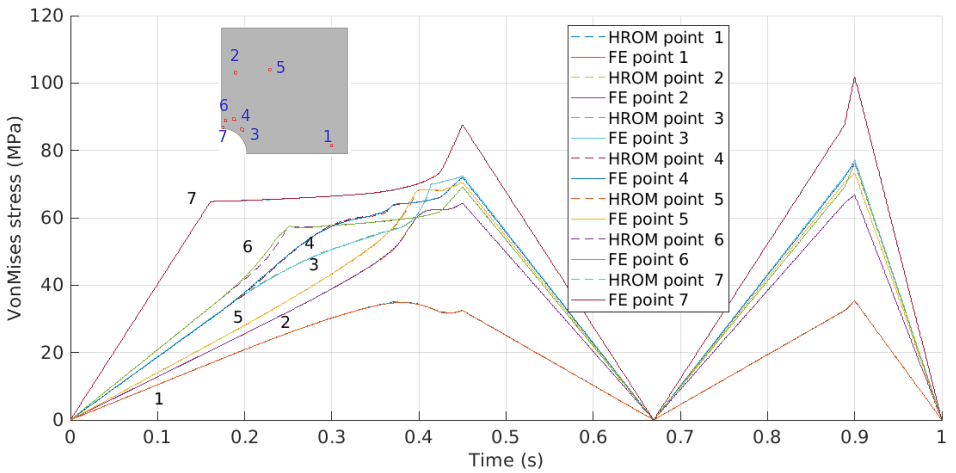}
\caption{Evolution of the von Mises equivalent stress  during the generalization test at the seven
Gauss points selected by MAW--ECM.
   The inset identifies the
location and numbering of the selected points.}
\label{fig:damage_maw_local_stresses}
\end{figure}


A more detailed assessment is provided in
Figure~\ref{fig:damage_maw_global}, which compares the latent
coordinates (the same comparison carried out in
Figure~\ref{fig:damage_mhrom_fixed_ecm}(b) for the fixed-weight
case) and the macroscopic stress--strain response along the
generalization trajectory. As in the fixed-weight case, \(\qMlin\) in
Figure~\ref{fig:damage_maw_latent} is virtually indistinguishable from
its encoded FE counterpart. The normalized nonlinear coordinate
\(\qMnonrel\), by contrast, exhibits a small but visible departure from
its encoded FE counterpart around \(t=0.2\)~s, at the onset of damage.
However, once damage evolution is fully established, the manifold-HROM
curve rejoins the reference FE curve, and the two solutions again
become indistinguishable.

The origin of the above-mentioned deviation appears to lie in the
weight regression: it occurs at the onset of damage, that is, when
\(\qMnonrel\approx0\), precisely where the sampled weights in
Figure~\ref{fig:damage_weights_samples}(a) exhibit their largest local
irregularities. Nevertheless, this local deviation has practically
no influence on the global displacement error, as confirmed by
Table~\ref{tab:damage_maw_accuracy}, nor on the macroscopic
stress--strain response shown in
Figure~\ref{fig:damage_maw_stress_strain}, where the complete
macroscopic cycle---including primary loading, unloading, reverse
loading in compression, and final unloading---is reproduced with no
visible discrepancy.

 Finally, Figure~\ref{fig:damage_maw_local_stresses} compares the von
Mises stress histories predicted by the FE model and the adaptive
manifold HROM at the seven selected Gauss points along the
generalization trajectory.  With the
exception of a small discrepancy at point \(6\) during the initial
damage transition, approximately over
\(0.2\lesssim t\lesssim0.3\),  the FE and manifold-HROM stress histories practically coincide.  The agreement is maintained not only
during the tensile branch represented in the training data, but also
throughout unloading, reverse loading in compression, and the final
return to the unloaded state.  Particularly noteworthy is point \(7\),
located in the vicinity of the hole and subjected to the largest stress
concentration.  Despite exhibiting the highest stresses and the
largest variation over the cycle, its complete history is reproduced
at plotting resolution.

\section{Concluding remarks}
%

The central premise examined in this work is that the conventional use
of fixed weights in sampling-and-weighting hyperreduction becomes
unnecessarily restrictive once the state is represented on a nonlinear
manifold. The two benchmarks of Sections~\ref{sec:metamaterial} and
\ref{sec:DamageProblem} support this premise. In both cases, replacing
the linear ROM by a nonlinear-manifold representation reduces the
number of integration points by a factor remarkably close to the
reduction in the number of independent state coordinates. In the
metamaterial problem, the latter decreases from \(52\) to \(2\), a
factor of \(26\), while the corresponding fixed cubature rule decreases
from \(1340\) to \(54\) Gauss points, a factor of approximately \(25\).
In the damage problem, the analogous factors are \(8\) and approximately
\(8.6\), respectively.
Allowing these weights to adapt to the latent state then provides a
further substantial reduction: MAW--ECM eliminates approximately
\(80\%\) of the remaining points in the metamaterial benchmark, while
in the damage problem the reduction exceeds \(97\%\), from \(277\)
points to only \(7\)--\(8\), with essentially no loss of accuracy
relative to the corresponding fixed-weight manifold HROM.

Of particular relevance is the weak dependence of these final
cardinalities on the fixed ECM rule used for initialization.  Initial supports
varying from \(54\) to \(187\) points in the first benchmark, and from
\(97\) to \(534\) points in the second, ultimately collapse to
approximately the same number of retained points. Moreover, the final
rules lie only a few points above the lower bounds imposed by the local
integration conditions.  These observations suggest that, once the weights are allowed to adapt,
the required support is controlled primarily by the dimension of the
local constraint space---which in the present examples is itself tied
to the intrinsic dimension of the solution manifold---rather than by
the dimension of the global span of all sampled integrands. This is perhaps the most significant practical
consequence of replacing global fixed weights by manifold-adaptive ones.

Adaptivity, however, is useful only if the resulting weight fields can
be represented reliably online. The numerical experiments confirm the role played by the proposed
graph regularization: without it, the sampled weights develop local
irregularities that substantially deteriorate the accuracy of their
regression. On the other hand, the strategy of performing the
redistribution locally whenever positivity is preserved and activating
the graph-coupled stage only when positivity enforcement has proved
successful is efficient. The experiments show that most points can be removed
in this inexpensive regime. In the metamaterial benchmark, for example,
this strategy reduces pruning times    from more than two hours
to less than one minute in one of the studied cases.

The present study is nevertheless exploratory in several respects.
Both examples involve only two independent latent coordinates. Extension
to substantially higher-dimensional parameter and latent spaces will
require revisiting the globally coupled graph operations entering the
regularized redistribution. Local graph constructions, partitions, or
clusters of the latent domain appear natural alternatives. Similarly,
the graph-based identification used for the history coordinate in
the damage problem has here been demonstrated only for an ordered
one-dimensional damage trajectory; its extension to genuinely
multidimensional loading histories remains to be investigated.

A second issue concerns the continuous representation of the adaptive
weights after pruning. Positivity and the invariant integration conditions are enforced at the
sampled latent states, but generic regression and extrapolation do not
preserve these properties automatically. Constraint-preserving
regression of the weight fields, including their extrapolation outside
the training domain, is therefore an immediate subject for further
work. The consistent linearization derived in
~\ref{app:MAW_linearization} also contains an additional tangent
contribution arising from the derivatives of the adaptive weights.
Although no convergence difficulties were encountered in the numerical
experiments reported here, the implications of this term for  robustness deserve a more
systematic investigation. Finally, extensions to dynamics and other
problem settings will require corresponding
invariant constraints to be identified with care. In particular, if adaptive cubature is applied to inertial terms,
preservation of the appropriate mass, first-moment, and inertia
properties should be examined so that state-dependent reweighting does
not introduce spurious inertial effects.

\section*{Acknowledgements}

Joaquín A. Hernández and S. Ares de Parga acknowledge  the support of Grant PID2024-158878OB-C21 funded by MICIU/AEI /10.13039/501100011033 and by ERDF/EU.

S. Ares de Parga also acknowledges partial support from the Department of Aeronautics and Astronautics at Stanford University.

\appendix
 \section{Consistent linearization of the MAW--ECM equilibrium equations}
\label{app:MAW_linearization}

This appendix derives the consistent Newton--Raphson linearization of
the MAW--ECM reduced equilibrium equations introduced in
Eq.~\refeq{eq:MAW_ROM_equilibrium}. The residual depends on the latent
coordinates through the constitutive response, the manifold tangent
basis,  and additionally
through the adaptive cubature weights. The resulting tangent operator
therefore comprises the standard and manifold-curvature contributions,
together with a new term arising from the variation of the weight
fields.

For convenience, Eq.~\refeq{eq:MAW_ROM_equilibrium} is first recast in
residual form. Defining
\begin{equation}
\label{eq:MAW_online_residual}
\Rdec(\qM;\muvec)
\defeq
\sum_{g\in\Zomega}
\rDec{g}(\qM;\muvec)\,
\omegag{g}(\qM)
-
\FextD(\qM;\muvec),
\end{equation}
where \(\rDec{g}\) and \(\FextD\) are given by
Eqs.~\refeq{eq:rDec_definition} and~\refeq{eq:FextD}, respectively,
the reduced equilibrium equations read
\begin{equation}
\label{eq:MAW_online_equilibrium_residual}
\Rdec(\qM;\muvec)=\bm{0}.
\end{equation}
The consistent tangent matrix,
\(\KtanD\in\RRn{\rD}{\rD}\), is obtained by differentiating
Eq.~\refeq{eq:MAW_online_residual} with respect to the latent
coordinates. Using the factorized decoder introduced in
Eq.~\refeq{eq:decoder_tau_generic}, together with
Eqs.~\refeq{eq:PhiD_tau_generic} and~\refeq{eq:FextD}, the result can
be written as
\begin{equation}
\label{eq:manifold_tangent_split_maw}
\KtanD
=
\KstdD
+
\KcurvD
+
\KweigD .
\end{equation}
The three terms account, respectively, for the constitutive
linearization, the curvature of the approximation manifold, and the
variation of the adaptive cubature weights.

Let
\[
\KFEg{g}
:=
\frac{\partial\rFE{g}}{\partial\dred}
\]
denote the consistent Gauss-point tangent operator. The standard contribution
is
\begin{equation}
\label{eq:manifold_standard_tangent_maw}
\KstdD
=
\JTAU^T
\left(
\sum_{g\in\Zomega}
\PhiROM^T
\KFEg{g}
\PhiROM\,
\omegag{g}
\right)
\JTAU.
\end{equation}
The curvature contribution is expressed componentwise as
\begin{equation}
\label{eq:manifold_curvature_tangent_maw}
\left(\KcurvD\right)_{ab}
=
\left[
\frac{\partial^2\TAU}
{\partial q_{M,a}\,\partial q_{M,b}}
\right]^T
\RTauD,
\qquad
a,b=1,\ldots,\rD,
\end{equation}
where
\begin{equation}
\label{eq:residual_tau_maw}
\RTauD
:=
\sum_{g\in\Zomega}
\PhiROM^T
\rFE{g}\,
\omegag{g}
-
\FextONEtau\muvec .
\end{equation}
Notice that the residual in Eq.~\eqref{eq:MAW_online_residual}
satisfies
\begin{equation}
\label{eq:residual_factorization_maw}
\Rdec
=
\JTAU^T\RTauD.
\end{equation}
Finally, the dependence of the cubature weights on the latent
coordinates gives
\begin{equation}
\label{eq:KweightD_def}
\begin{aligned}
\KweigD
&=
\JTAU^T
\sum_{g\in\Zomega}
\left(
\PhiROM^T\rFE{g}
\right)
\frac{\partial\omegag{g}}{\partial\qM}
\\
&=
\sum_{g\in\Zomega}
\rDec{g}
\frac{\partial\omegag{g}}{\partial\qM}.
\end{aligned}
\end{equation}
Here,
\(\partial\omegag{g}/\partial\qM\in\mathbb{R}^{1\times\rD}\)
is obtained by differentiating the regression model used to represent
the corresponding weight field.

For a fixed-weight rule,
\(\partial\omegag{g}/\partial\qM=\bm{0}\), and therefore
\(\KweigD=\bm{0}\). Likewise, \(\KcurvD\) vanishes when the decoder is
affine, since \(\JTAU\) is then constant and
\(\partial^2\TAU/\partial\qM^2=\bm{0}\). In the standard linear-ROM
case, \(\JTAU=\ident\).

\bibliographystyle{unsrt}
\bibliography{Bibliography_used}

@Book{Nocedal1999,
  author    = {J. Nocedal and S. J. Wright},
  title     = {Numerical optimization},
  publisher = {Springer},
  address   = {New York},
  year      =1999,
url={BIBLIOGRAPHY/Numerical_optimization_Nocedal__Springer_.pdf}}

@book{boyd2004convex,
  title={{Convex optimization}},
  author={Boyd, S.P. and Vandenberghe, L.},
  isbn={0521833787},
  year={2004},
  publisher={Cambridge Univ Pr},
 url={BIBLIOGRAPHY/Convex_Optimization_Boyd__Cambridge_.pdf}
}

@article{miehe2003computational,
  title={{Computational micro-to-macro transitions for discretized micro-structures of heterogeneous materials at finite strains based on the minimization of averaged incremental energy* 1}},
  author={Miehe, C.},
  journal={Computer Methods in Applied Mechanics and Engineering},
  volume={192},
  number={5-6},
  pages={559--591},
  issn={0045-7825},
  year={2003},
  publisher={Elsevier},
url={BIBLIOGRAPHY/miehe2003computational.pdf}
}

@article{an2009optimizing,
  title={{Optimizing cubature for efficient integration of subspace deformations}},
  author={An, S.S. and Kim, T. and James, D.L.},
  journal={ACM transactions on graphics},
  volume={27},
  number={5},
  pages={165},
  year={2009},
  publisher={NIH Public Access},
  url={BIBLIOGRAPHY/an2009optimizing.pdf}
}

@article{hernandez2017dimensional,
  title={Dimensional hyper-reduction of nonlinear finite element models via empirical cubature},
  author={Hern{\'a}ndez, J. A. and Caicedo, M. A. and Ferrer, A.},
  journal={Computer Methods in Applied Mechanics and Engineering},
  volume={313},
  pages={687--722},
  year={2017},
  publisher={Elsevier},
   url = {BIBLIOGRAPHY/hernandez2016dimensional.pdf}
}

@article{hernandez2014high,
  title={High-performance model reduction techniques in computational multiscale homogenization},
  author={Hern{\'a}ndez, J. A. and Oliver, J and Huespe, AE and Caicedo, MA and Cante, JC},
  journal={Computer Methods in Applied Mechanics and Engineering},
  volume={276},
  pages={149--189},
  year={2014},
  publisher={Elsevier},
   url = {BIBLIOGRAPHY/hernandez2014high.pdf}
}

@article{farhat2014dimensional,
  title={Dimensional reduction of nonlinear finite element dynamic models with finite rotations and energy-based mesh sampling and weighting for computational efficiency},
  author={Farhat, Charbel and Avery, Philip and Chapman, Todd and Cortial, Julien},
  journal={International Journal for Numerical Methods in Engineering},
  volume={98},
  number={9},
  pages={625--662},
  year={2014},
  publisher={Wiley Online Library},
   url = {BIBLIOGRAPHY/farhat2014dimensional.pdf} 
}

@article{farhat2015structure,
  title={Structure-preserving, stability, and accuracy properties of the energy-conserving sampling and weighting method for the hyper reduction of nonlinear finite element dynamic models},
  author={Farhat, Charbel and Chapman, Todd and Avery, Philip},
  journal={International Journal for Numerical Methods in Engineering},
  volume={102},
  number={5},
  pages={1077--1110},
  year={2015},
  publisher={Wiley Online Library},
url={BIBLIOGRAPHY/farhat2015structure.pdf}
}

@article{correa2015mechanical,
  title={Mechanical design of negative stiffness honeycomb materials},
  author={Correa, Dixon M and Seepersad, Carolyn Conner and Haberman, Michael R},
  journal={Integrating Materials and Manufacturing Innovation},
  volume={4},
  number={1},
  pages={1--11},
  year={2015},
  publisher={Springer},
  url={BIBLIOGRAPHY/correa2015mechanical.pdf}
}

@article{bremer2010nonlinear,
  title={A nonlinear optimization procedure for generalized Gaussian quadratures},
  author={Bremer, James and Gimbutas, Zydrunas and Rokhlin, Vladimir},
  journal={SIAM Journal on Scientific Computing},
  volume={32},
  number={4},
  pages={1761--1788},
  year={2010},
  publisher={SIAM},
   url={BIBLIOGRAPHY/bremer2010nonlinear.pdf}
}

@article{xiao2010numerical,
  title={A numerical algorithm for the construction of efficient quadrature rules in two and higher dimensions},
  author={Xiao, Hong and Gimbutas, Zydrunas},
  journal={Computers \& mathematics with applications},
  volume={59},
  number={2},
  pages={663--676},
  year={2010},
  publisher={Elsevier},
   url={BIBLIOGRAPHY/xiao2010numerical.pdf}
}

@article{qiu2004curved,
  title={A curved-beam bistable mechanism},
  author={Qiu, Jin and Lang, Jeffrey H and Slocum, Alexander H},
  journal={Journal of microelectromechanical systems},
  volume={13},
  number={2},
  pages={137--146},
  year={2004},
  publisher={IEEE},
   url = {BIBLIOGRAPHY/qiu2004curved.pdf}
}

@article{lee2019model,
  title={Model reduction of dynamical systems on nonlinear manifolds using deep convolutional autoencoders},
  author={Lee, Kookjin and Carlberg, Kevin T},
  journal={Journal of Computational Physics},
  pages={108973},
  year={2019},
  publisher={Elsevier},
   url = {BIBLIOGRAPHY/lee2019model.pdf}
}

@article{yano2019lp,
  title={An LP empirical quadrature procedure for reduced basis treatment of parametrized nonlinear PDEs},
  author={Yano, Masayuki and Patera, Anthony T},
  journal={Computer Methods in Applied Mechanics and Engineering},
  volume={344},
  pages={1104--1123},
  year={2019},
  publisher={Elsevier},
   url = {BIBLIOGRAPHY/yano2019lp.pdf}
}

@article{yano2019discontinuous,
  title={Discontinuous Galerkin reduced basis empirical quadrature procedure for model reduction of parametrized nonlinear conservation laws},
  author={Yano, Masayuki},
  journal={Advances in Computational Mathematics},
  pages={1--34},
  year={2019},
  publisher={Springer},
   url = {BIBLIOGRAPHY/yano2019discontinuous.pdf}
}

@article{hernandez2020multiscale,
  title={{A multiscale method for periodic structures using domain decomposition and ECM-hyperreduction}},
  author={Hern{\'a}ndez, J. A.},
  journal={Computer Methods in Applied Mechanics and Engineering},
  volume={368},
  pages={113192},
  year={2020},
  publisher={Elsevier},
   url = {BIBLIOGRAPHY/hernandez2020multiscale.pdf}
}

@article{pennec2006riemannian,
  title={A Riemannian framework for tensor computing},
  author={Pennec, Xavier and Fillard, Pierre and Ayache, Nicholas},
  journal={International Journal of computer vision},
  volume={66},
  number={1},
  pages={41--66},
  year={2006},
  publisher={Springer},
   url = {BIBLIOGRAPHY/pennec2006riemannian.pdf}
}

@article{fresca2022pod,
  title={POD-DL-ROM: enhancing deep learning-based reduced order models for nonlinear parametrized PDEs by proper orthogonal decomposition},
  author={Fresca, Stefania and Manzoni, Andrea},
  journal={Computer Methods in Applied Mechanics and Engineering},
  volume={388},
  pages={114181},
  year={2022},
  publisher={Elsevier},
   url = {BIBLIOGRAPHY/fresca2022pod.pdf}
}

@article{touze2021model,
  title={Model order reduction methods for geometrically nonlinear structures: a review of nonlinear techniques},
  author={Touz{\'e}, Cyril and Vizzaccaro, Alessandra and Thomas, Olivier},
  journal={Nonlinear Dynamics},
  volume={105},
  number={2},
  pages={1141--1190},
  year={2021},
  publisher={Springer},
   url = {BIBLIOGRAPHY/touze2021model.pdf}
}

@article{barnett2022quadratic,
  title={Quadratic approximation manifold for mitigating the Kolmogorov barrier in nonlinear projection-based model order reduction},
  author={Barnett, Joshua and Farhat, Charbel},
  journal={Journal of Computational Physics},
  volume={464},
  pages={111348},
  year={2022},
  publisher={Elsevier},
   url = {BIBLIOGRAPHY/barnett2022quadratic.pdf}
}

@article{ares2023hyper,
  title={{Hyper-reduction for Petrov--Galerkin reduced order models}},
  author={S. {Ares de Parga} and Bravo, J. R. and Hern{\'a}ndez, J. A. and Zorrilla, R and Rossi, R},
  journal={Computer Methods in Applied Mechanics and Engineering},
  volume={416},
  pages={116298},
  year={2023},
  publisher={Elsevier},
url = {BIBLIOGRAPHY/ares2023hyper.pdf}
}

@article{hernandez2024cecm,
  title={{CECM: A continuous empirical cubature method with application to the dimensional hyperreduction of parameterized finite element models}},
  author={Hernandez, J. A. and Bravo, J. R. and S. {Ares de Parga}},
  journal={Computer Methods in Applied Mechanics and Engineering},
  volume={418},
  pages={116552},
  year={2024},
  publisher={Elsevier},
url = {BIBLIOGRAPHY/hernandez2024cecm.pdf}
}

@article{barnett2023neural,
  title={Neural-network-augmented projection-based model order reduction for mitigating the Kolmogorov barrier to reducibility},
  author={Barnett, Joshua and Farhat, Charbel and Maday, Yvon},
  journal={Journal of Computational Physics},
  volume={492},
  pages={112420},
  year={2023},
  publisher={Elsevier},
url = {BIBLIOGRAPHY/barnett2023neural.pdf}
}

@article{bravo2024geometrically,
  title={Geometrically Parametrised Reduced Order Models for Studying the Hysteresis of the Coanda Effect in Finite Element-based Incompressible Fluid Dynamics},
  author={Bravo, JR and Stabile, G and Hess, M and Hernandez, J.A. and Rossi, R and Rozza, G},
  journal={Journal of Computational Physics},
  pages={113058},
  year={2024},
  publisher={Elsevier},
url = {BIBLIOGRAPHY/bravo2024geometrically.pdf}
}

@article{bravo2024subspace,
  title={A subspace-adaptive weights cubature method with application to the local hyperreduction of parameterized finite element models},
  author={Bravo, JR and Hern{\'a}ndez, JA and S. {Ares de Parga} and Rossi, Riccardo},
  journal={International Journal for Numerical Methods in Engineering},
  volume={125},
  number={24},
  pages={e7590},
  year={2024},
  publisher={Wiley Online Library},
url = {BIBLIOGRAPHY/bravo2024subspace.pdf}
}

@article{sperling2024comparative,
  title={A comparative study of enriched computational homogenization schemes applied to two-dimensional pattern-transforming elastomeric mechanical metamaterials},
  author={Sperling, SO and Guo, T and Peerlings, RHJ and Kouznetsova, VG and Geers, MGD and Roko{v{s}}, O},
  journal={Computational Mechanics},
  volume={74},
  number={1},
  pages={169--190},
  year={2024},
  publisher={Springer},
  url = {BIBLIOGRAPHY/sperling2024comparative.pdf}
}

@article{ares2025parallel,
title = {Parallel reduced-order modeling for digital twins using high-performance computing workflows},
journal = {Computers \& Structures},
volume = {316},
pages = {107867},
year = {2025},
issn = {0045-7949},
author = {S. {Ares de Parga} and J.R. Bravo and N. Sibuet and J.A. Hern\'andez and R. Rossi and Stefan Boschert and Enrique S. Quintana-Ortí and Andrés E. Tomás and Cristian Catalin Tatu and Fernando Vazquez-Novoa and Jorge Ejarque and Rosa M. Badia},
  url = {BIBLIOGRAPHY/ares2025parallel.pdf}
}

@article{jain2019hyper,
  title={Hyper-reduction over nonlinear manifolds for large nonlinear mechanical systems},
  author={Jain, Shobhit and Tiso, Paolo},
  journal={Journal of Computational and Nonlinear Dynamics},
  volume={14},
  number={8},
  pages={081008},
  year={2019},
  publisher={American Society of Mechanical Engineers},
  url       = {BIBLIOGRAPHY/jain2019hyper.pdf}
}

@article{ares2026nonlinear,
  title={Nonlinear projection-based model order reduction with machine learning regression for closure error modeling in the latent space},
  author={S. {Ares de Parga} and Tezaur, Radek and Hern{\'a}ndez, Carlos G and Farhat, Charbel},
  journal={Computer Methods in Applied Mechanics and Engineering},
  volume={448},
  pages={118443},
  year={2026},
  publisher={Elsevier},
  url       = {BIBLIOGRAPHY/ares2026nonlinear.pdf}
}

@article{chen2025negative,
  title={Negative stiffness mechanical metamaterials based on curved beams for reusable shock isolation},
  author={Chen, Shuai and Lian, Xu and Liu, Xin and Hu, Jiqiang and Wang, Bing and Xu, Jie and Wu, Linzhi},
  journal={International Journal of Smart and Nano Materials},
  volume={16},
  number={2},
  pages={397--418},
  year={2025},
  publisher={Taylor \& Francis},
  url       = {BIBLIOGRAPHY/chen2025negative.pdf}
}

@article{natarajan1995sparse,
  title={Sparse approximate solutions to linear systems},
  author={Natarajan, Balas Kausik},
  journal={SIAM journal on computing},
  volume={24},
  number={2},
  pages={227--234},
  year={1995},
  publisher={SIAM}
}

@article{belkin2003laplacian,
  title={Laplacian eigenmaps for dimensionality reduction and data representation},
  author={Belkin, Mikhail and Niyogi, Partha},
  journal={Neural computation},
  volume={15},
  number={6},
  pages={1373--1396},
  year={2003},
  publisher={MIT Press},
  url       = {BIBLIOGRAPHY/belkin2003laplacian.pdf}
}

@article{fresca2021comprehensive,
  author  = {Fresca, Stefania and Dede', Luca and Manzoni, Andrea},
  title   = {A Comprehensive Deep Learning-Based Approach to Reduced Order Modeling of Nonlinear Time-Dependent Parametric Partial Differential Equations},
  journal = {Journal of Scientific Computing},
  volume  = {87},
  number  = {20},
  year    = {2021},
  doi     = {10.1007/s10915-021-01462-7},
  url       = {BIBLIOGRAPHY/fresca2021comprehensive.pdf}
}

@book{spielman2019spectral,
  author    = {Daniel A. Spielman},
  title     = {Spectral and Algebraic Graph Theory},
  year       = {2019},
  publisher  = {Yale University},
  note       = {Lecture notes, available at http://cs-www.cs.yale.edu/homes/spielman/sagt/sagt.pdf}
}

@article{simo1987strain,
  title={Strain-and stress-based continuum damage models—I. Formulation},
  author={Simo, Juan C and Ju, Jiannwen W},
  journal={International journal of solids and structures},
  volume={23},
  number={7},
  pages={821--840},
  year={1987},
  publisher={Elsevier},
  url       = {BIBLIOGRAPHY/simo1987strain.pdf}
}

@book{rockafellar1996convex,
  author    = {Rockafellar, R. Tyrrell},
  title     = {Convex Analysis},
  publisher = {Princeton University Press},
  address   = {Princeton, NJ},
  year      = {1996},
  url       = {BIBLIOGRAPHY/Rockafellar_Convex_analysis_1996.pdf}
}

@article{koike2026sparse,
  title={Sparse POD Mode Selection and Manifold Dimensionality Reduction with Neural Networks},
  author={Koike, Tomoki and Mohan, Prakash and de Frahan, Marc T Henry and Qian, Elizabeth and Bessac, Julie},
  journal={arXiv preprint arXiv:2605.27756},
  year={2026},
  url       = {BIBLIOGRAPHY/koike2026sparse.pdf}
}

@article{zuo2026nonlinear,
  title={Nonlinear Correlation of POD and DMD Modal Coefficients in Reduced-Order Modeling of Flow Around a Cylinder in a Microchannel},
  author={Zuo, Bin and Yang, Xiaopei and Wang, Haichun and Xiao, Qianhao},
  journal={Micromachines},
  volume={17},
  number={7},
  pages={778},
  year={2026},
  publisher={MDPI},
  url       = {BIBLIOGRAPHY/zuo2026nonlinear.pdf}
}

@article{vizzaccaro2020comparison,
  title={Comparison of nonlinear mappings for reduced-order modelling of vibrating structures: normal form theory and quadratic manifold method with modal derivatives},
  author={Vizzaccaro, Alessandra and Salles, Loic and Touz{\'e}, Cyril},
  journal={Nonlinear Dynamics},
  year={2020},
  url       = {BIBLIOGRAPHY/vizzaccaro2020comparison.pdf}
}

@article{callaham2022role,
  title={On the role of nonlinear correlations in reduced-order modelling},
  author={Callaham, Jared L and Brunton, Steven L and Loiseau, Jean-Christophe},
  journal={Journal of Fluid Mechanics},
  volume={938},
  pages={A1},
  year={2022},
  url       = {BIBLIOGRAPHY/callaham2022role.pdf}
}

@book{wendland2005scattered,
  title={Scattered data approximation},
  author={Wendland, Holger},
  volume={17},
  year={2005},
  publisher={Cambridge university press Cambridge},
  url       = {BIBLIOGRAPHY/Wendland_Scattered_Data_Approximation_2004.pdf}
}

@article{poggio1990networks,
  title={Networks for approximation and learning},
  author={Poggio, Tomaso and Girosi, Federico},
  journal={Proceedings of the IEEE},
  volume={78},
  number={9},
  pages={1481--1497},
  year={1990},
  publisher={IEEE},
  url       = {BIBLIOGRAPHY/poggio1990networks.pdf}
}

@article{chmiel2025unified,
  title={Unified {LSPG} model reduction framework and assessment for hypersonic computational fluid dynamics},
  author={Chmiel, Matthew R and Barnett, Joshua and Farhat, Charbel},
  journal={AIAA Journal},
  volume={63},
  number={1},
  pages={72--90},
  year={2025},
  publisher={American Institute of Aeronautics and Astronautics},
}

@article{grimberg2021mesh,
  title={Mesh sampling and weighting for the hyperreduction of nonlinear {Petrov--Galerkin} reduced-order models with local reduced-order bases},
  author={Grimberg, Sebastian and Farhat, Charbel and Tezaur, Radek and Bou-Mosleh, Charbel},
  journal={International Journal for Numerical Methods in Engineering},
  volume={122},
  number={7},
  pages={1846--1874},
  year={2021},
  publisher={Wiley Online Library},
}

@article{zhang2026unified,
  author  = {Zhang, Qinghua and Ritzert, Stephan and Zhang, Jian and
             Kehls, Jannick and Reese, Stefanie and Brepols, Tim},
  title   = {A unified multi-perspective quadratic manifold for
             mitigating the {Kolmogorov} barrier in multiphysics damage},
  journal = {Journal of the Mechanics and Physics of Solids},
  volume  = {209},
  pages   = {106499},
  year    = {2026},
  doi     = {10.1016/j.jmps.2025.106499}
}

\end{document}